\documentclass[11pt, reqno]{amsart}  
\usepackage{amsmath, mathrsfs}
\usepackage{tikz}

\usepackage[all]{xy}
\SelectTips{cm}{}
\usepackage{amsthm, amsfonts}
\usepackage{geometry}
  \usepackage[overload]{empheq}
  
  \usepackage{tikz-cd}
  \usepackage{url}

  \usepackage{amssymb}
  
  \usepackage{graphicx}
  \usepackage{enumitem}
  \usepackage{color}
  \usepackage{enumitem}
  \usepackage[colorlinks=true]{hyperref}
  \hypersetup{urlcolor=cyan, citecolor=red, linkcolor=blue}
               
  \usepackage[capitalize]{cleveref}

  \usepackage{thmtools}
  \usepackage{float}

		\crefname{lemma}{Lemma}{Lemmas}
		\crefname{theorem}{Theorem}{Theorems}
		\crefname{prop}{Proposition}{Propositions}
		\crefname{cor}{Corollary}{Corollaries}
	
		\usepackage{xcolor}
\usepackage[most]{tcolorbox}
\usepackage{listings}

\definecolor{white}{rgb}{1,1,1}
\definecolor{mygreen}{rgb}{0,0.4,0}
\definecolor{light_gray}{rgb}{0.97,0.97,0.97}
\definecolor{mykey}{rgb}{0.117,0.403,0.713}

\tcbuselibrary{listings}
\newlength\inwd
\newtcblisting{pyin}[1][]{%
sharp corners,
enlarge left by=\inwd,
width=\linewidth-\inwd,
enhanced,
boxrule=0pt,
colback=light_gray,
listing only,
top=0pt,
bottom=0pt,
overlay={
  \node[
	anchor=north east,
	text width=\inwd,
	font=\footnotesize\ttfamily\color{mykey},
	inner ysep=1.5mm,
	inner xsep=0pt,
	outer sep=0pt
	] 
	at (frame.north west)
	{\refstepcounter{ipythcntr}\label{#1}In \theipythcntr:};
}
listing engine=listing,
listing options={
  aboveskip=1pt,
  belowskip=1pt,
  basicstyle=\footnotesize\ttfamily,
  language=Python,
  keywordstyle=\color{mykey},
  showstringspaces=false,
  stringstyle=\color{mygreen}
},
}

		\newtheoremstyle{1}
		{6pt} 
		{0pt} 
		{\itshape} 
		{} 
		{\bfseries} 
		{.} 
		{.5em} 
		{} 
		
		\newtheoremstyle{2}
		{6pt} 
		{0pt} 
		{} 
		{} 
		{\bfseries} 
		{.} 
		{.5em} 
		{} 
		
		\theoremstyle{1}
		\newtheorem{theorem}{Theorem}
		\newtheorem{lemma}[theorem]{Lemma}
		\newtheorem{prop}[theorem]{Proposition}
		\newtheorem{cor}[theorem]{Corollary}

    \theoremstyle{definition}
		\newtheorem{defn}[theorem]{Definition}
		\newtheorem{remark}[theorem]{Remark}
		\newtheorem{example}[theorem]{Example}

		\newtheorem{conjecture}[theorem]{Conjecture}

	\numberwithin{equation}{section}
		\numberwithin{lemma}{section}
		\numberwithin{theorem}{section}
    \numberwithin{prop}{section}
	\numberwithin{cor}{section}
		\numberwithin{defn}{section}
		\numberwithin{remark}{section}
		\numberwithin{example}{section}
		\numberwithin{assumption}{section}
    \numberwithin{conjecture}{section}
    \numberwithin{notation}{section}

        \newcommand{\F}{\mathbb{F}}
		\newcommand{\Q}{\mathbb{Q}}
		\newcommand{\Z}{\mathbb{Z}}
		
		\newcommand{\Sym}{\textnormal{Sym}}

        \newcommand{\Ker}{\textnormal{Ker}}
		\newcommand{\Gal}{\textnormal{Gal}}
		\newcommand{\Hom}{\textnormal{Hom}}

		\newcommand{\flag}{\mathcal{F}l}
		
		\newcommand{\KS}{\textnormal{KS}}
		\newcommand{\LL}{\mathcal{L}}
		
    \newcommand{\GL}{\textnormal{GL}}

    \newcommand{\Ver}{\textnormal{Ver}}
    \newcommand{\Rep}{\textnormal{Rep}}
    
    \newcommand{\crys}{\textnormal{crys}}
    \newcommand{\Std}{\textnormal{Std}}
    \newcommand{\BGG}{\textnormal{BGG}}
    \newcommand{\dR}{\textnormal{dR}}
    
    \newcommand{\et}{\textnormal{\'et}}
    \newcommand{\tor}{\textnormal{tor}}
    \newcommand{\can}{\textnormal{can}}
    \newcommand{\sub}{\textnormal{sub}}
    \newcommand{\Sh}{\textnormal{Sh}}
    \newcommand{\Shbar}{\overline{\textnormal{Sh}}}

	\newcommand{\Sbt}{\textnormal{Sbt}}
	\newcommand{\G}{\mathbb{G}}
	\newcommand{\SL}{\textnormal{SL}}
	\newcommand{\GSp}{\textnormal{GSp}}
	\newcommand{\Lie}{\textnormal{Lie}}
  
  \newcommand{\Fpbar}{\overline{\F}_p}
 
   \newcommand{\der}{\textnormal{der}}
   \renewcommand{\H}{\textnormal{H}}
   \newcommand{\m}{\mathfrak{m}}

\begin{document}

\title{A de Rham weight part of Serre's conjecture and generalized mod $p$ BGG decompositions}
\author{Martin Ortiz}

\maketitle

\begin{abstract}
  We propose a geometric approach to the weight part of Serre’s conjecture via the de Rham cohomology of special fibers of Shimura varieties.
  We conjecture that this formulation is equivalent to the one using Serre weights and the \'etale cohomology of Shimura varieties. We treat in detail the cases 
  of $G \in \{\GL_3,\GSp_4\}$, which are the simplest groups of semisimple rank at least $2$, where new phenomena occur. 
  First we prove this equivalence for generic weights and generic non-Eisenstein eigensystems for a compact $U(2,1)$ Shimura variety such that $G_{\Q_p}=\GL_3$. 
  Our approach involves proving the generic concentration in middle degree of mod $p$
  de Rham cohomology with coefficients. We achieve this by constructing generalized mod $p$ BGG resolutions, extending the work of Lan-Polo.
  After applying the results 
  from our companion paper \cite{Hodge-paper}, the problem reduces to the construction of 
  explicit BGG resolutions for $p$-restricted simple modules $L(\lambda)$ in a certain mod $p$ parabolic category $\mathcal{O}$.
  For $\GSp_4$, we also compute some explicit BGG decompositions, which help us 
  to upgrade the generic weak entailment from \cite{paper} to the expected generic entailment, as well as proving the equivalence above for one of the upper alcoves. 

  \end{abstract}

\maketitle  

\tableofcontents

\section{Introduction}
In this article we introduce a de Rham weight part of Serre's conjecture for Shimura varieties, and we conjecture that it is equivalent to the usual weight part of Serre's conjecture. 
We recall the statement of the latter conjecture in the case of Shimura varieties  
(although it can be formulated in much more generality). 
Let $(G,X)$ be a Shimura datum of Hodge type, $p$ a prime, and $K$ a level which is hyperspecial at $p$. Let
$\Sh/\mathcal{O}$ be Kisin's \cite{Kisin-integral-model} integral model associated to $(G,X,K)$. A Serre weight is an irreducible $\Fpbar$-valued representation of $G(\F_p)$. They can be described as follows. Fix $T \subseteq B \subseteq G$ a maximal torus and 
a Borel. The irreducible algebraic representations of $G_{\Fpbar}$ are labelled by dominant weights: $\{L(\lambda): \lambda \in X_{+}(T)\}$. Define the set of $p$-restricted weights as 
$$
X_1(T)=\{ \lambda \in X^*(T) : 0 \le \langle \lambda, \alpha^{\vee} \rangle<p \text{ for all } \alpha \in \Delta \},
$$
where $\Delta$ is the set of simple roots of $G$. For $\lambda \in X_{1}(T)$ let $F(\lambda):=L(\lambda)_{\mid G(\F_p)}$, then the set $\{F(\lambda): \lambda \in X_1(T)\}$ is precisely the set of Serre weights
\cite[Lem 9.2.4]{Gee-Herzig-Savitt}. 
Attached to each $F(\lambda)$ there is an \'etale local system $\underline{F}(\lambda)$ on $\Sh_{\overline{\Q}_p}$ by varying the level at $p$. Let $\mathbb{T}$ be the global spherical Hecke algebra away from $p$ and the bad places of 
$\Sh$, such that it acts on the \'etale cohomology of $\Sh_{\overline{\Q}_p}$. Let $\mathfrak{m} \subseteq \mathbb{T}$ be a maximal ideal, i.e. a mod $p$ system of Hecke eigenvalues. Define the
set of modular Serre weights as
$$
W(\mathfrak{m}):=\{ F(\lambda): \H^{\bullet}(\Sh_{\overline{\Q}_p},\underline{F}(\lambda))_{\mathfrak{m}} \neq 0 \},
$$
where $F(\lambda)$ runs along all Serre weights. In most cases, one can attach a mod $p$ Galois representation $\overline{r}_{\mathfrak{m}}$ to $\mathfrak{m}$, with $\overline{\rho}_{\mathfrak{m}}$ 
its local Galois representation(s) at $p$. In that case we will also use the more common notation 
$W(\overline{r}_{\mathfrak{m}})$ or $W(\overline{\rho}_{\mathfrak{m}})$. In this context, \cite{Gee-Herzig-Savitt} give a conjectural recipe
for $W(\overline{r}_\mathfrak{m})$ in terms of the Breuil--Mezard conjecture. This recipe is what's commonly known as the weight part 
of Serre's conjecture. 

On the other hand, one can consider coherent cohomology of Shimura varieties. Let $\Shbar=\Sh_{\Fpbar}$ and $\Shbar^\tor$ a choice of toroidal compactification. Let $P \subseteq G$
be a choice of Hodge parabolic, with Levi $M$. For $\lambda \in X_{M,+}(T)$ an $M$-dominant weight, one can define automorphic vector bundles $\omega(\lambda)$ over $\Shbar^\tor$. Then, for a certain class of Shimura varieties one can compare 
the eigensystems appearing in mod $p$ coherent cohomology and mod $p$ \'etale cohomology. For $\lambda \in X_{+}(T)$ let $V(\lambda) \in \Rep_{\check{\Z}_p}(G)$ 
be the dual Weyl module of $G$, then $L(\lambda)$ is the socle 
of $V(\lambda)_{\Fpbar}$. By restricting $V(\lambda)_{\Fpbar}$ to a $G(\F_p)$ representation it also defines a local system $\underline{V}(\lambda)_{\Fpbar}$ over $\Sh_{\overline{\Q}_p}$.
Let $\mathfrak{m}$ be generic and non-Eisenstein in the sense of \cite[Def 1.1.6]{Hodge-paper}. Then one expects
\begin{equation} \label{coherent-to-betti}
\H^0(\Shbar^\tor,\omega(\lambda+\eta))_{\mathfrak{m}} \neq 0 \iff \H^{d}_{\et}(\Sh_{\overline{\Q}_p},\underline{V}(\lambda)_{\Fpbar})_{\mathfrak{m}} \neq 0 
\end{equation}
to hold for $\lambda \in X_1(T)$ generic as in \cite[Def 1.1.1]{Hodge-paper}. Here $d$ is the dimension of $\Sh$ and  $\eta$ is the weight of $\Omega^d_{\Sh}$.
In particular, one can prove that \eqref{coherent-to-betti} holds for Harris--Taylor Shimura varieties or the Siegel threefold \cite[Prop 7.2.4]{Hodge-paper} \cite[Prop 5.1]{paper}. 

One motivating goal for us was to formulate an analogue of the weight part of Serre's conjecture 
in terms coherent cohomology. Crucially, when the rank of one of the factors of $G^{\der}$ is at least $2$, it is not always the case that $F(\lambda)=V(\lambda)_{\Fpbar}$.
Therefore, just considering 
$\H^0(\Shbar^\tor,\omega(\lambda+\eta))_{\mathfrak{m}}$ fails to distinguish which Jordan--Holder factors of $V(\lambda)_{\Fpbar}$ are modular. 
For Hilbert modular forms, precisely when the rank of all the factors of $G^{\der}$ is $1$, \cite{Diamond-Sasaki} formulate a weight part of Serre's conjecture for coherent cohomology in degree $0$.
 For $p$-restricted weights $\lambda$ their conjecture is conjectured to be equivalent to the weight part of Serre's conjecture in terms of 
Serre weights. This expectation is in line with \eqref{coherent-to-betti}, but we note that importantly, their conjecture includes  irregular weights for which \eqref{coherent-to-betti} does not apply, and weights 
that are not $p$-restricted.
For groups of higher rank no such conjecture exists.

Instead of considering coherent cohomology, we propose to 
look at de Rham cohomology, which provides a bridge between the two approaches above. 
It is sufficiently close to étale cohomology to capture the modularity of Serre weights, and we will construct generalized BGG filtrations 
that can be computed explicitly via coherent cohomology and theta linkage maps.
Given an algebraic representation $V \in \Rep_{\Fpbar}(G)$ there is an associated vector bundle with connection $\underline{V}$ on 
$\Shbar^{\tor}$ (the canonical extension), so we can define its de Rham cohomology $\H^{*}_{\dR}(\Shbar^\tor,\underline{V}):=\H^{*}(\underline{V} \otimes \Omega^{\bullet}_{\Shbar^\tor}(\log D))$. 
The spherical Hecke algebra $\mathbb{T}$ also acts on de Rham cohomology. 
Let $\mathfrak{m} \subseteq \mathbb{T}$ be a maximal ideal. We denote the set of de 
Rham Serre weights as 
$$
W_{\dR}(\mathfrak{m}):=\{ \lambda \in X_1(T) : 
\H^{\bullet}_{\dR}(\Shbar^{\tor},\underline{L}(\lambda))_{\mathfrak{m}} \neq 0 \}.
$$
For us, the de Rham weight part of Serre's conjecture would be the study of $W_{\dR}(\mathfrak{m})$. 
We conjecture that it is equivalent to the usual weight part of Serre's conjecture. 
Let $\epsilon \ge 0$. We define 
$W^{\epsilon}_{\dR}(\mathfrak{m}) \subseteq W_{\dR}(\mathfrak{m})$ as the 
subset consisting of the $\lambda$ which are $\epsilon$-generic as in \cite[Def 1.1.1]{Hodge-paper}. We define 
$W^{\epsilon}(\mathfrak{m}) \subseteq W(\mathfrak{m})$ analogously.

\begin{conjecture} \label{etale-dR-SW}
Let $(G,X)$ be a Shimura datum of Hodge type,
and let $\mathfrak{m} \subseteq \mathbb{T}$ be a maximal ideal.  
\begin{enumerate}
  \item (Strong form) Then 
$$
W_{\dR}(\mathfrak{m})=W(\mathfrak{m}).
$$
We use the surjection from $\lambda \in X_1(T)$ to Serre weights $F(\lambda)$ to make sense of this equality. 
\item (Weak form) Assume that $\Sh$ is a unitary Shimura variety or the Siegel threefold, and that 
$\mathfrak{m}$ is generic non-Eisenstein. 
Then there exists some $\epsilon \ge 0$ only depending on $G_{\Q}$ such that 
$$
W^{\epsilon}_{\dR}(\mathfrak{m})=W^{\epsilon}(\mathfrak{m}).
$$

\end{enumerate}
\end{conjecture}

Importantly, for generic non-Eisenstein $\mathfrak{m}$ there have been results of \cite{caraiani-scholze-compact}
\cite{caraiani-scholze-non-compact} \cite{Hamann-Lee} \cite{generic-vanishing-abelian} 
and others proving the concentration in middle degree
of \'etale cohomology for 
most abelian type Shimura varieties. When this generic 
concentration of mod $p$ \'etale cohomology holds, 
\Cref{etale-dR-SW}(2)
follows
from the following conjecture about concentration of de Rham cohomology, see \Cref{concentration-implies-dR-et}.
\begin{conjecture} \label{conjecture-concentration}
Let $\Shbar$ be the Siegel threefold or a unitary Shimura variety,
and let $\mathfrak{m}$ be generic and non-Eisenstein.
Then 
$$
\H^{\bullet}_{\dR}(\Shbar^\tor,\underline{L}(\lambda))_{\mathfrak{m}}$$
is concentrated in middle degree for generic $\lambda \in X^{*}(T)$ such that there exists $\mu \in X_1(T)$
with $L(\lambda) \in \textnormal{JH}[V(\mu)_{\Fpbar}]$. We use the term generic $\lambda$ in the sense of \cite[Def 1.1.1]{Hodge-paper}.
\end{conjecture}

In the cases when $G_{\Fpbar} \in \{\GL_3,\GSp_4\}$ the set of $\lambda$ considered above is contained in 
$X_1(T)$, but in 
general it is larger.
Using the formalism of prismatic cohomology (\Cref{prismatic-argument}) one can prove that for any $i$ the dimension of
$\H^{i}_{\dR}(\Shbar,\underline{L}(\lambda))_{\mathfrak{m}}$  for $\lambda \in X_{1}(T)$ is at least the one of  
$\H^{i}_{\et}(\Sh_{\overline{\Q}_p},F(\lambda))_{\mathfrak{m}}$.
Thus, for compact Shimura varieties 
we have that $W(\m) \subseteq W_{\dR}(\m)$. As evidence for \Cref{conjecture-concentration}, 
we prove it in the case that $G_{\Q_p}=\GL_3 \times \mathbb{G}_m$, which is the simplest unknown case. In the case 
that $G^{\der}_{\overline{\Q}}$ is almost-simple
\Cref{conjecture-concentration} is known for PEL Shimura varieties and $\lambda$ generic in the lowest
alcove by \cite{Lan-Suh-non-compact}.
\begin{theorem}[\Cref{concentration-GL3}] \label{concentration-GL3-intro}
Let $\Shbar$ be a compact unitary Shimura variety of signature $(2,1)$ with $p$ split in the quadratic imaginary field, 
so that $G_{\Q_p}=\GL_3 \times \mathbb{G}_m$.  
 Then \Cref{conjecture-concentration} holds in this case. If 
\Cref{concentration-coherent}(1) also holds for $p$ inert, then \Cref{conjecture-concentration} is also true in that case. 
\end{theorem}

We prove this by constructing a generalized BGG decomposition for weights in the upper alcove of $\GL_3$. For that we use
the de Rham realization functor from \cite{Hodge-paper}
\begin{equation} \label{realization-functor}
\Psi: D^{b}(\mathcal{O}_{P,\Fpbar}) \xrightarrow{f} D^{b}((\Shbar/\Fpbar)_{\crys}) \xrightarrow{g} D^{b}(C_{\Shbar/\Fpbar}).
\end{equation}
Here $\mathcal{O}_{P,\Fpbar}$ is a category of $(U\mathfrak{g},P)_{\Fpbar}$-modules modelled after the classical complex category $\mathcal{O}$, 
$(\Shbar/\Fpbar)_{\crys}$ is the crystalline topos over $\Fpbar$, and $C_{\Shbar/\Fpbar}$ is the category 
of $\mathcal{O}_{\Shbar}$-modules with $\Fpbar$-linear maps. All the maps are exact, and the image of $f$ lands in the subcategory 
of complexes of crystals on 
$(\Shbar/\Fpbar)_{\crys}$.
 The map $g$ is the pushforward 
from the crystalline topos 
to the Zariski topos. Furthermore, $\Psi$ appropriately extends to toroidal compactifications. 
For $W \in \Rep(P)$ let $\Ver_{P}(W)=U\mathfrak{g} \otimes_{U\mathfrak{p}} W$ be a parabolic Verma module.
The functor satisfies that $\Psi(L(\lambda)[0])=\underline{L}(\lambda)^{\vee} \otimes \Omega^{\bullet}_{\Shbar}$,
 and that $\Psi(\Ver_{P}(W)[0])=\mathcal{W}^{\vee}[0]$ where $\mathcal{W}:=P_{\dR} \times^{P} W$ is the
 vector bundle associated to $W$. Moreover,
  maps between parabolic Verma modules are sent by $\Psi$ to 
differential operators, and these are directly related to the theta linkage maps constructed in \cite{Hodge-paper}. Thus, to understand 
de Rham cohomology in terms of coherent cohomology of automorphic vector bundles
 we are reduced to studying generalized BGG
resolutions of $L(\lambda)$ in $D^{b}(\mathcal{O}_{P,\Fpbar})$ using Verma modules. This is the main content of this article. 
We obtain the following. For $\lambda \in X_{M,+}(T)$ let $L_{M}(\lambda) \subseteq W(\lambda)$ be the irreducible $M$-representation of highest weight $\lambda$ 
and the dual Weyl module for $M$, respectively. 
\begin{prop}[\Cref{GL3-BGG}] \label{intro-BGG}
Let $G=\GL_3/\Fpbar$, $P$ the $(2,1)$ block parabolic and $\lambda_1 \in C_1 \subseteq X_1(T)$ lying in the upper alcove. Let $\{\lambda_i,\mu_i,\nu_i : i=0,1\}$ be some of the affine 
Weyl reflections of $\lambda_1$, as drawn in \Cref{figure3}.
Then there exists a descending $3$ step filtration $F^{\bullet}$ of $L(\lambda_1)$ in $D^{b}(\mathcal{O}_{P,\Fpbar})$ such that 
\begin{enumerate}
\item $F^2=[0 \to \Ver_P W(\lambda_0) \xrightarrow{\phi} \Ver_P W(\lambda_1)]$
with $\phi$ a non-zero map.
\item $\textnormal{gr}^1=[ 0 \to \Ver_{P} L_{M}(\mu_1) \to 0]$. 
\item $\textnormal{gr}^{0}=[\Ver_{P} W(\nu_1) \xrightarrow{\psi} \Ver_{P} W(\nu_0) \to 0]$ with $\psi$ a non-zero map.
\end{enumerate}
In \Cref{BGG-lambda1-gsp4} and \Cref{BGG-lambda2-gsp4} we construct similar $4$ step filtrations on $L(\lambda)$ for $G=\GSp_4/\Fpbar$, $P$ the Siegel parabolic, and $\lambda$ inside one of the bottom $3$ (in the $\uparrow$ order)
$p$-restricted alcoves.
\end{prop}

Then applying \eqref{realization-functor} to the resolutions above, we get a filtration $F^{\bullet}_{\Shbar^{\tor}}L(\lambda)$ of 
$\underline{L}(\lambda)  \otimes \Omega^{\bullet}_{\Shbar^\tor}(\log D)$ whose graded pieces are
explicit complexes of automorphic vector bundles, 
and whose differentials are related to theta linkage maps as defined in \cite{Hodge-paper}. 
Here we use automorphic vector bundles to refer to the vector bundles associated to the irreducible representations of $M$.
Let us explain the cohomological consequences of this for $\GL_3$, we will assume that $\Sh$ is compact. 
 \begin{figure}
  \centering
  \begin{tikzpicture} [scale=0.8]
    \node [color=black] (1) at (0,0) (snlabel) {$\bullet$};
    \node[color=black] (4) at (-0.5,-0.1) (snlabel) {$-\rho$};
    \draw[dotted] (0,0)--(5,0);
    \draw[dotted] (0,0)--(2.5,3.95);
    \draw[dotted] (0,0)--(-2.5,-3.95);
    \draw[dotted] (0,0)--(2.5,-3.95);
    \draw (0,0) -- (4,0);
    \draw (4,0) -- (2,3.16);
    \draw (0,0) --(2,3.16);
    \draw (2,3.16) --(6,3.16);
    \draw (4,0)--(6,-3.16);
    \draw (2,3.16)--(6,-3.16);
    \draw (2,-3.16)--(0,-6.32);
    \draw (2,-3.16)--(6,-3.16);
    \node at (3.94,-1.094) {$\bullet$};
    \node at (3.94,-1.39) {$\mu_1$};
    \draw (-2,-3.16)--(0,-6.32);
    \draw (4,0) -- (6,3.16);
    \draw (0,0) --(2,-3.16);
    \draw (4,0) --(2,-3.16);
    \draw (0,0) --(-2,-3.16);
    \draw (-2,-3.16)--(2,-3.16);
    \node (lambda0) at (2.6,0.75) {$\lambda_0$};
    \node at (3,0.5) {$\bullet$};
    \node at (3,-0.5) {$\bullet$};
    \node (mu0) at (2.7,-0.3) {$\mu_0$};
    \node at (3.94,1.094) {$\bullet$};
    \node (lambda1) at (3.6,1.3) {$\lambda_1$};
    \node at (-0.476,-2.7) {$\bullet$};
    \node at (-0.476,-3.62) {$\bullet$};
    \node at (-0.25,-3.9) {$\nu_1$};
    \node at (-0.76,-2.5) {$\nu_0$};
    \draw[dashed,thick, red,->] (0,0) -- (0,2);
    \node at (-0.3,2.3) {$\textcolor{red}{\alpha_2}$};
    \draw[dashed,thick, blue,->] (0,0) -- (1.58,-1);
    \node at (1.85,-1.1) {$\textcolor{blue}{\alpha_1}$};
    \end{tikzpicture}
  \caption{}
  \label{figure3}
 \end{figure}
 Then $\H^2_{\dR}(\Shbar,\underline{L}(\lambda_1))_{\mathfrak{m}}$ has a $3$ step filtration 
 $0 \subset F^2 \subset F^1 \subset \H^2_{\dR}(\Shbar,\underline{L}(\lambda_1))_{\mathfrak{m}}$. Let 
 $\{\tilde{\lambda}_i,\tilde{\mu}_i,\tilde{\nu}_i\}$ be the same affine Weyl reflections as in \Cref{figure3}, now corresponding to $\tilde{\lambda}_1:=
 -w_0 \lambda_1$. 
 Then for $\lambda_0 \in C_0$ generic, 
 \begin{itemize}
 \item $F^2=\text{coker}[\H^0(\Shbar,\omega(\lambda_0+\eta)) \xrightarrow{\theta_{\lambda_0 \uparrow \lambda_1}}.
 \H^0(\Shbar,\omega(\lambda_1+\eta))]_{\mathfrak{m}}$.
 \item $\text{gr}^1=\H^1(\Shbar,\underline{L}_M(\tilde{\mu}_1)^{\vee})_{\mathfrak{m}}$
 \item $\text{gr}^2=\Ker[\H^2(\Shbar,\omega(\tilde{\lambda}_1)^{\vee}) \xrightarrow{\theta^{*}_{\tilde{\lambda}_0 \uparrow \tilde{\lambda}_1}}
 \H^2(\Shbar,\omega(\tilde{\lambda}_0)^{\vee})]_{\mathfrak{m}}$ is Serre dual to $F^2 \H^2_{\dR}(\Shbar,\underline{L}(\lambda_1)^{\vee})_{\mathfrak{m}}$
 after twisting the Hecke action. 
 \end{itemize}
 Here $\eta=(1,1,-2)$, and $\theta_{\lambda_0 \uparrow \lambda_1}$ is a theta linkage map, which is injective on $\H^0$.
 In contrast with the lowest alcove BGG
 decomposition of \cite{Lan-Polo}, the graded pieces are not a direct sum of coherent cohomology groups,
 instead there 
 are non-trivial interactions via theta linkage maps. Moreover, $F^{\bullet}_{\Shbar}L(\lambda_1)$ is not the filtration induced by the Hodge filtration of 
 $\underline{L}(\lambda_1)$, but it contains more information. E.g. $\omega(\lambda_1+\eta)[-2]$ and $\omega(\lambda_0+\eta)[-1]$ are (a posteriori) two consecutive graded pieces of the Hodge filtration, and $F^2_{\Shbar}L(\lambda_1)$ contains 
 the information that the connecting map between them is $\theta_{\lambda_0 \uparrow \lambda_1}$.
 Combining the concentration of de Rham cohomology in middle degree 
 with this explicit BGG decomposition shows that for $\mathfrak{m}$ generic non-Eisenstein and $\lambda_0 \in C_0$ 
 generic we have that
 $$
\text{coker}[\H^0(\Shbar,\omega(\lambda_0+\eta)) \xrightarrow{\theta_{\lambda_0 \uparrow \lambda_1}}
 \H^0(\Shbar,\omega(\lambda_1+\eta))]_{\mathfrak{m}}  \neq 0 \implies
 F(\lambda_1) \in W(\mathfrak{m}).
 $$
 This refines the identity \eqref{coherent-to-betti} to capture the notion of being modular for a non-lowest-alcove
 Serre weight in terms of coherent cohomology.

 The concentration in middle degree of de Rham cohomology is related to the study of kernels and cokernels of theta linkage maps, 
 and to vanishing results for higher coherent cohomology. 
 For the former, in general we only know how to control kernels of \textit{simple} theta 
 linkage maps on $\H^0$ (and its Serre dual statement), 
 which happen to be the ones appearing in $F^{\bullet}_{\Shbar}L(\lambda_1)$ for $\GL_3$. For the latter, we use
 the results of \cite{Deding-unitary} which prove the vanishing of $\H^i(\Shbar,\omega(\lambda))$ for $i \ge 1$ and generic $p$-restricted $\lambda$, and the results of \cite{general-cone-conjecture0} 
 showing the vanishing of $\H^0(\Shbar,\omega(\lambda))$ for generic non-dominant weights. 
 For $\GL_3$ these two are enough for us, since $\Shbar$ is a surface. 
 Conversely, \Cref{conjecture-concentration}(2) for $\GL_3$ implies the known concentration of 
 $\H^{\bullet}(\Shbar,\underline{L}_M(\tilde{\mu}_1)^{\vee})$ in degree $1$.
 In general, almost nothing is known about
 the vanishing of mod $p$ higher 
 coherent cohomology of automorphic vector bundles, and one might hope that through \Cref{conjecture-concentration}
 one could get some information.
 E.g. see \Cref{remark-lambda1-gsp4}.

 \begin{remark}
 Even though in \Cref{concentration-GL3-intro} we only prove concentration in middle degree for generic $\lambda \in X_1(T)$, it might 
 be plausible that it holds for all $p$-restricted $\lambda$ if we impose some additional genericity condition
 on $\overline{\rho}_{\mathfrak{m}}$, the local Galois representation at $p$. Namely, that it lies in a certain open subset of the Emerton--Gee stack.
 See \Cref{non-generic-concentration} for an example. On the other hand,
 concentration in middle degree fails completely if $\lambda$ is allowed to be arbitrarily large,
 already in the case of the modular curve, see \Cref{non-concentration-modular-curve}.
 The same example shows that the naive extension of \Cref{etale-dR-SW} to arbitrary dominant 
 weights by considering \'etale cohomology with coefficients in $L(\lambda)_{\mid G(\F_p)}$ does not hold. 
 \end{remark}

 While proving \Cref{conjecture-concentration} for $\GSp_4$ remains out of reach with current methods, the explicit nature of the BGG decompositions is already powerful enough to obtain new results.
  In fact,
 we can upgrade the generic weak entailment for $\GSp_4$ of \cite[Thm 5.2]{paper}
  to the entailment predicted by the explicit description of Breuil--Mezard cycles in \cite{Le-Hung-Lin}. 
  The weak entailment states that for $\m$ non-Eisenstein and $\lambda_0 \in C_0$ generic 
  $F(\lambda_0) \in W(\m)$ implies $F(\lambda_1) \in W(\m)$ or $F(\lambda_2) \in W(\m)$, where $\lambda_{i}$ are two $p$-restricted affine Weyl reflections of $\lambda_0$.
  
 \begin{theorem}[\Cref{correct-entailment}]
 Let $\Shbar$ be the Siegel threefold. Let $\lambda_0 \in C_0$,
  and $\mathfrak{m}$ a non-Eisenstein 
 maximal ideal.  
   Then 
 $$
 F(\lambda_0) \in W(\m) \implies F(\lambda_2) \in W(\m)
 $$
 holds for generic $\lambda_0 \in C_0$. Moreover,
 \Cref{etale-dR-SW}(2) holds for such $\m$ and generic weights in the alcove $C_2$. 
 \end{theorem}

 This result illustrates why the mod $p$ de Rham framework is useful. On one hand, the BGG 
 decomposition underlying this entailment is entirely invisible at the level of étale cohomology.
On the other hand, coherent cohomology alone fails to distinguish between the modularity of $F(\lambda_1)$ and $F(\lambda_2)$. Moreover, compared to similar results like \cite{le2025serreweightconjecturesmathrmgsp4} our method 
does not rely on the Taylor--Wiles method, so no extra assumptions are needed.

 \subsection{Methods}
 Let $\lambda \in X_1(T)$.
 The main body of the article concerns the explicit computation of certain BGG resolutions of $L(\lambda)$ in $D^{b}(\mathcal{O}_{P,\Fpbar})$ for 
 $G \in \{\GL_3, \GSp_4\}$, in \Cref{intro-BGG}(1).
 To do so, we follow the same procedure to construct the BGG resolution in characteristic $0$.
We define $\BGG_{L(\lambda)}:=\Ver_{P}(\wedge^{\bullet} \mathfrak{g}/\mathfrak{p} \otimes L(\lambda))_{\chi}$ 
as the localization of the standard complex (which corresponds to the de Rham complex on $G/P$) by 
the mod $p$ character $\chi$ of the Harish--Chandra center of $U\mathfrak{g}$ corresponding to $\lambda$.
This is a resolution of $L(\lambda)$, 
and for $p$ greater than the Coxeter number of any factor of $G^{\text{der}}$, 
each term has a (non-canonical) descending filtration whose graded pieces are of the form $\Ver_{P} L_{M}(w \cdot \lambda)$
for $w \in W_{\text{aff}}$. In our examples we choose these filtrations so that the weights 
which are the largest in the $\uparrow$ order appear with smaller grading.
That way we can define a filtration of $\BGG_{L(\lambda)}$ whose graded pieces 
lie in a single Weyl chamber, since maps of Verma modules are compatible 
with the $\le$ order. 
The heuristic is that then the graded pieces of the resulting complex on $\Shbar$ would have terms 
whose coherent cohomology is concentrated in a single degree, as it is the case 
in characteristic $0$.
The construction of this filtration relies on 
the simple geometry of the $p$-restricted alcoves in our low rank examples, specifically on the interaction of 
the linkage order $\uparrow$ with the $\le$ order.  

The hardest part is proving that the graded pieces of this filtration are quasi-isomorphic to a simpler, 
explicit complex. As defined in the abelian category these graded pieces contain Verma modules with very large 
multiplicities, so we have to be clever about it. Ultimately our method is ad-hoc and relies 
on the knowledge of the formal characters
of $L(\lambda)$ for $p$-restricted $\lambda$, thanks to the BGG resolution in characteristic $0$ and the decomposition 
of dual Weyl modules into irreducibles. This suggests that in general one might be able to describe $\BGG_{L(\lambda)}$
using roughly as many Verma modules as the irreducible constituents appearing in $V(\lambda)_{\Fpbar}$.
Finally, to get 
the cohomological results we plug these generalized BGG resolutions into the de Rham realization functor, 
and then we use all the available results from \cite{Hodge-paper}. Namely, we 
use all the vanishing results of coherent cohomology available, as well as results 
on the injectivity of theta linkage maps. We also use the compatibility of the filtration on $\BGG_{L(\lambda)}$ with 
Serre duality on $\Shbar^\tor$, so that in the $\GL_3$ case we only 
need to prove concentration in middle degree of $2$ out of the $3$ graded pieces of $\BGG_{\Shbar}L(\lambda_1)$.
In order to establish \Cref{correct-entailment} we need to prove some new results of independent interest. In \Cref{random-vanishing} and \Cref{vanishing1} we prove the vanishing of a coherent cohomology group in degree $2$ 
for a non-dominant weight, which falls right outside the scope of all previous vanishing results. 
For that we carefully use the Ekedahl--Oort stratification of the flag Shimura variety, along with the weight elimination 
results of \cite{heejong-lee}. Finally, we coin the notion of a $p$-translated lowest alcove BGG complex in \Cref{p-translation}, and together with a result in the vein of \cite{Deligne-Illusie} \cite{Illusie-coefficients} 
we use it to prove that a certain theta linkage map is injective in cohomological degree $1$ in \Cref{injectivity3}. This new idea is responsible for the second version of the article, as before 
proving such injectivity result seemed out of reach. 

\subsection{Notation}
Let $G/\Fpbar$ be a reductive group. Fix $T \subset B \subset P \subset G$
a maximal torus, Borel subgroup and a parabolic subgroup.
Let $P=MU_{P}$ be the Levi decomposition of 
$P$. Let $B^{-}$ be the opposite Borel, uniquely determined by $B \cap B^{-}=T$.
\begin{itemize}
\item Let $X^{*}(T)=\Hom_{\Fpbar}(T,\G_m)$ and $X_{*}(T)=\Hom_{\Fpbar}(\G_m,T)$ be the group of characters and cocharacters, respectively.
\item 
Let $\Phi^{+} \subset \Phi$ be the set of positive roots determined by $B$, and let
$\Delta \subset \Phi^{+}$ be the set of simple roots. 
\item Let $\lambda,\mu \in X^*(T)$. We say that $\lambda \le \mu$ if $\mu-\lambda \in \Z^{\ge 0}\Phi^{+}$.
\item Let $W$ be the Weyl group. For $\alpha \in \Phi$ let $s_{\alpha} \in W$ be the corresponding reflection. It acts on $X^*(T)$ as 
$\lambda \mapsto \lambda-\langle \lambda, \alpha^{\vee} \rangle \alpha$. Let $w_0 \in W$ be the longest element 
with respect to the Bruhat order. 
\item Let
 $W^{M}$ be the 
set of minimal length representatives of $W_M \backslash W$ in $W$.
\item (Dot action) Let $\rho \in X^*(T)$ be the half sum of positive roots for $G_{\Fpbar}$. We define the 
 dot action of $W$ on $X^*(T)$ by $w \cdot \lambda:=w(\lambda+\rho)-\rho$.
 \item (Affine Weyl group) Let $W_{\text{aff}}=p\mathbb{\Z}\Phi \rtimes W$ be the affine Weyl group.
 The dot action of $W$ on $X^*(T)$ extends to an action of $W_{\text{aff}}$ where $p\mathbb{\Z}\Phi$ acts by 
 translation. Then $W_{\text{aff}}$ is generated by reflections $s_{\gamma,n}$ for $\gamma \in \Phi$ and $n \in \Z$,
 which are characterized by acting as 
 $s_{\gamma,n} \cdot \lambda=\lambda+(pn-\langle \lambda+ \rho,\gamma^{\vee} \rangle) \gamma$.
 \item (Linkage order) We say that $\lambda \uparrow_{\gamma} \mu$ if there exists $n \in \Z$ such that $\mu=s_{\gamma,n} \cdot \lambda$,
$\lambda \le \mu$, and $n$ is the smallest integer such that $\lambda \le s_{\gamma,n} \cdot \lambda$. 
We say that $\lambda \uparrow \mu$ if there exists a chain 
$\lambda \uparrow_{\gamma_1} \lambda_1 \uparrow_{\gamma_2} \ldots \uparrow_{\gamma_n} \mu$. 
 \item We divide $X^{*}(T)-\rho$ into $\rho$-shifted alcoves, defined as the 
 interior of regions defined by the hyperplanes $H_{n,\gamma}=\{\langle \lambda+\rho, \gamma^{\vee} \rangle=np\}$
 for $\gamma \in \Phi$ and $n \in \Z$. The lowest alcove $C_0$ is defined by 
 $\{\lambda: 0 < \langle \lambda+\rho,\gamma^{\vee} \rangle <p \text{ for all } \gamma \in \Phi^{+}\}$.
 \item  We say that 
 $\lambda \in X^*(T)$ is $p$-small if $|\langle \lambda+ \rho, \gamma^{\vee} \rangle| <p$ for 
 all $\gamma \in \Phi$.
 The affine Weyl group acts simply transitively on $\rho$-shifted alcoves 
  and the stabilizer of the region of $p$-small weights is the finite Weyl group $W$.
\item Let $X_{+}(T) \subset X^{*}(T)$ be the set of dominant weights $\lambda$, satisfying 
$\langle \lambda, \alpha^{\vee} \rangle \ge 0$ for every $\alpha \in \Delta$. Similarly, let $X_{M,+}(T)$
be the set of $M$-dominant weights, consisting of the $\lambda$ satisfying 
$\langle \lambda, \alpha^{\vee} \rangle \ge 0$ for all $\alpha \in \Delta_{M}$. 
\item For a linear algebraic group $H$ over $\Fpbar$ let $\Rep_{\Fpbar}(H)$ be the category of algebraic representations 
of $H$ 
which are finite dimensional $\Fpbar$-vector spaces. 
\item (Dual Weyl modules) For $\lambda \in X_{+}(T)$ let 
$
V(\lambda):=\text{Ind}^{G}_{B} w_{0} \lambda=
\H^0(G/B,B \times w_{0}\lambda/B),
$
where $B$ acts as $b(g,v)=(gb,w_{0}\lambda(b)v)$. It is an algebraic representation of $G$ defined over $W(\Fpbar)$. We say that  
$V(\lambda)$ is the \textit{dual Weyl module} for $G$ of highest weight $\lambda$.
Analogously, for $\lambda \in X_{M,+}(T)$ we define 
$W(\lambda)$ as the dual Weyl module for $M$, with respect to the Borel $B \cap M$. We also consider it as a representation of $P$ by inflation. 
\item Let $L(\lambda) \subset V(\lambda)_{\Fpbar}$ be the socle as an $G$-representation, 
and $L_{M}(\lambda) \subset W(\lambda)_{\Fpbar}$ the socle as a $M$-representation.
The
$L(\lambda)$ for $\lambda \in X_{+}(T)$ parametrize all irreducible representations of $G_{\Fpbar}$ \cite[Cor II.2.7]{Janzten-book}.
\item Assume that $G$ descends to a reductive group $G/\F_p$. The set of $p$-restricted weights $X_1(T) \subseteq X_{+}(T)$ consists of those $\lambda$ such that 
$0 \le \langle \lambda, \alpha^{\vee} \rangle<p$ for all $\alpha \in \Delta$. Given $\lambda \in X_1(T)$
let $F(\lambda)\in \Rep_{\Fpbar}G(\F_p)$ be the restriction of $L(\lambda)$ to $G(\F_p)$. 
Then $\{F(\lambda): \lambda \in X_1(T)\}$ contains all the \textit{Serre weights}, i.e. irreducible representations in 
$\Rep_{\Fpbar}G(\F_p)$. Moreover, $F(\lambda)=F(\lambda')$ only if $\lambda-\lambda' \in 
(\sigma-\text{id})X^0(T):=(\sigma-\text{id})\{ \lambda \in X^*(T): 
\langle \lambda, \alpha^{\vee} \rangle=0 \text{ for all } \alpha \in \Delta\}$, where $\sigma$ is the Frobenius acting on $X^{*}(T)$
\cite[Lem 9.2.4]{Gee-Herzig-Savitt}. 
 
\end{itemize} 

\subsection{Acknowledgments}
The contents of this article are part of the author's PhD thesis. For his guidance and 
conversations during this period I thank my PhD advisor George Boxer. This work
was supported by the Engineering and Physical Sciences Research Council [EP/S021590/1]. 
The EPSRC Centre for Doctoral Training in Geometry and Number Theory 
(The London School of Geometry and Number Theory), University College London, King's College London and 
Imperial College London. Part of its writing was done while the author was a member of the Max Planck Institute for Mathematics 
in Bonn, which I thank for their hospitality and funding.

\section{The general conjecture}
We explain in more detail the relation between the main two conjectures in the introduction. 
For simplicity, we work with the Shimura varieties appearing in \Cref{etale-dR-SW}(2).
For $\Sh$ a unitary Shimura variety we say that a mod $p$ eigensystem $\mathfrak{m} \subseteq \mathbb{T}$
is non-Eisenstein if $\overline{r}_{\mathfrak{m}}$ is irreducible. For $G=\GSp_4$ we say that $\mathfrak{m}$
is non-Eisenstein if the composition of $\overline{r}_{\mathfrak{m}}: \Gal_{\Q} \to \GSp_{4}(\Fpbar)$ with 
the standard representation $\GSp_4 \hookrightarrow \GL_4$ is irreducible. 
We will use the notion of genericity of $\mathfrak{m}$
as in \cite[Def 1.1.6]{Hodge-paper}, originating in \cite{Hamann-Lee}. We recall that there are functors $F_{Q}: \Rep_{Q} \to \text{Coh}(\Sh^{\tor})$ for 
$Q \in \{G,P\}$, and for $V \in \Rep(Q)$ we let $\mathcal{V}:=F_{Q}(V)$.

\begin{prop} \label{concentration-implies-dR-et}
Let $\mathfrak{m}$ be generic non-Eisenstein. Assume that $\H^{\bullet}_{\et}(\Sh_{\overline{\Q}_p},\underline{V})_{\mathfrak{m}}$ is concentrated in middle degree for any 
$V \in \Rep_{\Fpbar}G(\F_p)$. Then
\begin{enumerate}
\item \Cref{conjecture-concentration} implies \Cref{etale-dR-SW}(2). 
\item (Integral version) Assume that $G^{\der}_{\overline{\Q}}$ is almost-simple and
 that both $\H^{\bullet}_{\dR}(\Sh^\tor,\mathcal{V}(\lambda))_{\mathfrak{m}}/\check{\Z}_p$
and $\H^{\bullet}_{\dR}(\Sh^\tor,\mathcal{V}(\lambda)^{\vee})_{\mathfrak{m}}$ are concentrated
 in middle degree and torsion free 
 for all generic $\lambda \in X_1(T)$. Then \Cref{conjecture-concentration} holds. Conversely, \Cref{conjecture-concentration} implies the above concentration and torsion-freeness in middle degree. 
\end{enumerate}
\begin{proof}
For (1) we 
proceed by induction on the $\le$ order on the weights.
 By standard homological algebra considerations,
and the rational \'etale-crystalline comparison theorem we know that the Euler characteristic of 
$\H^{\bullet}_{\dR}(\Shbar^{\tor},\mathcal{V}(\lambda))$ and 
$\H^{\bullet}_{\text{\'et}}(\Sh_{\overline{\Q}_p},V(\lambda)_{\overline{\F}_p})$ are the same. 
Then for $\lambda \in C_0$, \Cref{etale-dR-SW}(2) follows from the observation above, the concentration in middle degree of 
both cohomologies, and the fact that $L(\lambda)=V(\lambda)_{\Fpbar}$. 
For general $\lambda \in X_1(T)$,
 we have an exact sequence 
$L(\lambda) \to V(\lambda) \to V$, where the highest weights of $V$ are lower than $\lambda$
in the $\uparrow$ order (hence also in the $\le$ order) by the linkage principle \cite[Prop II.6.13]{Janzten-book}.
Therefore, by the induction hypothesis 
we get concentration of $\H^{\bullet}_{\dR}(\Shbar^{\tor},\mathcal{V}(\lambda))_{\mathfrak{m}}$
in middle degree.
Since both are concentrated in middle degree, this implies that the dimensions of 
$\H^{\bullet}_{\dR}(\Shbar^{\tor},\mathcal{V}(\lambda))_{\mathfrak{m}}$ and 
$\H^{\bullet}_{\text{\'et}}(\Sh_{\overline{\Q}_p},V(\lambda)_{\overline{\F}_p})$
are the same.
Moreover, the same equality of dimensions holds 
with coefficients $L(\lambda)$ and $F(\lambda)=L(\lambda)_{\mid G(\F_p)}$ by the induction hypothesis and 
the long exact sequence corresponding to $L(\lambda) \to V(\lambda) \to V$. 

For (2) by assumption we get that $\H^{\bullet}_{\dR}(\Shbar^\tor,\mathcal{V}(\lambda))_{\mathfrak{m}}$ is concentrated in middle degree. 
By considering the embedding $L(\lambda) \hookrightarrow V(\lambda)_{\Fpbar}$
 and induction on the weight we also get that 
$\H^{\bullet}_{\dR}(\Shbar^\tor,L(\lambda))_{\mathfrak{m}}$ is concentrated in degrees $[0,d]$ for
generic $\lambda \in X_1(T)$. The case of $\lambda \in C_0$ is handled 
by the integral BGG decomposition of \cite[Thm 8.2.3]{Hodge-paper} and \cite[Prop 8.12]{Lan-Suh-non-compact}. For the latter we use that 
$G^{\der}_{\overline{\Q}}$ is almost-simple so that it applies to all $\lambda_0 \in C_0$ generic. By doing the same procedure with
$\H^{\bullet}_{\dR}(\Shbar^\tor,\mathcal{V}(\lambda)^{\vee})_{\mathfrak{m}}$
we get concentration of $\H^{\bullet}_{\dR}(\Shbar^\tor,L(\lambda))_{\mathfrak{m}}$ in degree $d$. 
The converse direction follows easily. 
\end{proof}
\end{prop}

Importantly, the generic concentration of \'etale cohomology is now known in many cases, due to 
\cite{caraiani-scholze-compact} \cite{caraiani-scholze-non-compact} \cite{Hamann-Lee} \cite{generic-vanishing-abelian}
among others. \footnote{One needs to combine their results 
with the fact that for the appropriate notion of non-Eisensteinness, usual and compactly supported \'etale cohomology localized at $\mathfrak{m}$ agree.}
We note that in (2) we expect that the concentration of one should imply the concentration of the other by Poincare duality.
 One would need to adapt 
 the usual Poincare duality for smooth proper varieties 
 \cite[\href{https://stacks.math.columbia.edu/tag/0FW3}{Tag 0FW3}]{stacks-project} 
 to the case of coefficients and the relative log setting, and then prove that usual and compactly supported 
 de Rham cohomology are the same after localizing at $\mathfrak{m}$.

\begin{remark} \label{prismatic-argument}
We \footnote{The author thanks Bao Le Hung and Daniel Le for pointing out this argument to us.} 
sketch an argument that shows that for compact Shimura varieties 
$$
W(\mathfrak{m}) \subseteq W_{\dR}(\mathfrak{m}).
$$ 
Let $\mathfrak{X}$ be the $p$-adic formal completion of $\Sh$ over $W=W(\Fpbar)$. Then 
\cite{prismatic-realization} construct a realization functor 
$\Rep_{W}(G) \to \text{Vect}(\mathfrak{X}^{\text{syn}})$ where $\mathfrak{X}^{\text{syn}}$ is the 
syntomification of $\mathfrak{X}$ \cite{F-gauges}. 
The relevance of these is that there are de Rham and \'etale specializations \cite[\S 5.3.13, 6.3.2]{F-gauges}
$$
T_{\dR}: \text{Perf}(\mathfrak{X}^{\text{syn}}) \to D_{\dR}(\mathfrak{X})\;\;\;
T_{\et} : \text{Perf}(\mathfrak{X}^{\text{syn}}) \to D_{\text{lisse}}(\Sh_{\overline{\Q}_p},W)
$$
where $D_{\dR}(\mathfrak{X})$ is the derived category of vector bundles with a flat connection and a Griffiths-transversal 
filtration. They prove that these specializations when applied to $V \in \Rep_{W}(G)$
 yield $(\mathcal{V},\nabla)$ with its Hodge filtration, 
 and the \'etale local system $\underline{V}$ respectively. 
 Both specialization maps come from pullback to a special locus on $\mathfrak{X}^{\text{syn}}$, at least locally, 
 and the de Rham specialization is "deeper". That is, locally on $\text{Spf}(R) \to \mathfrak{X}^{\text{syn}}$
 there is an ideal $I \subseteq R$
 such that $T_{\et}$ corresponds to the base change to $R[1/I]$, and $T_{\dR}$ corresponds 
 to the base change to $R/I$ after twisting by Frobenius.
The key point is that $L(\lambda)$ is a quotient of the Weyl module $V(\lambda)^{\vee}$ over 
$W$ by an element of $\Rep_{W}(G)$, so it is naturally a perfect complex in $\Rep_{W}(G)$. Thus, we get an element 
$M \in \text{Perf}(\mathfrak{X}^{\text{syn}})$, which is the cone of a map between objects coming from $\Rep_{W}(G)$. 
Since both specializations commute with colimits/cones, we see that $T_{\dR}(M)=(L(\lambda),\nabla)$ and 
$T_{\et}(M)=\underline{L}(\lambda)$. Therefore, $M$ computes both the desired \'etale and de Rham cohomology. From the fact 
that the de Rham specialization is deeper we get by semicontinuity that 
$$
\text{dim} \H^i_{\dR}(\Shbar,L(\lambda)) \ge \text{dim} \H^i_{\et}(\Sh_{\overline{\Q}_p},\underline{L}(\lambda)).
$$
One can also check that both specialization functors are compatible with Hecke operators away from $p$, 
so we also obtain the result after localizing at $\mathfrak{m}$. 
\end{remark}

Note that the above remark would imply that we don't need to assume concentration in middle degree of \'etale cohomology for 
\Cref{concentration-implies-dR-et}(1) to hold. 

\begin{example}(The modular curve) \label{non-concentration-modular-curve}
The strong form of \Cref{etale-dR-SW} is true in the case of the modular curve. We use weights of $\SL_2$
for simplicity. For $L(k-2)$ and $2 \le k \le p$ this follows from the integral \'etale-de Rham comparison theorem. 
This comparison theorem holds since $\H^{*}_{\dR}(\Sh^\tor,L(k-2))$ can be realized as a subspace of $\H^{*}_{\dR}(E^{k-2}/\Z_p)$ as in \cite[\S 5]{Lan-Suh-1},
where the product is taken over $\Sh^{\tor}$, and the dimension of $E^{k-2}$ is less than $p$. For $k=p+1$ the integral BGG decomposition still holds, and 
$\H^{*}_{\dR}(\Shbar^\tor,L(p-1))$ is computed by a complex of the form $[\omega^{1-p}\to \omega^{p+1}]$ in degree $0,1$. 
We have that $\omega^{1-p}$
and $\omega^{p+1}$ have cohomology concentrated in degree $1$ and $0$ respectively, 
so that de Rham cohomology is concentrated in degree $1$. Then the log analogue of \Cref{prismatic-argument} implies that \'etale cohomology must also be concentrated in degree $1$, so we conclude by 
\Cref{concentration-implies-dR-et}(1).

The example of the modular curve also serves to illustrate that when the weight $\lambda$ is too big
 de Rham cohomology will not be
concentrated in middle degree, even after localizing at a generic global parameter.  Let 
$2 \le k \le p-1$ and consider the vector bundle with connection $\mathcal{H}^{(p^2)} \otimes \Sym^{k-2} \mathcal{H}$ associated to $L(p^2+k-2)=L(1)^{(p^2)} \otimes L(k-2)$.
 There is an exact sequence as 
$(U\mathfrak{g},P)$-modules $0 \to L_1 \to L(p^2+k-2) \to L_2 \to 0$, where $L_1=p^2 \otimes L(k-2)$
and $L_2=-p^2 \otimes L(k-2)$. By "$p$-translating" the integral BGG complex for $L(k-2)$
we get a filtration on $L(p^2+k-2) \otimes \Omega^{\bullet}_{\Shbar^\tor}$ with subobject 
$C^{\bullet}_0=[\omega^{p^2+2-k} \xrightarrow{\tilde{\theta}} \omega^{p^2+k}]$ and quotient 
$C^{\bullet}_1=[\omega^{-p^2+2-k} \to \omega^{-p^2+k}]$, where $\tilde{\theta}$ is the theta linkage map.
Then $\Ker(\tilde{\theta}: \H^0(\Shbar^\tor,\omega^{p^2+2-k}) \to \H^0(\Shbar^\tor,\omega^{p^2+k}))$ is contained 
in $\H^0_{\dR}(\Shbar^\tor,L(p^2+k-2))$. Let $H$ be the Hasse invariant. By the relations $\theta^p=H^{p+1}\theta$ and $\tilde{\theta}=\theta^{k-1}/H^{k-1}$ the former is equal to the kernel of $\theta$,
which one can compute (\cite{katz-theta}) to be given by the embedding 
$$
\H^0(\Shbar^\tor,\omega^{p+2-k}) \to \H^0(\Shbar^\tor,\omega^{p(p+2-k)}) \xrightarrow{H^{k-2}} 
\H^0(\Shbar^\tor,\omega^{p^2+2-k}),
$$
where the first map is the $p$th power map.
Thus, we see that 
$\H^0(\Shbar^\tor,\omega^{p+2-k})_{\m} \subseteq \H^0_{\dR}(\Shbar^\tor,L(p^2+k-2))_{\mathfrak{m}}$, which can be non-zero even after localizing at a generic global eigensystem. A similar analysis shows that 
$\H^{\bullet}_{\dR}(\Shbar^\tor,L(\lambda))_{\mathfrak{m}}$ is always concentrated for $\lambda \in X_1(T)$ and $\mathfrak{m}$ non-Eisenstein. 

This example also shows that \Cref{etale-dR-SW}(2) cannot hold for arbitrary dominant weights.
That is,
$\H^{\bullet}_{\dR}(\Shbar^\tor,L(\lambda))_{\mathfrak{m}} \neq 0$ does not in general imply that
$\H^{\bullet}_{\et}(\Sh_{\overline{\Q}_p},L(\lambda)_{\mid G(\F_p)})_{\mathfrak{m}} \neq 0$ for an arbitrary $\lambda \in X_{+}(T)$.
For $\lambda=p^2+k-2$ 
we have $L(\lambda)_{\mid G(\F_p)}=F(1)\otimes F(k-2)$ by Steinberg's tensor product, which has $F(k-1)$ and $F(k-3)$
as irreducible constituents. We can choose $\overline{\rho}_{\mathfrak{m}}$ to be reducible, non-split,
and such that
 $\H^0(\Shbar^\tor,\omega^{p+2-k})_{\mathfrak{m}} \neq 0$, so that $\H^0_{\dR}(\Shbar^\tor,L(p^2+k-2))_{\mathfrak{m}} \neq 0$.
  However, neither $F(k-1)$ nor $F(k-3)$ are in $W(\m)=\{F(p-k)\}$ for a generic $2 \le k \le p-1$. 
\end{example}
For an example where de Rham cohomology of a non-generic $p$-restricted weight (or heuristically, when 
$\overline{\rho}_{\mathfrak{m}}$ lies in some special locus of the Emerton--Gee stack)
is not concentrated in middle degree see \Cref{non-generic-concentration}.

\section{Generalized mod $p$ BGG decompositions}
Let $G$ be a reductive group over $\Fpbar$ and $P \subseteq G$ a parabolic subgroup. We also fix a maximal torus $T$
and a Borel $B$ such that $T \subseteq B \subseteq P$. Our goal is to construct certain generalized mod $p$ BGG resolutions of 
simple $G$-representations $L(\lambda)$ in $D^{b}(\mathcal{O}_{P,\Fpbar})$.
We start by recalling our definition of category $\mathcal{O}_{P,\Fpbar}$.

\begin{defn}(Verma modules and a mod $p$ category $\mathcal{O}$) \label{verma-defn}
A $(U \mathfrak{g},P)$-module over $\Fpbar$ is a $\Fpbar$-module $M$ equipped with an action of $U\mathfrak{g}$
and an algebraic action of $P$, i.e. $M$ is a filtered union of elements in $\Rep_{\Fpbar}(P)$.
Moreover, $M$ 
 satisfies the following additional properties.
\begin{enumerate}
\item $M$ is a finitely generated $U\mathfrak{g}$-module. 
\item Since the action of $P$ is algebraic, the derivative of $P$ on $M$ is well-defined.
It agrees with the restriction of the
$U\mathfrak{g}$ action to $U\mathfrak{p}$.
\end{enumerate}
We denote this category of $(U \mathfrak{g},P)$-modules by $\mathcal{O}_{P,\Fpbar}$.
Let $V \in \Rep_{\Fpbar}(P)$. We define $\Ver_{P}(V)\coloneqq U \mathfrak{g} \otimes_{U \mathfrak{p}} V
\in \mathcal{O}_{P,\Fpbar}$, with 
$U \mathfrak{g}$ acting by left multiplication on $U\mathfrak{g}$, and $P$ by the adjoint action on $U\mathfrak{g}$, and by its action on $V$ on the right. 
\end{defn}

Our definition of $\mathcal{O}_{P,\Fpbar}$ mirrors the definition of category $\mathcal{O}$ over 
the complex numbers, with two main changes. One is that we work with a general parabolic instead of a Borel. 
The second is that when working in characteristic $p$, one cannot just work with representations of Lie algebras, 
since the Lie algebra does not distinguish between weights
which which differ by multiples of $p$. This category also appears in \cite{modular-O} and 
\cite{Quan-O} for $P=B$, denoted as modular category $\mathcal{O}$.

We will use the notion of "Serre-duality" $V \mapsto V^*$ on $D^{b}(\mathcal{O}_{P,\Fpbar})$
from \cite[Def 6.3.2]{Hodge-paper}. It is characterized by $(\Ver_{P}W)^{*}:=\Ver_{P}W^*$, where $W \in \Rep_{\Fpbar}(P)$ and 
$W^*=W^{\vee} \otimes W(2\rho-2\rho_{M})$; and by the property that dual maps of Verma modules produce dual differential operators on $G/P$.
We recall a lemma that we will use to compute Serre duals of certain complexes. 
We endow $D^{b}(\mathcal{O}_{P,\Fpbar})$ with the natural $t$-structure, which enables to see $\mathcal{O}_{P,\Fpbar}$ as the subcategory of objects whose cohomology is concentrated in degree $0$. 
 
\begin{lemma} \cite[Lem 6.3.3]{Hodge-paper} \label{lemma-serre-dual-acyclic}
Let $V \in \mathcal{O}_{P,\Fpbar}$ have a finite filtration by Verma modules. Then $[V]^*$ lies in $\mathcal{O}_{P,\Fpbar}$, and we simply denote it by $V^*$. Similarly, 
if $C=[V_1 \to V_2 \to \ldots V_n] \in D^{b}(\mathcal{O}_{P,\Fpbar})$ is a complex with each $V_i$ having a finite filtration by Verma modules, then 
$$
C^*=[V^{*}_n \to V^{*}_{n-1}\to \ldots \to V^*_1].
$$
\end{lemma}

We will use the following lemmas often in the next subsection.
\begin{lemma} \label{splitting-lemma}
Let $\lambda \in X_{+}(T)$. Then 
\begin{enumerate}
\item (The linkage principle) $\textnormal{Ext}^1_{G}(L(\lambda),L(\mu)) \neq 0$ implies that
 $\mu \in W_{\textnormal{aff}} \cdot \lambda$. Moreover, if $L(\mu) \in \textnormal{JH}[V(\lambda)_{\Fpbar}]$ then $\mu \uparrow \lambda$.
\item $\textnormal{Ext}^1_{G}(L(\lambda),L(\lambda))=0$.
\item Let $\lambda \in X_{M,+}(T)$. Then $\textnormal{Ext}^1_{P}(L_{M}(\lambda),L_{M}(\lambda))=0$.
\item $\textnormal{Ext}^1_{(U\mathfrak{g},P)}(\Ver_P L_{M}(\lambda),\Ver_P L_{M}(\lambda))=0$.
\end{enumerate}
\begin{proof}
Part $(1)$ follows directly from \cite[Cor II.6.17, Prop II.6.13]{Janzten-book}.
Part $(2)$ is \cite[II.2.12]{Janzten-book}. For part $(3)$ such an extension $V$ will be split as a $M$-module 
by part $(2)$, 
generated by vectors $v_1, v_2$ of weight $\lambda$, with $v_1$ corresponding to the sub $P$-representation. 
Let $g \in U_P$, then since $U_{P}$ acts trivially on the quotient 
$gv_2 \in v_2+M\langle v_1 \rangle$, but since $\Lie(U_P)$ only contains positive root spaces, we conclude that
$gv_2=v_2$. Therefore, $U_P$ acts trivially on both of the generators. Since $P$ normalizes $U_P$, we conclude that 
$U_P$ acts trivially on the whole $V$, hence $V$ splits as a $P$-representation. 
For $(4)$ let $W$ be such an extension. It contains two vectors of weight $\lambda$, one $v_1$
corresponding to the subobject, and 
$v_2$ corresponding to the quotient. Let $W'$ be the $P$-submodule of $W$ generated by $v_1$ and $v_2$. Then 
there is an exact sequence $0 \to L_{M}(\lambda) \to W' \to L_{M}(\lambda)$. If the last map is zero, it means that 
$v_2$ is in the span of $Pv_1$, but this is not possible by considering their weights. 
Therefore, $W'$ is a self 
extension of $L_{M}(\lambda)$, hence it is split by part $(2)$.
 That splitting defines a map $\Ver_P L_{M}(\lambda) \to W$
mapping $v_{\lambda}$ to $v_2$,
by the universal property of Verma modules. 
Since $\Ver_P L_M(\lambda)$ is generated by $v_{\lambda}$ as a $(U\mathfrak{g},P)$-module
it means the map is a section. 
\end{proof}
\end{lemma}

\begin{lemma} \label{splitting-lemma-2}
Let $\lambda, \mu \in X^{*}(T)$. 
\begin{enumerate}
\item Assume that $W(\mu)$ is irreducible, then 
$$
\textnormal{Ext}^1_{P}(L_{M}(\lambda),W(\mu)) \neq 0
$$
implies that $\lambda > \mu$.
\item Assume that $W(\mu)$ is irreducible, then 
$$
\textnormal{Ext}^1_{\mathcal{O}_{P,\Fpbar}}(\Ver_{P}L_{M}(\lambda),\Ver_P W(\mu)) \neq 0
$$
implies that $\lambda > \mu$.
\item Suppose that $\lambda \notin W_{M,\textnormal{aff}} \cdot \mu$. Then 
$$
\textnormal{Ext}^1_{\mathcal{O}_{P,\Fpbar}}(\Ver_{P}L_{M}(\lambda),\Ver_P L_{M}(\mu)) \neq 0
$$
implies that $\lambda > \mu$.
\end{enumerate}
\begin{proof}
For $(1)$ let $W \in \textnormal{Ext}^1_{P}(L_{M}(\lambda),W(\mu))$, and assume that $\lambda \ngeq \mu$.
 If $\mu=\lambda$ the claim is \Cref{splitting-lemma}(3). Otherwise, let $v_{\mu} \in W$ be
the unique vector of weight $\mu$. Since $W(\mu)$ is irreducible we have 
$W(\mu)=W(-w_{0,M}\mu)^{\vee}$, then by \cite[II.2.13(b)]{Janzten-book} we have 
$\Hom_{P}(W(\mu),W)=\Hom_{P}(W^{\vee},W(-w_{0,M}\mu))=\Hom_{B}(W^{\vee},-\mu)$. Since $v_{\mu}$ is a highest weight vector, we get a map $W(\mu) \to W$
which sends the highest weight vector of $W(\mu)$ to $v_{\mu}$. Since $W(\mu)$ is irreducible 
it is injective, so it defines a section of 
$W \to W(\mu)$. Part $(2)$ follows from $(1)$ and \Cref{splitting-lemma}(4). That is, let $W' \in \textnormal{Ext}^1_{\mathcal{O}_{P,\Fpbar}}(\Ver_{P}L_{M}(\lambda),\Ver_P W(\mu))$ such that 
$\lambda \ngeq \mu$. Then the $P$-subrepresentation of $W'$ generated by
$L_{M}(\lambda)$ and the unique vector of weight $\mu$ is an element of $\textnormal{Ext}^1_{P}(L_{M}(\lambda),W(\mu))$, which provides a splitting.
For $(3)$, let $W' \in \textnormal{Ext}^1_{\mathcal{O}_{P,\Fpbar}}(\Ver_{P}L_{M}(\lambda),\Ver_P L_{M}(\mu))$ and suppose that $\lambda \ngeq \mu$. 
Let $W \in \textnormal{Ext}^1_{P}(L_{M}(\lambda),L_{M}(\mu))$ be the sub $P$-representation of $W'$ generated by $L_{M}(\lambda)$ and the unique vector of weight $\mu$. Then $W$ is split as an $M$-representation 
by the linkage principle for $M$, and $U_{P}$ must act trivially on the vector of weight $\mu$. Thus, $W$ and $W'$ are split. 
\end{proof}
\end{lemma}

In the next subsections we generalize the lowest alcove BGG decomposition of  \cite{Tilouine-Polo} \cite{Lan-Polo}
 to some higher alcoves in the case of 
$\GL_3$ and $\GSp_4$. We make our definition of generic weights very explicit for our two particular examples. 
\begin{defn}
Let $n \ge 1$, $G \in \{ \GL_n, \GSp_{2n}\}$ and
let $X=\{X_{p}\}_{p}$ be a family of statements concerning $X^*(T)$ for each rational prime $p$. We say that $X$ holds for generic $\lambda \in X^*(T)$ if there exists some constant $\epsilon \ge 0$ depending only on $n$ 
such that each $X_p$ holds for all $\lambda$ that are $\epsilon$-generic. Here $\lambda$ is $\epsilon$-generic
if for every $\gamma \in \Phi$ there exists an integer $N_{\gamma}$ such that 
$N_{\gamma}p+\epsilon <\langle \lambda, \gamma^{\vee} \rangle< (N_{\gamma}+1)p-\epsilon$.
\end{defn}

\subsection{The case of $\GL_3$}
Let $G=\GL_3$ in this subsection, with $T$ the diagonal maximal torus and $B$ the upper triangular Borel. A weight 
$\lambda$ will be given by a tuple $(k_1,k_2,k_3): \text{Diag}(t_1,t_2,t_3) \mapsto t^{k_1}_1 t^{k_2}_2 t^{k_3}_3$.
 Let $\alpha_1=(1,-1,0)$,
 $\alpha_{2}=(0,1,-1)$ be the simple roots, and $\rho=(2,1,0)$.
 Let $P$ be the standard parabolic of signature $(2,1)$, so that 
$\alpha_1 \in \Phi_{M}$. There are two $p$-restricted alcoves: the lowest alcove $C_0$, and the upper alcove $C_1$.
 Let $\lambda_0=(a,b,c) \in C_0$ and let $\lambda_1=s_{\alpha_1+\alpha_2,1} \cdot \lambda_0 \in C_1$ be 
 its affine reflection.
  Let $\mu_0=s_{\alpha_2,0} \cdot \lambda_0$, $\mu_1=s_{\alpha_1,1}
 \cdot \mu_0$, $\nu_0=s_{\alpha_1+\alpha_2,0} \cdot \mu_0$, $\nu_1=s_{\alpha_2,-1} \cdot \nu_0$ be some of the affine 
 reflections, as 
 shown in \Cref{figure5}. 
 \begin{figure}
  \centering
  \begin{tikzpicture} [scale=0.8]
    \node [color=black] (1) at (0,0) (snlabel) {$\bullet$};
    \node[color=black] (4) at (-0.5,-0.1) (snlabel) {$-\rho$};
    \draw[dotted] (0,0)--(5,0);
    \draw[dotted] (0,0)--(2.5,3.95);
    \draw[dotted] (0,0)--(-2.5,-3.95);
    \draw[dotted] (0,0)--(2.5,-3.95);
    \draw (0,0) -- (4,0);
    \draw (4,0) -- (2,3.16);
    \draw (0,0) --(2,3.16);
    \draw (2,3.16) --(6,3.16);
    \draw (4,0)--(6,-3.16);
    \draw (2,3.16)--(6,-3.16);
    \draw (2,-3.16)--(0,-6.32);
    \draw (2,-3.16)--(6,-3.16);
    \node at (3.94,-1.094) {$\bullet$};
    \node at (3.94,-1.39) {$\mu_1$};
    \draw (-2,-3.16)--(0,-6.32);
    \draw (4,0) -- (6,3.16);
    \draw (0,0) --(2,-3.16);
    \draw (4,0) --(2,-3.16);
    \draw (0,0) --(-2,-3.16);
    \draw (-2,-3.16)--(2,-3.16);
    \node (lambda0) at (2.6,0.75) {$\lambda_0$};
    \node at (3,0.5) {$\bullet$};
    \node at (3,-0.5) {$\bullet$};
    \node (mu0) at (2.7,-0.3) {$\mu_0$};
    \node at (3.94,1.094) {$\bullet$};
    \node (lambda1) at (3.6,1.3) {$\lambda_1$};
    \node at (-0.476,-2.7) {$\bullet$};
    \node at (-0.476,-3.62) {$\bullet$};
    \node at (-0.25,-3.9) {$\nu_1$};
    \node at (-0.76,-2.5) {$\nu_0$};
    \draw[dashed,thick, red,->] (0,0) -- (0,2);
    \node at (-0.3,2.3) {$\textcolor{red}{\alpha_2}$};
    \draw[dashed,thick, blue,->] (0,0) -- (1.58,-1);
    \node at (1.85,-1.1) {$\textcolor{blue}{\alpha_1}$};
    \end{tikzpicture}
  \caption{}
  \label{figure5}
  \end{figure}
We state the formal character of $L(\lambda_1)$, which follows from the fact that 
$0 \to L(\lambda_1) \to V(\lambda_1)_{\Fpbar} \to L(\lambda_0) \to 0$ is exact and the BGG decomposition in characteristic $0$. 
\begin{lemma} \label{character-GL3}
Let $\lambda_1 \in C_1$. We have an exact sequence of $P$-representations 
$$
0 \to L_{M}(\mu_1) \to W(\mu_1) \to W(\mu_0)\to 0
$$
and
the formal character of $L(\lambda_1)$ in $\mathcal{O}_{P,\Fpbar}$ is given by
$$
[L(\lambda_1)]=\Ver_{P} W(\lambda_1)-\Ver_{P} W(\lambda_0)-\Ver_{P} L_{M}(\mu_1)+\Ver_{P}W(\nu_1)-
\Ver_{P} W(\nu_0).
$$
\begin{proof}
We use the BGG resolution in characteristic $0$: for $i=0,1$ there is a resolution 
$$
0\to \Ver_{P} W(\nu_i) \to \Ver_{P} W(\mu_i) \to \Ver_{P} W(\lambda_i) \to V(\lambda_i)_{\overline{\Q}_p} \to 0.
$$
\end{proof}
\end{lemma}
Note that all the other dual Weyl modules for $M$ appearing in the character are irreducible. 
Using this result we can construct an explicit BGG-like decomposition for weights in the upper alcove $C_1$.
\begin{prop} \label{GL3-BGG}
Let $\lambda_1 \in C_1$. There exists a resolution $\BGG_{L(\lambda_1)}$ of $L(\lambda_1)$ 
in the category $\mathcal{O}_{P,\Fpbar}$
with a $3$-step filtration $0 \subset F^{2} \subset F^{1} \subset F^{0}=\BGG_{L(\lambda_1)}$ in $D^{b}(\mathcal{O}_{P,\Fpbar})$ such that 
\begin{enumerate}
\item $$
F^2=[0 \to \Ver_{P} W(\lambda_0) \xrightarrow{\phi} \Ver_{P}W(\lambda_1)]
$$
where $\phi$ is $\pi_* \phi(\lambda_0 \uparrow \lambda_1)$ defined in \cite[Prop 5.1.3(2),Lem 5.2.1]{Hodge-paper} \footnote{It can also be characterized as the unique (up to scalar) non-trivial map between the two Verma modules.}
up to a scalar in $\overline{\F}^{\times}_p$.
\item $$
\textnormal{gr}^1=[0 \to L_{M}(\mu_1) \to 0]
$$
\item $\textnormal{gr}^{0}=(F^{2}\BGG_{L(\lambda_1)^{\vee}})^{*}[-2]$, where $(-)^*: D^{b}(\mathcal{O}_{P,\Fpbar}) \cong D^{b}(\mathcal{O}_{P,\Fpbar})$ is the 
Serre duality functor. 
\end{enumerate}
\begin{proof}
We define $\BGG_{L(\lambda_1)}$ as the complex $\Std_{P}(L(\lambda_1))_{\chi_{\lambda_1}}$, as in 
\cite[Thm 8.2.2]{Hodge-paper}. Namely, it is the isotypic component of the standard complex $\Std_{P}(L(\lambda_1))=L(\lambda_1)\otimes \wedge^{\bullet}\mathfrak{g}/\mathfrak{p}$
where the Harish--Chandra center of $U\mathfrak{g}$ acts by the character corresponding to $\lambda_1$.
 For each term $\wedge^{k} \mathfrak{g}/\mathfrak{p} \otimes L(\lambda_1)$ fix a descending filtration 
of $P$-representations
whose graded pieces are irreducible. Therefore,  
$\BGG^{k}_{L(\lambda_1)}=\Ver_{P}(\wedge^{k} \mathfrak{g}/\mathfrak{p} \otimes L(\lambda_1))_{\chi_{\lambda_1}}$
has a filtration whose graded pieces are of the form $\Ver_{P}L_{M}(\lambda)$, for 
$\lambda \in \{\lambda_{i},\mu_{i},\nu_{i} : i=0,1\}$.
This follows from the mod $p$ Harish--Chandra isomorphism \cite{mod-p-Harish-Chandra} \cite[Lem 8.2.1]{Hodge-paper}
and the fact that these are the only $M$-dominant weights in $W_{\text{aff}} \cdot \lambda_0$
that occur as a highest weight vector of a Jordan--Holder factor
of $\wedge^{k} \mathfrak{g}/\mathfrak{p} \otimes L(\lambda_1)$ as a $P$-representation. 
This is even true for $\wedge^{k} \mathfrak{g}/\mathfrak{p} \otimes V(\lambda_1)$,
where it can be checked as follows. Let $\mu \in W_{\text{aff}} \cdot \lambda_0$ be such a weight
satisfying $L_{M}(\mu) \in \text{JH}[\wedge^{k} \mathfrak{g}/\mathfrak{p} \otimes L(\lambda_1)_{\Fpbar}]$.
We have $\wedge^k \mathfrak{g}/\mathfrak{p}=\oplus_{w \in W^{M}(k)} W(w \cdot 0)$ by \cite[\S 4.4]{Tilouine-Polo}, 
and over $\overline{\Q}_{p}$ any Jordan--Holder factor of $W(w \cdot 0) \otimes V(\lambda_1)$
is of the form $W(w \cdot 0 +\nu)$ for some weight $\nu$ of $V(\lambda_1)$ 
\cite[Cor 4.7]{Andersen-random}. 
Thus, $L_{M}(\mu) \in \text{JH}[W(w \cdot 0+\nu)_{\Fpbar}]$ for some $w \in W^{M}$ and $\nu$ as above. 
By the linkage principle for $M$ we have that $\mu \uparrow_{M} w \cdot 0+\nu$. In particular, there exists 
$x \in W$
such that $\chi:=x^{-1} \cdot (w \cdot 0+\nu)$ is $\rho$-dominant, and in $W_{\text{aff}} \cdot \lambda_0$.
 We claim that $w \cdot 0 + \nu \in \{\lambda_{i},\mu_{i},\nu_{i}\}$.
From the definition of $\chi$, we have that $x^{-1}w\rho-\rho=\chi-x^{-1}\nu$. The left-hand side 
is in $\Z^{\le 0}\Phi^{+}$, since 
$x^{-1}w\rho$ is a weight of $V(\rho)$. 
Therefore, for the right-hand side to be in $\Z^{\le 0}\Phi^{+}$ we must have $\chi \in \{\lambda_0,\lambda_1\}$,
from the explicit geometry of the alcoves. 
Thus, since $\mu \uparrow_{M} \lambda$ where 
$\lambda \in \{\lambda_{i},\mu_{i},\nu_{i} : i=0,1\}$, we must also have that $\mu \in \{\lambda_{i},\mu_{i},\nu_{i} : i=0,1\}$.
This proves the claim. 
Moreover, we may assume that in the filtration of each
$\BGG^{k}_{L(\lambda_1)}$, $\Ver_{P}W(\lambda_1)$ appears first (with lower grading), next $\Ver_{P}W(\lambda_0)$, next
$\Ver_{P}L_{M}(\mu_1)$ or $\Ver_{P}W(\mu_0)$, next $\Ver_{P}W(\nu_0)$, and finally $\Ver_{P}W(\nu_1)$. 
This can be proved by splitting extensions that are not in that order, using 
 \Cref{splitting-lemma-2}(2,3). Namely, the only case that is not covered by \Cref{splitting-lemma-2}(2) is if $W$ is an extension of 
$\Ver_{P} W(\nu_i)$ by $\Ver_{P} L_{M}(\mu_1)$ for $i=0,1$, which is covered by \Cref{splitting-lemma-2}(3).

 We define $F^{2}$ to be the subcomplex of $\BGG^{\bullet}_{L(\lambda_1)}$ by the property that the only graded pieces of the filtration appearing in $F^{2,\bullet} \subseteq \BGG^{\bullet}_{L(\lambda_1)}$
are $\Ver_{P} W(\lambda_{i})$ for $i=0,1$. This is well-defined by the ordering of the filtration 
of each term and the following observation.  If $V_1$ is a $(U\mathfrak{g},P)$-module with a 
filtration whose graded
pieces are 
$\Ver_{P} L_{M}(\chi_{i})$, and $V_2$ with a filtration with 
$\Ver_{P} L_{M}(\chi'_{j})$
as graded pieces such that $\chi'_j \ngeq \chi_i$ for all $i,j$, then there are no non-zero maps
$V_1 \to V_2$.
 Similarly, define $F^{1}$ to be the subcomplex of 
$\BGG^{\bullet}_{L(\lambda_1)}$ such that the only graded pieces appearing are $\Ver_{P}W(\lambda_i)$ and 
$\Ver_{P} L_{M}(\mu_i)$. It is well-defined by the same reason. Note that $\lambda_{1}, \mu_1, \nu_1$
appear precisely in one of $\wedge^{k} \mathfrak{g}/\mathfrak{p} \otimes V(\lambda_1)$ (hence also in one of 
$\wedge^{k} \mathfrak{g}/\mathfrak{p} \otimes L(\lambda_1)$)
with multiplicity $1$,
by the proof of lowest alcove BGG decomposition in the parabolic case \cite{Lan-Polo}. 
 Then 
\begin{enumerate}
\item $F^2=[M_0 \to M_1 \to M_2]$ where $M_0$ and $M_1$ are direct sums of $\Ver_{P} W(\lambda_0)$
by \Cref{splitting-lemma}(4), 
and $M_2$ is an extension of $\Ver_{P} W(\lambda_1)$ by some copies of $\Ver_{P} W(\lambda_0)$.
 This follows by \Cref{splitting-lemma-2} and the fact that $\lambda_1$ 
only appears as a weight of $\wedge^{k} \mathfrak{g}/\mathfrak{p} \otimes V(\lambda_1)$ for $k=0$, 
where it does with multiplicity $1$.
 We now show that $F^2$ is quasi-isomorphic to the complex in the proposition. 
 Note that the first differential is linear, since $M_0$ and $M_1$ have the same 
highest weight vectors. 
Using the long exact sequence associated to $F^2 \to \BGG \to \BGG/F^2$ 
we have $\mathcal{H}^{0}(F^2)=0$ and $\mathcal{H}^1(F^2) \cong \mathcal{H}^0(\BGG/F^2)$. 
We are using cohomological grading, so that $\mathcal{H}^{\bullet}\BGG_{L(\lambda_1)}$ is concentrated in degree $2$. 
From the first equality we can assume that $M_0=0$, since $M_1/M_0$ will 
still be some copies of $\Ver_{P} W(\lambda_0)$.
Since the map from $M_1$ to the projection of $M_2$ consisting of $\Ver_{P} W(\lambda_0)$s is linear, we may also 
assume that $M_1$ maps to $\Ver_{P} W(\lambda_1)$. 
All the copies of $\Ver_{P} W(\lambda_0)$ in $M_1$ will map to
$\Ver_{P} W(\lambda_1)$ via the unique non-trivial map or as $0$, by \cite[Cor 5.2.4(2)]{Hodge-paper} and Serre duality.
Therefore, we see that if $M_1$ contains at least two copies of $\Ver_{P} W(\lambda_0)$ then $\mathcal{H}^1(F^2)$
contains $\Ver_{P} W(\lambda_0)$. Since by construction $\mathcal{H}^0(\BGG/F^2)$ does not contain 
the weight $\lambda_0$, we may assume that $M_1 \subseteq \Ver_{P} W(\lambda_0)$.
Moreover, by looking at \Cref{character-GL3} and the fact that $\BGG/F^2$ does not contain the weight
$\lambda_0$ we see that the formal Euler characteristic of $F^2$ must contain $-\Ver_{P} W(\lambda_0)$, so that
$M_1\neq 0$, and hence $M_1=\Ver_{P} W(\lambda_0)$.
 Again by 
considering the Euler characteristic, 
we must then have that $M_2=\Ver_{P} W(\lambda_1)$. 
Thus,
$$
F^2=[0 \to \Ver_{P} W(\lambda_0) \xrightarrow{\phi} \Ver_{P} W(\lambda_1)]
$$
for some $\phi$ which cannot be zero, by the equality $\mathcal{H}^1(F^2) \cong \mathcal{H}^0(\BGG/F^2)$.
 Thus, it is the unique
non-trivial map. 
\item We claim that $\text{gr}^{0} F^{\bullet}$ is Serre dual of $F^2 \BGG_{L(\lambda_1)^{\vee}}$ up to a shift.
 We observe that $\Std_{P}L(\lambda_1)^{\vee}$
is Serre dual to $\Std_{P}L(\lambda_1)$ up to a shift. Here we are using that for $\GL_3$ 
all the weights appearing in $\wedge^k \mathfrak{g}/\mathfrak{p}=\oplus_{w \in W^M(k)} W(w \cdot 0)$ are automatically
$p$-restricted 
for $M$. 
Therefore, $\BGG_{L(\lambda_1)^{\vee}}$ is Serre dual to 
$\BGG_{L(\lambda_1)}$, since they are quasi-isomorphic.
We can and do choose the filtrations of $\wedge^{\bullet} \mathfrak{g}/\mathfrak{p} \otimes L(\lambda_1)$ and $
\wedge^{\bullet} \mathfrak{g}/\mathfrak{p} \otimes L(\lambda_1)^{\vee}$ to be Serre dual to each other. 
We use the notation $\lambda^{\vee}:=-w_{0}\lambda$.
Observe that Serre duality sends $W(\lambda^{\vee}_1)$ to $W(\nu_1)$, 
$W(s_{\alpha_1+\alpha_2,1} \cdot \lambda^{\vee}_1)$ to $W(\nu_0)$ and so on. 
Write $\text{gr}^{0} F^{\bullet}_{L(\lambda_1)^{\vee}}=[V_1 \to V_2 \to V_3]$, then its Serre dual shifted by $2$ is 
$[V^*_3 \to V^*_2 \to V^*_1]$ by \Cref{lemma-serre-dual-acyclic}, and $V^*_i$ is a subobject of $\BGG^{3-i}_{L(\lambda_1)^{\vee}}$. 
Since $(-)^*$ is exact, it preserves the abelian subcategory of objects that are extensions of Verma modules, and the fact that Serre duality exchanges the Verma modules appearing in
the respective complexes, we see that $V^*_i$ is identified with 
$F^{2} \BGG^{3-i}_{L(\lambda_1)^{\vee}}$. This crucially uses that $\BGG_{L(\lambda_1)}$ and its filtration are defined as genuine complexes. 
This argument shows that the filtrations $F^{\bullet}$ are compatible with Serre duality.
In particular,
$$
\text{gr}^{0} F^{\bullet}=[\Ver_{P} W(\nu_1) \xrightarrow{\phi_{0}^{*}} \Ver_{P} W(\nu_0) \to 0]
$$
for the unique non-zero map $\phi_0: \Ver_{P} W(s_{\alpha_1+\alpha_2,1} \cdot \lambda^{\vee}_1) \to \Ver_{P} W(\lambda^{\vee}_1)$
appearing in $F^2 \BGG_{L(\lambda^{\vee}_1)}$. 
\item By a similar reasoning as above, we get 
$$
\text{gr}^1 F^{\bullet}=[V_0 \xrightarrow{\psi} V_1 \xrightarrow{\psi_1} V_2]
$$
where $V_1$ is an extension of some copies of $\Ver_{P} W(\mu_0)$, followed by 
$\Ver_{P} L_{M}(\mu_1)$, and then more copies of $\Ver_{P} W(\mu_0)$, and $V_{0}, V_2$
are copies of $\Ver_{P} W(\mu_0)$. Note that all the maps are linear. By the same argument as in $(2)$ 
we have that $\text{gr}^{1} F^{\bullet}$ is Serre dual to $\text{gr}^{1} F^{\bullet}_{L(-w_{0}\lambda_1)}$ up to a shift. 
We have the exact sequence $\mathcal{H}^0(F^1) \to \mathcal{H}^{0}(\text{gr}^1) \to \mathcal{H}^1(F^2)$. 
The first term is zero since it embeds into $\mathcal{H}^0(\BGG)=0$, and the last is the kernel of the non-trivial map 
$\Ver_{P}W(\lambda_0) \to \Ver_{P}W(\lambda_1)$. Suppose that $\mathcal{H}^{0}(\text{gr}^1) \neq 0$, then $\Ver_{P}W(\mu_0)$
embeds into $\Ver_{P}W(\lambda_0)$. However, there is a unique non-zero map 
$\Ver_{P}W(\mu_0) \to \Ver_{P}W(\lambda_0)$. Moreover, it is not injective, 
since its kernel is $\Ver_{P}W(\nu_0)$, by the lowest alcove BGG resolution \cite[Thm 5.2]{Lan-Polo}.
Therefore, $\psi$ is injective. Since the projection of $V_0$ to the copies of $\Ver_P W(\mu_0)$ in $V_1$ 
is linear, we can assume that $V_0$ maps to $\Ver_P L_{M}(\mu_1)$, but such a map must be zero, so 
$V_0$ maps to the first copies of $\Ver_{P} W(\mu_0)$ in $V_1$, and that map is injective. Thus, we 
can assume that $V_0=0$. By duality with $\text{gr}^{1} F^{\bullet}_{L(-w_{0}\lambda_1)}$ we also have that
$\mathcal{H}^2(\text{gr}^1)=0$, so we can assume that $V_2=0$. Finally, by comparing the Euler characteristic 
of $F^{\bullet}\subseteq \BGG_{L(\lambda_1)}$ and \Cref{character-GL3} we conclude 
that $V_1=\Ver_{P} L_{M}(\mu_1)$.
\end{enumerate}
\end{proof} 
\end{prop}

\subsection{The case of $\GSp_4$}
Now let $G=\GSp_4$, with $P$ the Siegel parabolic parametrizing Lagrangians, and let $\lambda_0 \in C_0$.
Let $\beta$ be the long root and $\alpha$ the short root, characterized by $\alpha \in \Phi_{M}$. 
Let $\lambda_1=s_{\alpha+\beta,1} \cdot \lambda_0$, $\lambda_2=s_{2\alpha+\beta,1} \cdot \lambda_1$, 
$\lambda_3=s_{\alpha+\beta,2} \cdot \lambda_2$, $\mu_i=s_{\beta,0} \cdot \lambda_i$, 
$\nu_i=s_{\alpha+\beta,0} \cdot \mu_i$, and $\epsilon_i=s_{2\alpha+\beta,0} \cdot \nu_i$ be the corresponding 
reflections in various $M$-dominant alcoves, as in \Cref{figure6}. 
\begin{figure}[h]
  \centering
  \begin{tikzpicture}[scale=1.2]
    \coordinate (A) at (0,0);
    \coordinate (B) at (2,0);
    \coordinate (C) at (2,2);
    \coordinate (D) at (4,2);
    \coordinate (P) at (1,1);
    \coordinate (Q) at (3,1);
    
    \draw[->, red] (0,0) -- (0,1.2) node[right] {$\beta$};
    \draw[->, red] (0,0) -- (-0.6,0.6) node[above] {$-\alpha$};
    \draw[dotted] (0,0)--(4,0);
    \draw[dotted] (0,0)--(0,-4);
    \draw[dotted] (0,0)--(-3,-3);
    \draw[dotted] (0,0)--(3,-3);
    \draw[dotted] (0,0)--(3,3);

    \draw (A) -- (C); 
    \draw (A) -- (B); 
    \draw (B) -- (D); 
    \draw (C) -- (D); 
    \draw (B) -- (P); 
    \draw (C) -- (Q); 
    \draw (B) -- (C); 
    \draw (A) --(2,-2);
    \draw (B) --(2,-2);
    \draw (B) --(1,-1);
    \draw (B) --(4,-2);
    \draw (2,-2)--(4,-2);
    \draw (2,-2)--(3,-1);
    \draw (A)--(0,-2);
    \draw (0,-2)--(2,-4);
    \draw (2,-2)--(2,-4);
    \draw (1,-1)--(0,-2);
    \draw (1,-3)--(2,-2);
    \draw (0,-2)--(2,-2);
    \draw (A) --(-2,-2);
    \draw (-2,-2)--(-2,-4);
    \draw (-2,-4)--(0,-2);
    \draw (-1,-1)--(0,-2);
    \draw (-2,-2)--(-1,-3);
    \draw (-2,-2)--(0,-2);
    
    \node at (0,0) {$\bullet$};
    \node at (-0.35,0) {$-\rho$};
    \node at (1.5,0.25) {$\cdot$};
    \node[left] at (1.5,0.25) {$\scriptstyle \lambda_0$};
    \node at (1.75,0.5) {$\cdot$};
    \node at (1.6,0.65) {$\scriptstyle \lambda_1$};
    \node at (2.25,0.5) {$\cdot$};
    \node at (2.15,0.69) {$\scriptstyle \lambda_2$};
    \node at (3.5,1.75) {$\cdot$};
    \node[left] at (3.5,1.75) {$\scriptstyle \lambda_3$};
    \node at (1.5,-0.25) {$\cdot$};
    \node[left] at (1.5,-0.25) {$\scriptstyle \mu_0$};
    \node at (1.75,-0.5) {$\cdot$};
    \node[below] at (1.75,-0.5) {$\scriptstyle \mu_1$};
    \node at (2.25,-0.5) {$\cdot$};
    \node[below] at (2.25,-0.5) {$\scriptstyle \mu_2$};
    \node at (3.5,-1.75) {$\cdot$};
    \node[left] at (3.5,-1.75) {$\scriptstyle \mu_3$};
    \node at (0.25,-1.5) {$\cdot$};
    \node[above] at (0.25,-1.5) {$\scriptstyle \nu_0$};
    \node at (0.5,-1.75) {$\cdot$};
    \node[right] at (0.5,-1.75) {$\scriptstyle \nu_1$};
    \node at (0.5,-2.25) {$\cdot$};
    \node[right] at (0.5,-2.25) {$\scriptstyle \nu_2$};
    \node at (1.75,-3.5) {$\cdot$};
    \node[above] at (1.75,-3.5) {$\scriptstyle \nu_3$};
    \node at (-0.25,-1.5) {$\cdot$};
    \node[above] at (-0.25,-1.5) {$\scriptstyle \epsilon_0$};
    \node at (-0.5,-1.75) {$\cdot$};
    \node[left] at (-0.5,-1.75) {$\scriptstyle \epsilon_1$};
    \node at (-0.5,-2.25) {$\cdot$};
    \node[left] at (-0.5,-2.25) {$\scriptstyle \epsilon_2$};
    \node at (-1.75,-3.5) {$\cdot$};
    \node[above] at (-1.75,-3.5) {$\scriptstyle \epsilon_3$};

\end{tikzpicture} 
  \caption{}
  \label{figure6}
  \end{figure}

\begin{lemma} \label{character-GSP4}
We have the following relations in the Grothendieck group of $\Rep_{\Fpbar}G$
$$
[V(\lambda_0)]=L(\lambda_0), \; [V(\lambda_1)]=L(\lambda_1)+L(\lambda_0), \;  [V(\lambda_2)]=L(\lambda_2)+L(\lambda_1).
$$
We also have the following characters for $L(\lambda_i)$ for $\lambda_{i} \in X_{1}(T)$ according to their alcoves. 
\begin{enumerate}
\item $[L(\lambda_0)]=\Ver_{P} W(\lambda_0) - \Ver_{P} W(\mu_0)+ \Ver_{P} W(\nu_0)-\Ver_{P} W(\epsilon_0)$
\item $[L(\lambda_1)]=\Ver_{P} W(\lambda_1)- \Ver_{P} W(\lambda_0)- \Ver_{P} L_{M}(\mu_1)+
 \Ver_{P} L_{M}(\nu_1)-\Ver_{P} W(\epsilon_1)+ \Ver_{P} W(\epsilon_0)$
 \item $[L(\lambda_2)]=\Ver_{P} W(\lambda_2)- \Ver_{P} W(\lambda_1)
 -\Ver_{P} L_{M}(\mu_2)+\Ver_{P} L_{M}(\mu_1)+\Ver_{P} L_{M}(\nu_2)-\Ver_{P} L_{M}(\nu_1) \break 
 -\Ver_{P} W(\epsilon_2)+ \Ver_{P} W(\epsilon_1).$
\end{enumerate}
\begin{proof}
The decomposition of dual Weyl modules is from \cite[\S 7]{Jantzen-decomp} \cite[Prop 2.10]{paper}. To get the characters in terms of Verma modules we use the parabolic BGG decomposition in characteristic $0$
$$
0 \to \Ver_{P} W(\epsilon_i) \to \Ver_{P} W(\nu_i) \to \Ver_{P} W(\mu_i) \to \Ver_{P} W(\lambda_i) \to V(\lambda_i)_{\overline{\Q}_p} \to 0,
$$
and that $[W(\mu_2)]=W(\lambda_0)+L_{M}(\mu_2), 
[W(\mu_1)]=L_M(\mu_1)+W(\mu_0)$ in the Grothendieck group. 
\end{proof}
\end{lemma}

\begin{prop} \label{BGG-lambda1-gsp4}
Let $\lambda_1 \in C_1$.
There exists a resolution $\BGG_{L(\lambda_1)}$ of $L(\lambda_1)$ in $\mathcal{O}_{P,\Fpbar}$
 with a decreasing $4$-step filtration $F^{\bullet}$ in $D^{b}(\mathcal{O}_{P,\Fpbar})$ such that 
\begin{enumerate}
\item $F^3=[0 \to 0 \to \Ver_P W(\lambda_0) \xrightarrow{\phi} \Ver_P W(\lambda_1)]$
with $\phi$ a non-zero map.
\item $\textnormal{gr}^2=[0 \to 0 \to \Ver_{P} L_{M}(\mu_1) \to 0]$
\item For $i=0,1$, $\textnormal{gr}^{i}=(\textnormal{gr}^{3-i} \BGG_{L(-w_0\lambda_1)})^*[-3]$.
\end{enumerate}
\begin{proof}
We define $\BGG_{L(\lambda_1)}$ as a localization as in the proof of \Cref{GL3-BGG}, and we follow the same method.
The Verma modules of highest weight $\{\lambda_i,\mu_i,\nu_i,\epsilon_i\}$ for $i=0,1$ are the only ones appearing in $\BGG_{L(\lambda_1)}$, 
by the same argument as in the proof of \Cref{GL3-BGG}, using that all dominant affine Weyl reflections of $\lambda_1$ 
which are not $\lambda_0$ 
are $\ge \lambda_1$. 
We define descending filtrations of $P$-representations on each 
$\wedge^k \mathfrak{g}/\mathfrak{p} \otimes L(\lambda_1)$ with irreducible graded pieces. Then we can assume that in the filtration of each $\BGG^k_{L(\lambda_1)}$
$\Ver_{P}W(\lambda_{1})$ appear first (with lower grading), then $\Ver_{P}W(\lambda_{0})$, 
then $\Ver_{P}W(\mu_{0})$ or $\Ver_{P}L_{M}(\mu_1)$, then $\Ver_{P}W(\nu_{0})$ or $\Ver_{P}L_{M}(\nu_{1})$, then $\Ver_{P}W(\epsilon_{0})$, 
and finally $\Ver_{P}W(\epsilon_{1})$.
This can be ensured as in the proof of \Cref{GL3-BGG}. 
We define $F^3 \subseteq \BGG_{L(\lambda_1)}$ to be the filtration 
containing all the $\Ver_{P} W(\lambda_i)$ as graded pieces, $F^2$ the one containing all the $\Ver_{P} W(\lambda_i)$
and $\Ver_{P} L_M(\mu_i)$, and so on. 
This is a well-defined complex since the highest weights of each graded piece 
are not weakly increasing in the $\ge$ order. The graded pieces of $F^{\bullet} \BGG_{L(\lambda_1)}$ and 
 $F^{\bullet} \BGG_{L(\lambda_1)^{\vee}}$ satisfy the same Serre duality relations as before. There we are using 
 that the standard complex is Serre dual to itself up to shift for $p \ge 3$, which is automatically implied by the 
 existence of $\lambda_1 \in C_1$.
\begin{enumerate}
 \item We get $F^3=[M_0 \to M_1 \to M_2
\to M_3]$, where $M_{0}, M_1,M_2$ are copies of  $\Ver_{P} W(\lambda_0)$, and $M_3$ is an extension of 
$\Ver_{P} W(\lambda_1)$ by copies of  $\Ver_{P} W(\lambda_0)$, 
by \Cref{splitting-lemma-2} and \Cref{splitting-lemma}. We want to prove 
that it is quasi-isomorphic to the one in the proposition.
Note that the first two maps are linear, in the sense that each map of the form
 $\Ver_{P}W(\lambda_0) \to \Ver_P W(\lambda_0)$ is a scalar.
By taking the exact sequence associated to $F^3 \to \BGG \to \BGG/F^3$ we get 
$\mathcal{H}^0(F^3)=0$, and $\mathcal{H}^{i}(F^3)\cong \mathcal{H}^{i-1}(\BGG/F^3)$ for $i=1,2$.
From $\mathcal{H}^0(F^3)=0$ we may assume that $M_0=0$, since $M_1/M_0$ will again be a sum of  
$\Ver_P W(\lambda_0)$. If $\mathcal{H}^{1}(F^3) \neq 0$ then it must contain $\Ver_{P} W(\lambda_0)$ as a graded piece.
However, $\mathcal{H}^{0}(\BGG/F^3)$ only contains lower weights in the $\le$ order, 
so this forces $\mathcal{H}^{1}(F^3)=0$, and we can assume that $M_1=0$ 
as before. Consider a summand $\Ver_P W(\lambda_0)$
of $M_2$. If the projection to the quotient of $M_3$ consisting of $\Ver_{P} W(\lambda_0)$ is non-zero,
it must be 
injective, so that it defines a splitting in $M_3$, and we can remove both summands. Therefore, 
we may assume that $M_2$ maps to $\Ver_{P} W(\lambda_1)$. 
At the same time there is an exact sequence 
$$
0 \to \mathcal{H}^2(\BGG/F^3) \to \mathcal{H}^3(F^3) \to L(\lambda_1) \to \mathcal{H}^3(\BGG/F^3) \to 0.
$$
By weight considerations, and since $L(\lambda_1)$ is irreducible 
as a $(U\mathfrak{g},P)$-module, the last term is $0$, and $\mathcal{H}^2(\BGG/F^3)$ does not contain 
a vector of weight $\lambda_0$. 
By considering the Euler characteristic of $\BGG$
and \Cref{character-GSP4} we deduce that $M_2 \neq 0$.
The map from the subobject of $\mathcal{H}^3(F^3)$ which is a quotient of $\Ver_P W(\lambda_1)$ to $L(\lambda_1)$ is 
induced by the unique non-zero map $\Ver_{P} W(\lambda_1) \to L(\lambda_1)$ (by the universal 
property of Verma and Weyl modules). This non-zero map is surjective, since 
$\lambda_1$ is $p$-restricted, which implies that $L(\lambda_1)=U\mathfrak{g}\langle v_{\lambda_1}\rangle$.
Therefore, $M_3=\Ver_P W(\lambda_1)$, otherwise there would be a vector of weight $\lambda_0$
in the kernel of $\mathcal{H}^3(F^3) \to L(\lambda_1)$.
Again, by comparing the Euler characteristic of $F^3$ with the one of $L(\lambda_2)$ using \Cref{character-GSP4}
we conclude that $M_2=\Ver_P W(\lambda_0)$. Moreover, the map $\phi$ cannot be zero, otherwise 
$\mathcal{H}^{2}(F^3) \cong \Ver_{P} W(\lambda_0) \cong \mathcal{H}^1(\BGG/F^3)$. 

\item We have $\text{gr}^2=[V_0 \to V_1 \to V_2 \to V_3]$ where 
$V_{0}, V_{1},V_{3}$ are direct sums of $ \Ver_{P} W(\mu_0)$, and $V_2$ is an extension of $\Ver_{P} L_M(\mu_1)$
by some copies of $\Ver_{P} W(\mu_0)$, in some order. 
 Note that all the maps are linear. 
By considering the long exact sequence for $\text{gr}^2 \to \BGG/F^3 \to \BGG/F^2$ and 
the fact that $\mathcal{H}^{0}(\BGG/F^3)=0$ from above, we get that $\mathcal{H}^{0}(\text{gr}^2)=0$,
so we can assume that $V_0=0$. From the 
same long exact sequence and point $(1)$ above we have the exact sequence 
$$
\mathcal{H}^2(\BGG/F^2) \to \mathcal{H}^3(\text{gr}^2) \to \mathcal{H}^3(\BGG/F^3)=0.
$$ 
Therefore, $\mathcal{H}^3(\text{gr}^2)=0$, since if it were non-zero it would contain the weight $\mu_0$, 
and $\mathcal{H}^2(\BGG/F^2)$ does not. By quotienting by some copies of $\Ver_{P}W(\mu_0)$
we can assume that $V_3=0$. We can also assume that $V_1 \to V_2$ is zero. We have an exact sequence 
$$
0 \to \mathcal{H}^0(\BGG/F^2) \to \mathcal{H}^1(\text{gr}^2) \to \mathcal{H}^1(\BGG/F^3) \cong \mathcal{H}^2(F^3),
$$
and the last term embeds into $\Ver_{P}W(\lambda_0)$, with non-trivial cokernel. There is a unique non-zero map 
$\psi: \Ver_{P}W(\mu_0) \to \Ver_{P}W(\lambda_0)$
by \cite[Cor 5.2.4]{Hodge-paper}, so we can assume that $V_1$ contains at most one copy of 
$\Ver_{P}W(\mu_0)$, since $\mathcal{H}^0(\BGG/F^2)$ does not contain the weight $\mu_0$.
Suppose that $V_1 \neq 0$.
By the same reason the map $\mathcal{H}^1(\text{gr}^2) \to \mathcal{H}^1(\BGG/F^3) \subseteq \Ver_{P}W(\lambda_0)$
is non-zero,
so it is induced by $\psi$. 
Thus, $\mathcal{H}^0(\BGG/F^2)$ is the kernel of $\psi$. By the lowest alcove 
BGG decomposition it 
is also equal to the image of a non-zero map $\Ver_{P}W(\nu_0) \to \Ver_{P}W(\mu_0)$, in particular it contains 
the weight $\nu_0$. On the other hand, we have the exact sequence 
$0=\mathcal{H}^0(\text{gr}^1) \to \mathcal{H}^0(\BGG/F^2) \to \mathcal{H}^0(\text{gr}^0)$, where the first term vanishes by Serre duality and the knowledge of $\text{gr}^2$.
This is a contradiction 
since the last term does not contain the weight $\nu_0$. This proves that $V_1=0$. Moreover,
$V_2=\Ver_{P}L_{M}(\mu_1)$, since by considering the Euler characteristic of $\text{gr}^2$,
$V_2$ contains as many copies of $\Ver_{P}W(\mu_0)$
as $V_1$.
\end{enumerate}
\end{proof}
\end{prop}

Finally, we consider a resolution for the $C_2$ alcove. The shape of the complexes in $(1)$ and $(2)$ will be the key for \Cref{correct-entailment} ahead.
\begin{prop} \label{BGG-lambda2-gsp4}
	There exists a resolution $\BGG_{L(\lambda_2)}$ of $L(\lambda_2)$
	in  $\mathcal{O}_{P,\Fpbar}$ with a $4$-step  decreasing filtration $F^{\bullet}$  in $D^b(\mathcal{O}_{P,\Fpbar})$ such that for $\lambda_2$ generic
	\begin{enumerate}
		\item $F^3=[0 \to 0 \to \Ver_P W(\lambda_1) \xrightarrow{\phi} M_3]$, 
	where $M_3$ is an extension of $\Ver_{P} W(\lambda_2)$ by  $\Ver_{P} W(\lambda_0)$, 
  and $\phi$ lands in $\Ver_P W(\lambda_2)$,
	where it is equal to $\pi_* \phi(\lambda_1 \uparrow \lambda_2)$ defined as in \cite[Prop 5.1.3, Lem 5.1.2]{Hodge-paper}.
	\item $\textnormal{gr}^2=[0 \to  \Ver_{P} W(\mu_1) \to V_2 \to 0]$,
	where $V_2$ is an extension of $\Ver_{P} W(\mu_2)$  by $\Ver_{P} W(\mu_0)$.  
  Furthermore, the map $ \Ver_{P} W(\mu_1)  \to V_2 \to \Ver_{P} W(\mu_0)$ is trivial, and the induced map $\Ver_{P} L_{M}(\mu_1) \to \Ver_{P} L_{M}(\mu_2)$ is the unique non-trivial map. 
	\item For $i=0,1$, $\textnormal{gr}^{i}=(\textnormal{gr}^{3-i} \BGG_{L(\lambda^{\vee}_2)})^*[-3]$.
  \end{enumerate}
	\begin{proof}
	We use the same method as in the proof of \Cref{GL3-BGG}.
  As in \Cref{BGG-lambda1-gsp4} the only Verma modules appearing in $\BGG_{L(\lambda_2)}$
  are the ones with highest weights $\{\lambda_i,\mu_i,\nu_i,\epsilon_i\}$ for $i=0,1,2$.
   Here we see that $\mu_2 \ge \lambda_0$,
	so we need to be more careful when defining the filtration. Define $\BGG_{L(\lambda_2)}$ as a complex 
	as in the proof of \Cref{GL3-BGG}. Fix a descending filtration of $P$-representations for 
  each $\wedge^{k} \mathfrak{g}/\mathfrak{p} \otimes L(\lambda_2)$ whose graded pieces are irreducible. 
  We can assume that on the associated filtration for each $\BGG^k_{L(\lambda_2)}$ the graded pieces appear in the following order. 
  First $\Ver_{P}W(\lambda_2)$ (with lower grading), 
  next $\Ver_{P}W(\lambda_1)$, next $\Ver_{P}W(\lambda_0)$ or $\Ver_{P}L_{M}(\mu_2)$, next
  $\Ver_{P}L_{M}(\mu_1)$ or $\Ver_{P}W(\mu_0)$, next $\Ver_{P}W(\nu_0)$ or $\Ver_{P}L_{M}(\nu_1)$, 
  next $\Ver_{P}L_{M}(\nu_2)$ or $\Ver_{P}W(\epsilon_0)$, next $\Ver_{P}W(\epsilon_1)$, and finally 
  $\Ver_{P}W(\epsilon_2)$. 
  This is ensured by \Cref{splitting-lemma-2}(2,3) and the explicit geometry of the affine Weyl reflections appearing. 
   Recall that $\lambda_2$ only appears on $\wedge^{k} \mathfrak{g}/\mathfrak{p} 
	\otimes L(\lambda_2)$ for $k=0$ with multiplicity $1$, and similarly for $\mu_2$ with $k=1$. 
  Define $F^3$ to be the subcomplex $[M_0 \to M_1 \to M_2
	\to M_3]$ of
  $\BGG_{L(\lambda_2)}$ as follows. We let $M_{i}$ for $i=0,1,2$ be the subobject of 
	$\BGG^{i}$ containing all of the $\Ver_{P} W(\lambda_1)$.
	We define $M_3$ to be the subobject of $\BGG^3$ containing 
	all the $\Ver_{P} W(\lambda_i)$, so it is an extension of $\Ver_{P}W(\lambda_2)$ followed by
	$N_3$ consisting of copies of $\Ver_{P}W(\lambda_1)$, and then by $T_3$ consisting of copies of $\Ver_{P}W(\lambda_0)$. 
	Then $F^3$ is a well-defined complex by 
	\Cref{splitting-lemma} and the construction of the filtration on $\BGG^k$.
	Define $F^2$ to be the subcomplex generated by $F^3$ and containing all of the remaining 
	$\Ver_{P} W(\lambda_0)$ and $\Ver_P L_{M}(\mu_i)$.
   Finally, dually let 
	$F^1$ be generated by $F^2$, all the possible $\Ver_{P} L_{M}(\nu_i)$, and all the $\Ver_{P} W(\epsilon_0)$
	appearing in the last $3$ terms. We show that the graded pieces of this filtration are quasi-isomorphic to the complexes in the proposition. 
  \par

	As in \Cref{BGG-lambda1-gsp4} we get $\mathcal{H}^{0}(F^3)=0$ and 
	$\mathcal{H}^{i}(F^3) \cong \mathcal{H}^{i-1}(\BGG/F^3)$ for $i=1,2$. From the first one we 
	can assume that $M_0=0$.
	If $\mathcal{H}^{1}(F^3) \cong \mathcal{H}^{0}(\BGG/F^3)$ is non-zero it must contain the weight $\lambda_1$, 
	 which does not appear in $\mathcal{H}^{0}(\BGG/F^3)$, therefore we may also assume that $M_1=0$. 
   We have that $M_2$ maps trivially to $T_3$.
	 If a summand $\Ver_{P} W(\lambda_1)$ of $M_2$ maps isomorphically to a copy of 
	 $\Ver_{P} W(\lambda_1)$ in $N_3$, then it defines a splitting of the extension of that copy with 
	 $\Ver_{P} W(\lambda_2)$, so in that case $F^3$ is quasi-isomorphic to the same complex 
   with the two copies of 
	 $\Ver_{P} W(\lambda_1)$ removed. Therefore, we may assume that $M_2$ maps to $\Ver_{P} W(\lambda_2)$. There is 
     only one non-zero map $\Ver_P W(\lambda_1) \to \Ver_P W(\lambda_2)$ by Serre duality
     and the uniqueness of maps between their Serre duals \cite[Cor 5.2.4]{Hodge-paper}.
	 Thus, if $M_2$ has more than one copy 
	 of $\Ver_{P}W(\lambda_1)$, the weight $\lambda_1$ would appear in 
	 $\mathcal{H}^{2}(F^3) \cong \mathcal{H}^{1}(\BGG/F^3)$. Thus, we may assume $M_2 \subseteq \Ver_{P} W(\lambda_1)$, 
   and comparing the Euler characteristic of $F^3\hookrightarrow \BGG$ with \Cref{character-GSP4} implies that 
   $M_2=\Ver_{P} W(\lambda_1)$ and $N_3=0$. 
	 Moreover, the map from $M_2$ to $\Ver_{P} W(\lambda_2)$ is $\phi$.
    Let $T_3=\Ver_{P} W(\lambda_0)^{\oplus n}$ for some $n$. 
	 We can determine $n$ by computing the dimension of the $\lambda_0$-weight space in the Euler 
	 characteristic of $F^3$. Let $w_{\lambda}(\mu)$ be the dimension of the $\lambda$ weight space in 
	 $\Ver_{P} W(\mu)$. By construction of $F^3$ (as a complex in the abelian category!)
	 \begin{align} \label{big-computation}
	 n+w_{\lambda_0}(\lambda_2)-&w_{\lambda_0}(\lambda_1)=\sum^{2}_{i=0}[L(\lambda_2):W(\lambda_{i})]
   w_{\lambda_0}(\lambda_i)-
	 [\mathfrak{g}/\mathfrak{p} \otimes L(\lambda_2):W(\lambda_1)]w_{\lambda_0}(\lambda_1) \\
	 & +[\wedge^{2} \mathfrak{g}/\mathfrak{p} \otimes L(\lambda_2):W(\lambda_1)]w_{\lambda_0}(\lambda_1)
	 - [\wedge^{3} \mathfrak{g}/\mathfrak{p} \otimes L(\lambda_2):W(\lambda_1)]
   w_{\lambda_0}(\lambda_1).
	 \end{align}
    By a laborious computation 
    that we relegate to \Cref{euler-char-computation} we conclude that $n=1$. This proves $(1)$.

	 By its construction we have 
	 $\text{gr}^2=[V_0 \to V_1 \to V_2 \to V_3]$, where $V_{0,1}$ is an extension of $N_{0,1}$,
    consisting of copies of $\Ver_{P} W(\lambda_0)$, 
	 and $T_{0,1}$, consisting of extensions of $\Ver_{P} L_{M}(\mu_{1})$ and $\Ver_{P}W(\mu_0)$.
  Similarly, $V_2$ is an extension of
	 $N_2$ consisting of copies of $\Ver_{P}W(\lambda_0)$, then $\Ver_{P} L_{M}(\mu_2)$,
  then $N'_2$, 
	 consisting of copies of $\Ver_{P} W(\lambda_0)$, and then $T_2$ as before. 
	  Finally, $V_3$ is an extension 
	of copies of $\Ver_{P} W(\mu_{0})$ and $\Ver_{P} L_{M}(\mu_{1})$. By the long exact sequence for 
  $\text{gr}^2 \to \BGG/F^3 \to \BGG/F^2$ and the properties of $F^3$ obtained above
	we get that $\mathcal{H}^{0}(\text{gr}^2)=0$ and an exact sequence
	$0 \to \mathcal{H}^{0}(\BGG/F^2) \to \mathcal{H}^{1}(\text{gr}^2) \to \mathcal{H}^2(F^3)$. 
  On the other hand $\mathcal{H}^2(F^3)=0$, since one can check that by the uniqueness of such a map of Verma modules, 
  the map $\phi$ in $F^3$ is obtained by "$p$-translating" the non-trivial map 
  $\Ver_{P}W(\epsilon_0) \to \Ver_{P}W(\nu_0)$ by $p(\alpha+\beta)$, which is injective by the lowest alcove BGG 
  decomposition. By $p$-translating we mean that $W(\lambda_1)=W(\epsilon_0) \otimes W(p(\alpha+\beta))$, so we send $v_{\lambda_1}$
  to the image of $v_{\epsilon_0}$ times the vector generating the character $W(p(\alpha+\beta))$.
  From $\mathcal{H}^{0}(\text{gr}^2)=0$ we may assume that $N_0=0$, since it is automatic that $N_0$ maps trivially to 
  $T_1$. Moreover, we may assume 
	that $V_0$ maps to $N_1$ and that 
	$V_0$ does not contain any copies of $\Ver_{P} L_{M}(\mu_1)$, since $\mu_1 \nleq \lambda_0$,
  so any vector of 
  weight $\mu_1$ would map to $0$. By the lowest alcove BGG decomposition the unique non-trivial map 
  $\Ver_{P}W(\mu_0) \to \Ver_{P}W(\lambda_0)$ is not injective, so we may also assume that $V_0=0$.

  We have that $N_1$ maps trivially to $T_2$.
	If a summand $\Ver_{P} W(\lambda_0)$ of 
	$N_1$ maps isomorphically to a summand of $N'_2$ it defines a splitting of the extension with 
	$\Ver_{P}L_{M}(\mu_2)$, in which case we can quotient by such a map. Thus, we can assume that $N_1$ maps 
	to $N_2$. Using that $\lambda_0$ does not appear in 
	$\mathcal{H}^{1}(\text{gr}^2) \cong \mathcal{H}^{0}(\BGG/F^2)$ the map $N_1 \to N_2$ is injective, 
	so we may assume that $N_1=0$. 
	By comparing the Euler characteristic of $F^2$ and $\BGG_{L(\lambda_2)}$
	using \Cref{character-GSP4} we deduce that there is exactly one copy of $\Ver_{P} W(\lambda_0)$
	in $V_2$. By going over the reduction step to assume that 
	$V_0=0$, the multiplicity of $\Ver_{P} W(\mu_0)$ in $V_1$ can be 
	computed as the multiplicity of $W(\mu_0)$ in $\wedge^{2} \mathfrak{g}/\mathfrak{p} \otimes L(\lambda_2)$
	minus the multiplicity of $W(\mu_0)$ in $\wedge^{3} \mathfrak{g}/\mathfrak{p} \otimes L(\lambda_2)$.
	By \Cref{euler-char-computation}(2) this multiplicity is $1$,
   and the multiplicity of 
	$\Ver_{P}L_{M}(\mu_1)$ in $V_1$ 
	is also $1$. On the other hand, using the previous long exact sequence we get 
  an exact sequence $\mathcal{H}^2(\BGG/F^2) \to \mathcal{H}^3(\text{gr}^2) \to \mathcal{H}^3(\BGG/F^3)=0$, 
  and the extension of $\Ver_{P}L_{M}(\mu_2)$ and 
  $\Ver_{P}W(\lambda_0)$ in $V_2$ maps trivially to $V_3$. Since $\mu_{0,1}$ don't appear in 
  $\mathcal{H}^2(\BGG/F^2)$ we have $\mathcal{H}^3(\text{gr}^2)=0$. Suppose that $V_1$ maps trivially to $T_2$. 
  Then by 
  modifying $T_2$ we can assume that $V_3=0$.
  As before, by examining the reduction step to assume that $V_3=0$ and \Cref{euler-char-computation}(3), the multiplicity of 
	$\Ver_{P} W(\mu_0)$ in $V_2$ is $1$, and $V_2$ does not contain a copy of $\Ver_{P}L_{M}(\mu_1)$.
  Now suppose that $V_1$ maps non-trivially to $T_2$. First assume that $V_1 \to T_2$ is injective, 
  then we can assume that $V_1=0$. By the reasoning as above 
  we can also assume that $V_3=0$, but this is a contradiction by considering that the Euler characteristic contains $\Ver_{P}L_{M}(\mu_1)$ with negative multiplicity. 
  Thus, the other option remaining is that only the copy of $\Ver_{P}W(\mu_0)$ in $V_1$ 
  maps non-trivially, so that we can assume that $V_1=\Ver_{P}L_{M}(\mu_1)$, and $V_3=0$ as before. Computing 
  the Euler characteristic, we obtain $T_2=0$. In both cases we may 
  assume that $V_1$ is an extension of $\Ver_{P}L_{M}(\mu_1)$ and $\Ver_{P} W(\mu_0)$ in some order, 
	and  $V_2$ is a triple extension of $\Ver_{P} L_{M}(\mu_2)$,  $\Ver_{P} W(\lambda_0)$ and
  $\Ver_{P} W(\mu_0)$ in some order. 
  
  To go further we use the following lemma.
  \begin{lemma} \label{1-dim-ext}
  We have that 	 $\textnormal{Ext}^1_{(U\mathfrak{g},P)}(\Ver_{P}L_{M}(\mu_2), \Ver_{P} W(\lambda_0))$ is $1$-dimensional, generated 
  by $\Ver_{P} W(\mu_2)$, similarly  
  $\textnormal{Ext}^1_{(U\mathfrak{g},P)}(\Ver_{P}L_{M}(\mu_1), \Ver_{P} W(\mu_0))$
  is generated by $\Ver_{P} W(\mu_1)$, and \\
  $\textnormal{Ext}^1_{(U\mathfrak{g},P)}(\Ver_{P}W(\lambda_0), \Ver_{P} L_{M}(\mu_2))$ is generated by 
  $\Ver_{P} W(\nu_2)^{*}$. 
  \begin{proof}
  First we prove the statement for $\text{Ext}^1_{P}(L_{M}(\mu_2),W(\lambda_0))$. Let $W$ be such a non-split 
  extension. By adjunction $\Hom_{P}(W,W(\mu_2)) \cong \Hom_{B^{-}_{M}}(W,\mu_2)$, so there is a non-zero map 
  $W \to W(\mu_2)$. The socle of $W$ is $L_{M}(\mu_2)$, so the map must be injective, and by dimension counting it is 
  an isomorphism. Now let $N \in \textnormal{Ext}^1_{(U\mathfrak{g},P)}(\Ver_{P}L_{M}(\mu_2), \Ver_{P} W(\lambda_0))$ 
  non-trivial, as in the proof of \Cref{splitting-lemma-2} (using that $\Ver_{P} L_{M}(\mu_2)$ 
  does not contain the weight $\lambda_0$) we get subextension 
  $W \in \text{Ext}^1_{P}(L_{M}(\mu_2),W(\lambda_0))$. If $W$ is trivial we would get a splitting of $N$, so 
  $W=W(\mu_2)$, and we get a map $\Ver_{P} W(\mu_2) \to N$. Since it is an isomorphism on graded pieces by construction 
  it is an isomorphism. The other ones are proved in the same way. 
  \end{proof}
  \end{lemma}

   Now suppose that $N'_2=0$, so that by \Cref{1-dim-ext}
    $V_2$ contains $\Ver_{P} W(\nu_2)^{*}$. 
    From \cite[Prop 3.1.21]{Hodge-paper} we also have a de Rham realization functor 
    $\Psi_{G/P}: D^{b}(\mathcal{O}_{P,\Fpbar})\to D^{b}(C_{[G/P]/\Fpbar})$ computing de Rham cohomology. 
    In particular from the filtration of $\BGG_{L(\lambda_1)}$
    we get a descending filtration of $L(\lambda_2) \otimes 
    \Omega^{\bullet}_{G/P}$ which we denote by $F^{i}_{G/P}$. We use the convention that $F^{0}_{G/P}=L(\lambda_2) \otimes 
    \Omega^{\bullet}_{G/P}$ and $F^3_{G/P}$ is the last non-trivial step. 
     We will use the following result 
    on the concentration of coherent cohomology of automorphic vector bundles on $G/P$. Denote by 
    $\omega_{G/P}(\lambda)=F_{G/P}W(\lambda)$ the vector bundle on $G/P$ corresponding to $W(\lambda)$.
    \begin{lemma} \label{concentration-parabolic-flag}
    Let $\lambda \in X_{1}(T)$ generic, and $d \in \{0,1,3\}$. Then for $w \in W^{M}$ we have that 
    \begin{enumerate}
    \item $
     \H^{d}(G/P,\omega_{G/P}(w_{0,M} (ww_{0} \cdot \lambda))) \neq 0
    $
    implies $l(w)=d$, which determines $w$. 
    \item For generic $\lambda_0 \in C_0$
    and $w \in W^{M}$, 
    $\H^{\bullet}(G/P,\omega_{G/P}(w_{0,M}(ww_{0} \cdot \lambda_0)))$ is concentrated in degree $l(w)$, and it is equal 
    to $V(\lambda_0)$ as a $G$-representation. 
    \item 
    $\H^{\bullet}(G/P,\omega(w_{0,M}(w_0 \cdot \lambda)))$ is concentrated in degree $0$. 
    \end{enumerate}
    \begin{proof}
    Let $G/B \xrightarrow{\pi} G/P$. There is 
    a spectral sequence $\H^i(G/P,R^j \pi_* \LL(\lambda)) \implies \H^{i+j}(G/B,\LL(\lambda))$.
    We know that $\H^0(G/B,\LL(\lambda)) \neq 0$ implies that $\lambda$ is antidominant, i.e. $w_0 \cdot \lambda 
    \in X_{+}(T)$. This is consistent with the identity $V(\lambda)=\H^0(G/B^{-},\LL(\lambda))=\H^0(G/B,\LL(w_{0}\lambda))$.
    Also, for $\lambda \in W\cdot X_1(T)$ generic $\H^1(G/B,\LL(\lambda)) \neq 0$ implies that $\lambda \in s_{\beta}w_0 \cdot 
    X_{1}(T) \cup s_{\alpha}w_0 \cdot  X_{1}(T)$ \cite[Thm 3.6.a]{Andersen2} \cite[Thm 4.8]{AndersenH1}.
     Moreover, if $\langle \lambda +\rho, \alpha^{\vee} \rangle<0$, 
     $R^0 \pi_* \LL(\lambda)=W(w_{0,M}\lambda)$ and $R^1 \pi_* \LL(\lambda)=0$, and if 
     $\langle \lambda + \rho, \alpha^{\vee} \rangle>0$, $R^0 \pi_* \LL(\lambda)=0$ and 
     $R^1 \pi_* \LL(\lambda)=W(-w_{0,M}\lambda)^*$.
     Thus, the above spectral sequence 
     degenerates on its first page for any $\lambda$. 
     Therefore for $\lambda \in X_1(T)$, $\H^{0}(G/P,\omega_{G/P}(w_{0,M} (ww_{0} \cdot \lambda)))=\H^0(G/B,\LL(ww_{0} \cdot \lambda))$
     which proves the case $d=0$, and 
     $\H^{1}(G/P,\omega_{G/P}(w_{0,M} (ww_{0} \cdot \lambda)))=\H^1(G/B,\LL(ww_{0} \cdot \lambda))$
     which proves $d=1$. Similarly, 
     $\H^3(G/P,\omega_{G/P}(w_{0,M} (ww_{0} \cdot \lambda)))^{\vee} 
     \cong \H^{0}(G/P,\omega_{G/P}(w_{0,M} (ww_{0} \cdot \lambda))^{*})
     =\H^{1}(G/B,\LL(-ww_{0} \cdot \lambda))$, which proves $d=3$. The concentration in part $(2)$
     follows from $(1)$ and Serre duality, since in the Weyl translates of the lowest alcove all dual Weyl modules are also Weyl modules. 
     The fact that this cohomology is $V(\lambda_0)$ follows from the BGG resolution for $p$-small weights, and the concentration of each of the pieces, since 
     $\H^{\bullet}_{\dR}(G/P,L(\lambda_0))$ contains $4$ copies of $V(\lambda_0)$.
      Part $(3)$ follows 
     from Kempf vanishing, which states that $\H^{\bullet}(G/B,\LL(w_{0}\cdot \lambda))$ is concentrated 
     in degree $0$. 
      \end{proof}
    \end{lemma}
     Then \Cref{concentration-parabolic-flag} implies that 
    $\H^{0,1}(F^2_{G/P})=0$, since $\H^0(G/P,L_{M}(\mu_1+\eta)) \hookrightarrow \H^0(G/P,\omega(\mu_1+\eta))=0$.
    Also, the first non-zero term of $\text{gr}^1 F_{G/P}$ is an extension of 
    $\omega_{G/P}(\nu_0+\eta)$ and $\omega_{G/P}(\nu_2+\eta)$, so that all the terms 
    in it have trivial $\H^0$. Therefore, $\H^1(\text{gr}^1)=0$, and a fortiori $\H^1(F^1_{G/P})=0$. 
    We have that 
    $\H^{\bullet}_{\dR}(G/P,L(\lambda_2))=\H^{\bullet}_{\dR}(G/P)^{\text{dim}L(\lambda_2)}$ is concentrated 
    in even degree, 
    being isomorphic in each degree to $L(\lambda_2)$ as a $G$-representation. This is because 
    $F_{G/P}L(\lambda_2) \cong L(\lambda_2) \otimes \mathcal{O}_{G/P}$ $G$-equivariantly, 
    equipped with the trivial connection; and 
    $\H^{\bullet}_{\dR}(G/P_{\overline{\F}_p})$ has the same dimension 
    as over the complex numbers (which is concentrated in even degree) by semicontinuity and considering the Euler characteristic. 
    Therefore, from the long exact sequence for $F^1_{G/P} \to \BGG_{G/P} \to \text{gr}^0$
    we obtain the exact sequence 
    $$
    0 \to L(\lambda_2) \to \H^0(\text{gr}^{0}) \to \H^1(F^1_{G/P})=0,
    $$
    However, from looking at its shape $\H^0(\text{gr}^{0})$ contains 
    $V(\lambda_0)$, which is a contradiction. We deduce that $N_2=0$ instead. 
    Now, take $V_1$ to be an extension of $\Ver_{P}L_{M}(\mu_1)$ and $\Ver_{P} W(\mu_0)$ in some order, 
	and  $V_2$ is a triple extension of $\Ver_{P} L_{M}(\mu_2)$ by $\Ver_{P} W(\lambda_0)$, 
  and then by $\Ver_{P} W(\mu_0)$. We would like to show that the composition $V_1 \to V_2 \to \Ver_{P} W(\mu_0)$ is zero. Suppose otherwise, 
    then $\Ver_{P} W(\mu_0)$ in $V_1$ maps isomorphically to the quotient of $V_2$. Thus, we would have 
    $$
    \text{gr}^1=
    [0 \to \Ver_{P} L_{M}(\mu_1) \to V \to 0], 
    $$
    where $V$ is an extension of $\Ver_{P} L_{M}(\mu_2)$ by $\Ver_{P} W(\lambda_0)$ and $\Ver_{P} L_{M}(\mu_1)$ necessarily maps trivially to the quotient. We compute on $G/P$ again. 
    By taking the long exact sequence for $F^3_{G/P}\to \BGG \to \BGG/F^3_{G/P}$ we get a surjection 
    $\H^4(\BGG/F^3_{G/P}) \twoheadrightarrow \H^5(F^3_{G/P})=\Ker[\H^6(G/P,\omega(\lambda_1)) \to \H^6(G/P,\omega(\lambda_2))]$. By \Cref{concentration-parabolic-flag}(1)
    $\H^{\bullet}(G/P,\omega(\lambda_i))$ is concentrated in degree $6$. Therefore, by comparing its Euler characteristic to the generic fiber and using the Borel--Weil--Bott theorem 
    we see that it is equal to $V(\lambda_i)^{\vee}$ in the Grothendieck group. Thus, $\H^5(F^3_{G/P})$ contains $L(\lambda_0)$ as a Jordan--Holder factor. We will show for a contradiction that $L(\lambda_0)$ is not a Jordan--Holder factor of $\H^4(\BGG/F^3_{G/P})$. 
    By considering the long exact sequence for its filtration it is enough to prove this about $\H^4(\text{gr}^2_{G/P})$ and $\H^4(\text{gr}^1_{G/P})$, since $\H^4(\text{gr}^0_{G/P})=0$. By taking 
    the stupid filtration for the explicit shape of 
    $\text{gr}^2_{G/P}$ and its dual $\text{gr}^1_{G/P}$ and using the vanishing results of \Cref{concentration-parabolic-flag}(2),
     it is enough to show that $L(\lambda_0)$ is not a factor of $\H^2(G/P,L_{M}(\mu_2))$, $\H^3(G/P,L_M(\mu_1))$, $\H^2(G/P,L_M(\nu_1))$ or $\H^3(G/P,L_{M}(\nu_2))$. 
     \begin{itemize}
     \item Taking the long exact sequence for 
     $L_{M}(\mu_2) \to W(\mu_2) \to W(\lambda_0)$ and using that $\H^{\bullet}(G/P,\omega(\lambda_0))$ is concentrated in degree $3$ we get that $\H^2(G/P,L_{M}(\mu_2))=\H^2(G/P,\omega(\mu_2))$. By 
     \Cref{concentration-parabolic-flag}(1) $\H^{\bullet}(G/P,\omega(\mu_2))$ is concentrated in degree $2$, so arguing as before, in the Grothendieck group $\H^2(G/P,\omega(\mu_2))=V(\lambda_2)^{\vee}$, \
      which does not contain $L(\lambda_0)$ as a factor. 
     \item By Serre duality $\H^3(G/P,L_M(\mu_1))^{\vee}=\H^0(G/P,L_{M}(\nu_1)) \hookrightarrow \H^0(G/P,\omega(\nu_1))=0$, by \Cref{concentration-parabolic-flag}(1).
     \item Taking the long exact sequence for $W(\nu_0) \to W(\mu_1)^{*} \to L_{M}(\nu_1)$ and using \Cref{concentration-parabolic-flag}(2)  we get that $\H^2(G/P,L_M(\nu_1))=\H^2(G/P,\omega(\mu_1)^*)=\H^1(G/P,\omega(\mu_1))^{\vee}=0$, 
     by \Cref{concentration-parabolic-flag}(1).
     \item Taking the long exact sequence for $L_{M}(\nu_2) \to W(\nu_2) \to W(\epsilon_0)$ we get that $\H^3(G/P,L_{M}(\nu_2))=0$, since $\H^2(G/P,\omega(\epsilon_0))$ and $\H^3(G/P,\omega(\nu_2))=0$ by 
     \Cref{concentration-parabolic-flag}(1).
     \end{itemize}
 So far we have proved that $\text{gr}^2=[0 \to V_1 \xrightarrow{f} V_2 \to 0]$, where 
    $V_1$ is an extension of $\Ver_{P}L_{M}(\mu_1)$ and $\Ver_{P} W(\mu_0)$ in some order, 
	and  $V_2$ is a triple extension of $\Ver_{P} L_{M}(\mu_2)$ by $\Ver_{P} W(\lambda_0)$, 
  and then by $\Ver_{P} W(\mu_0)$.  
  Furthermore, the map $V_1 \to V_2 \to \Ver_{P} W(\mu_0)$ is trivial. 
  We have also obtained that $\Ker(f)=\mathcal{H}^1(\text{gr}^2)=\mathcal{H}^2(\BGG/F^2) \hookrightarrow \mathcal{H}^0(\text{gr}^0)$.
  To finish we proceed by ruling out the different possibilities. 
  \begin{enumerate}
  \item First assume for a contradiction that the first extension of $V_2$ is split.
   Suppose that $\Ver_{P} W(\mu_0)$ is a subobject of $V_1$. Then by Serre duality 
  $\Hom(\Ver_{P} W(\mu_0),\Ver_{P} L_{M}(\mu_2))=\Hom(\Ver_{P} L_{M}(\nu_2),\Ver_{P} W(\nu_0))$, and the latter is zero since $\nu_0 \ngeq \nu_2$. There is a unique non-trivial map 
  $\Ver_{P} W(\mu_0) \to \Ver_{P} W(\lambda_0)$, and its kernel is the image of $\Ver_{P} W(\nu_0) \to \Ver_{P} W(\mu_0)$.
   In any case we see that $\Ker(f)$ contains a vector of weight $\nu_0$, 
  which is $>$ than any vector appearing in $\mathcal{H}^0(\text{gr}^0)$. Thus, we conclude that $V_1=\Ver_{P} W(\mu_1)$.
  Since $\lambda_0 \ngeq \mu_1$, $\Ver_P L_{M}(\mu_1)$ maps trivially to $\Ver_{P} W(\lambda_0)$. 
  By the argument as above the map $\Ver_{P} L_{M}(\mu_1) \to \Ver_{P} L_{M}(\mu_2)$ must be non-trivial. Furthermore, there is only one such map, and it is injective. This is because 
  $\Ver_{P} L_{M}(\mu_1)=W(1,0)^{(p)} \otimes \Ver_{P} W(\mu_1-(p,0))$, and the non-trivial map $\Ver_{P} W(\mu_1-(p,0)) \to \Ver_{P} W(\mu_2-(p,0))$ is the first differential of the lowest alcove BGG complex 
  for the dominant translate of $\mu_2-(p,0)$. Further, the cokernel of $\Ver_{P} L_{M}(\mu_1) \to \Ver_{P} L_{M}(\mu_2)$ is $W(1,0)^{(p)} \otimes \text{Im}[\Ver_{P} W(\mu_2-(p,0))\to \Ver_{P} W(\tilde{\mu}_0)]$,
  where $\tilde{\mu}_0$ is the affine Weyl translate of $\mu_2-(p,0)$ in the same alcove as $\mu_0$. Therefore, we get a map $\Ver_{P} W(\mu_0) \to \text{Coker} \to W(1,0)^{(p)} \otimes \Ver_{P} W(\tilde{\mu}_0)$. 
  However, $\tilde{\mu}_0$ and $\mu_0$ do not have the same infinitesimal character, being in the same alcove but generically not the same weight. Hence, the map above is zero. 
  This implies that $h: V_1 \to \Ver_{P} L_{M}(\mu_2)$ factors through $\Ver_{P} L_{M}(\mu_1)$. As before we can assume 
  that the induced map $g: \Ver_{P} W(\mu_0) \to \Ver_{P} W(\lambda_0)$ is non-zero. 
  Then we can find a lift $v \in V_1$ of a vector of weight $\nu_0$ in the kernel of $g$ such that $h(v)=0$.
   This implies that $v \in \Ker(f)$, which is a 
  contradiction as before. 
  \item Now we can assume that $V_2$ contains $\Ver_{P} W(\mu_2)$. We now suppose that $V_1$ contains $\Ver_{P} W(\mu_0)$.
   We prove that the induced map 
  $\Ver_{P} W(\mu_0) \to \Ver_{P} W(\mu_2)$ vanishes, which would produce a contradiction arguing as above. By Serre duality we get a map $\Ver_{P} W(\mu_2)^* \to \Ver_{P} W(\nu_0)$. Now, although $W(\mu_2)^*$ is not 
  an irreducible $P$-representation, it is generated by the unique vector of weight $\nu_2$, since $W(\mu_2)^*$ is a non-split extension of $W(\epsilon_0)$ by $L_{M}(\nu_2)$. Since $\nu_0 \ngeq \mu_2$ such vector is mapped to zero, 
  so the map is zero. 
  \end{enumerate}
Finally, the induced map $\Ver_{P} L_{M}(\mu_1) \to \Ver_{P} L_{M}(\mu_2)$ must be the unique non-trivial one, since $\mathcal{H}^0(\text{gr}^0)$ does not contain a vector of weight $\mu_1$.

	\end{proof}
	\end{prop}
  
  We write the first case of the next lemma in detail, since the fact that the shape of the BGG complex 
  suggests the entailment 
  $F(\lambda_0) \in W(\mathfrak{m}) \implies F(\lambda_2) \in W(\mathfrak{m})$
  can be traced to this computation. 
  \begin{lemma} \label{euler-char-computation}
  Let $\lambda_0 \in C_0$, then
  \begin{enumerate}
   \item We have $n=1$ in \eqref{big-computation}.
   \item We have $[\wedge^2 \mathfrak{g}/\mathfrak{p} \otimes L(\lambda_2): W(\mu_0)]-
   [\wedge^3 \mathfrak{g}/\mathfrak{p} \otimes L(\lambda_2): W(\mu_0)]=1$ and 
   $[\wedge^2 \mathfrak{g}/\mathfrak{p} \otimes L(\lambda_2): L_{M}(\mu_1)]-
   [\wedge^3 \mathfrak{g}/\mathfrak{p} \otimes L(\lambda_2): L_{M}(\mu_1)]=1$.
   \item We have $[\mathfrak{g}/\mathfrak{p} \otimes L(\lambda_2): W(\mu_0)]-
   [L(\lambda_2): W(\mu_0)]=1$ and 
   $[\mathfrak{g}/\mathfrak{p} \otimes L(\lambda_2): L_{M}(\mu_1)]-
   [L(\lambda_2): L_{M}(\mu_1)]=0$
  \end{enumerate}
  \begin{proof}
  For $(1)$ we
  first compute the multiplicities appearing in \eqref{big-computation}. 
  	  It is a simple fact (using that $M^{\der} \cong \SL_2$) that for $V \in \Rep(P)$
    with the highest weights of its Jordan--Holder factors being Weyl translates of weights in $X_{1}(T)$,
	  $[V:W(\lambda)]-[V:L_{M}(\tilde{\lambda})]=w_{\lambda}V-w_{\lambda+\alpha}V$. Here $w_{\lambda} V$ denotes 
	  the dimension of the weight space, and $\tilde{\lambda}$ is the reflection of $\lambda$ across the hyperplane 
	  $\{\langle \mu+\rho, \alpha^{\vee} \rangle=p\}$. 
 Therefore, 
  $[L(\lambda_2):W(\lambda_0)]=w_{\lambda_0}L(\lambda_2)-w_{\lambda_i+\alpha}L(\lambda_2)$.
  This is because $L_{M}(\mu_2)=L_{M}(\tilde{\lambda_0})$  is not a
  constituent of $\wedge^k \mathfrak{g}/\mathfrak{p} \otimes L(\lambda_2)$
  unless $k=1$. Combining this with the fact that $[L(\lambda_2):W(\lambda_2)]=1$ we reduce the claim to two smaller claims 
  \begin{enumerate}
    \item $\sum^{3}_{i=0} (-1)^i[\wedge^{i} \mathfrak{g}/\mathfrak{p} \otimes L(\lambda_2):W(\lambda_1)]=-1$
  \item $w_{\lambda_0} L(\lambda_2)-w_{\lambda_0+\alpha} L(\lambda_2)=1$,
  \end{enumerate}
  where the first one is the coefficient of $w_{\lambda_0}(\lambda_1)$ in the left-hand side of 
  \eqref{big-computation}. Then $(1)$ follows from the decomposition of dual Weyl modules in \Cref{character-GSP4} and 
  the way the BGG resolution is constructed in characteristic $0$. For the second one  let 
  $\lambda'_i:=s_{\alpha} \cdot \lambda_i$, and $w_{\lambda}(\mu):=w_{\lambda} \Ver_{B}(\mu)$.
  We use the (Borel) BGG resolution, the fact that only 
  $\lambda^{'}_{1,2}$, $\lambda_{1,2}$ and $\mu_2$ out of the $\{w \cdot \lambda_i : w \in W, i=0,1,2\}$
  can be $\ge \lambda_0$, and $w_{\lambda_0}(\mu_2)=1$
  to obtain  
  \begin{equation} \label{eq11}
  w_{\mu} L(\lambda_2)=w_{\mu}(\lambda_2)-
  w_{\mu}(\lambda'_2)-w_{\mu}(\lambda_1)+
  w_{\mu}(\lambda'_1)-1 \;\;\; \textnormal{for } \mu=\lambda_0+\alpha
  \end{equation}
  and 
   $$
   w_{\mu} L(\lambda_2)=w_{\mu}(\lambda_2)-
  w_{\mu}(\lambda'_2)-w_{\mu}(\lambda_1)+w_{\mu}(\lambda'_1) \;\;\;
  \text{for } \mu=\lambda_0.
  $$ 
  Write $\lambda_0=(a,b)$, with $\alpha=(1,-1)$ and $\beta=(0,2)$.
   Then $\lambda_1=(p-b-3,p-a-3)$, $\lambda'_1=(p-a-4,p-b-2)$,
  $\lambda_2=(p+b-1,p-a-3)$, and $\lambda'_2=(p-a-4,p+b)$. We compute the dimensions of the weight spaces 
 above case by case. 
  \begin{enumerate}[label=(\alph*)]
  \item Let $\mu \le  \lambda \in X^{*}(T)$ and write $\lambda-\mu=N \alpha +M \beta$. 
  From the PBW theorem we can compute that 
  $w_{\mu}(\lambda)=\lvert \{ (x,y) \in \Z^{\ge 0} : x+2y \le N, x+y \le M\} \rvert$, where $x$ 
  corresponds to $x_{-\alpha-\beta}$ and $y$ to $x_{-2\alpha-\beta}$. Then 
  $$
  w_{\mu}(\lambda)=
  \begin{cases}
  \sum^{M}_{i=0}  \text{min}(\lfloor (N-i)/2 \rfloor,M-i) +1 \; \textnormal{ if } N \ge M \\
   \sum^{N}_{i=0} \lfloor (N-i)/2 \rfloor +1 \; \textnormal{ if } M \ge N.
  \end{cases}
  $$
  Now specify $N,M$ to $\mu=\lambda_0$ and $\lambda=\lambda_2$. We always have $N>M$. 
  Therefore, 
  \begin{align*}
    w_{\lambda_0}(\lambda_2)-w_{\lambda_0+\alpha}(\lambda_2)&=
    \sum^{M}_{i=0} \text{min}(\lfloor (N-i)/2 \rfloor,M-i)-\text{min}(\lfloor (N-i-1)/2 \rfloor,M-i) \\
    &=\lfloor N/2 \rfloor-(N-M-1)=\lfloor (p+b-a-1)/2 \rfloor-b
  \end{align*}
  since $\lfloor (N-i)/2 \rfloor \le M-i$ for $0 \le i \le 2M-N+1=p-a-b-2$.
  \item For $w_{\lambda_0}(\lambda_1)$ write $\lambda_1-\lambda_0=N(\alpha+\beta)$ where 
  $N=p-a-b-3$. We get  
  $$
    w_{\lambda_0}(\lambda_1)-w_{\lambda_0+\alpha}(\lambda_1)=
    \lfloor N/2\rfloor+1=\lfloor (p-a-b-3)/2 \rfloor+1.
  $$
  \item For $w_{\lambda_0}(\lambda'_i)$ we divide the $C_0$ alcove into $C'_{0}$ defined by 
  $2a<p-4$, and $C^{''}_{0}$ defined by $2a>p-4$. First assume that $\lambda_0 \in C'_{0}$. 
  In this case we have that $0 \le N \le M$, so that
  $$
    w_{\lambda_0}(\lambda'_1)-w_{\lambda_0+\alpha}(\lambda'_1)= 
    \lfloor N/2\rfloor+1=\lfloor (p-2a-4)/2 \rfloor+1.
  $$
  If $\lambda_0 \in C^{''}_{0}$ then all of the $N$ appearing are negative, so the 
  expression above is $0$. If $\lambda_0$ is in the boundary of both, then $N=0$ and 
  $w_{\lambda_0}(\lambda'_1)-w_{\lambda_0+\alpha}(\lambda'_1)=1$.
  \item For $\lambda \in C'_{0}$ and $\lambda=\lambda^{'}_2$
   we have $0 \le N=p-2a-4 \le M$, 
  so  that 
  $$
  w_{\lambda_0}(\lambda'_2)-w_{\lambda_0+\alpha}(\lambda'_2)=\lfloor N/2 \rfloor+1=\lfloor (p-2a-4)/2\rfloor+1.
  $$
  For $\lambda_0 \in C^{''}_{0}$ we again have $N<0$ and the expression above is zero. 
  For $\lambda_0$ in the boundary of both $N=0$ and the expression above is $1$. 
  We conclude that
  $$
    w_{\lambda_0}(\lambda'_2)-w_{\lambda_0+\alpha}(\lambda'_2)=  w_{\lambda_0}(\lambda'_1)-w_{\lambda_0+\alpha}(\lambda'_1).
  $$
  \end{enumerate}
  Plugging the above into the equations \eqref{eq11} we get
  $$
  w_{\lambda_0} L(\lambda_2)-w_{\lambda_0+\alpha} L(\lambda_2)=\lfloor (p+b-a-1)/2 \rfloor-b- \lfloor (p-a-b-3)/2 \rfloor=1
  $$
  as desired,
  by computing the two options for $a+b \bmod{2}$.
  Part $(2)$ is proved in a similar way. 
  Finally, part $(3)$ follows from $(2)$ and \Cref{character-GSP4}. 

  \end{proof}
  \end{lemma}

  \section{Cohomological consequences}
  To get some consequences about de Rham cohomology of Shimura varieties we will need the results from \cite{Hodge-paper}. In particular, let us recall 
  the de Rham realization functor. Let $(G,X)$ be a Shimura datum of Hodge type, $K$ a level hyperspecial at $p$ and $\Sh/\mathcal{O}$ its integral model. Let $P$ be a choice of Hodge parabolic, 
  with Levi $M$. In \cite{Hodge-paper} we constructed the exact functor 
  $$
  \Psi: D^{b}(\mathcal{O}_{P,\Fpbar}) \xrightarrow{f} D^{b}((\Shbar/\Fpbar)_{\crys}) \xrightarrow{g} D^{b}(C_{\Shbar/\Fpbar}),
  $$
  which extends to functors $\Psi^{\can}$ and $\Psi^{\sub}$ on the toroidal compactification, where now one considers the log-crystalline site and log differential operators.
   We record some of its properties that we will use.
  \begin{theorem} \label{properties-functor}
  \begin{enumerate}
  \item For $V \in \Rep_{\Fpbar}(G)$, $\Psi^{\can}(V[0])=\mathcal{V}^{\vee} \otimes \Omega^{\bullet}_{\Shbar^\tor}(\log D)$. 
  \item For $W \in \Rep_{\Fpbar}(P)$, $\Psi^{\can}(\Ver_{P}(W)[0])=\mathcal{W}^{\vee}[0]$.
  \item Let $\lambda, \mu \in X^*(T)$ and $\gamma \in \Phi^{+}-\Phi^{+}_{M}$ such that $\lambda \uparrow_{\gamma} \mu$. Let $\pi_* \phi(- w_{0,M}\mu-\eta \uparrow  -w_{0,M}\lambda-\eta): \Ver_{P} W(\mu+\eta)^{\vee} \to 
  \Ver_{P} W(\lambda+\eta)^{\vee}$ defined by combining \cite[Prop 5.1.3, Lem 5.1.2]{Hodge-paper}. Then by considering the above map as a two term complex
  $$
  \Psi(\pi_* \phi(- w_{0,M}\mu-\eta \uparrow  -w_{0,M}\lambda-\eta))=[\omega(\lambda+\eta) \xrightarrow{\theta_{\lambda \uparrow \mu}} \omega(\mu+\eta)],
  $$
  where $\theta_{\lambda \uparrow \mu}$ is a theta linkage map as defined in \cite[Notation 5.1.7]{Hodge-paper}.
\end{enumerate}
  \end{theorem}
   Now we state the rest of the inputs from \cite{Hodge-paper} required to prove our cohomological results. In this section 
  $\Sh$ will be one of the following. 
  \begin{enumerate}
  \item $(G,X)$ is the Shimura datum of a compact PEL unitary Shimura variety of signature $(2,1)$ with respect to an quadratic imaginary field.
  \item $(G,X)$ is the Shimura datum for the Siegel threefold. 
  \end{enumerate}

  First we state the known results on the vanishing of coherent cohomology.
  \begin{theorem} \label{concentration-coherent}
  Let $\lambda \in X^{*}(T)$. Then
  \begin{enumerate}
  \item If $\Sh$ is the Siegel threefold or a $U(2,1)$ Shimura variety with $p$ split in the reflex field, then $\H^i(\Shbar^{\tor},\omega^{\sub}(\lambda))=0$ for $i>0$ and $\lambda \in X_1(T)$ generic. 
  \item If $\mathfrak{m}$ is non-Eisenstein then $\H^i(\Shbar^{\tor},\omega^{\sub}(\lambda))_{\m}=\H^i(\Shbar^{\tor},\omega^{\can}(\lambda))_{\m}$ for all $i$.
  \item Let $w \in W$, and $\lambda \in X_1(T)$. Then $\H^0(\Shbar^\tor,\omega(w \cdot \lambda)) \neq 0$ implies that $w=1$ for generic $\lambda \in X_1(T)$.
  \item Let $w \in W^{M}$, and $\mathfrak{m}$ non-Eisenstein. Then $\H^{\bullet}(\Shbar^{\tor},\omega(w \cdot \lambda)^{\vee})_{\m}$ is concentrated in degree $d-l(w)$
  for generic $\lambda \in C_0$ in the lowest alcove, 
  where $d$ is the dimension of $\Shbar$.
  \end{enumerate}
  \begin{proof}
  The first one follows by  \cite[Ex 4.30, Thm 5.10]{alexandre} and \cite{Deding-unitary} respectively. Part (2) is \cite[Prop 5.1]{paper}.
  Part (3) is in \cite{cone-conjecture-gsp4}, and part $(4)$ follows from (2) and \cite[Prop 8.12]{Lan-Suh-non-compact}. 
  \end{proof}
  \end{theorem}
  Then we need an injectivity statement for a certain theta linkage map. We remark that in the $\GSp_4$ this was already proven in \cite{paper}.
  \begin{theorem} \label{injectivity} \cite[Thm 5.1.11]{Hodge-paper}
  Let $\lambda \in X_1(T)$ and $\mu \in X^*(T)$ such that $\lambda \uparrow_{\delta} \mu$, where $\delta \in \Phi^{+}$ is the longest root. Then the theta linkage map 
  $$
  \theta_{\lambda \uparrow \mu}: \H^0(\Shbar^{\tor},\omega(\lambda+\eta)) \to \H^0(\Shbar^{\tor},\omega(\mu+\eta))
  $$
  is injective for generic $\lambda \in X_1(T)$.
  \end{theorem}

  \subsection{The case of $\GL_3$}
	It turns out that in the case of $\GL_3$, \Cref{concentration-coherent} and \Cref{injectivity} provide just enough
  information to deduce \Cref{conjecture-concentration} from \Cref{GL3-BGG}.

\begin{theorem} \label{concentration-GL3}
Let $\Shbar$ be a compact Shimura variety for $G=GU(2,1)$ with respect to a quadratic imaginary 
number field $E$. If $p$ is inert in $E$ assume that
\Cref{concentration-coherent}(1) also holds in that case. Then for a
 generic $\lambda_0=(a,b,c) \in C_0$ 
$$
\H^{\bullet}_{\dR}(\Shbar,L(\lambda_1))
$$
is concentrated in middle degree. Let $\{\tilde{\lambda}_i,
\tilde{\mu}_i,\tilde{\nu}_i\}$ be the affine reflections of $\tilde{\lambda}_1:=-w_0 \lambda_1$, and $\eta=(1,1,-2)$.
Then $\H^2_{\dR}(\Shbar,L(\lambda_1))$ has a filtration $\H^2(F^{\bullet}_{\Shbar})$ whose
 graded pieces are 
\begin{align*}
\H^2(F^2_{\Shbar})=\textnormal{coker}&[\H^0(\Shbar,\omega(\lambda_0+\eta)) 
\xrightarrow{\theta_{\lambda_0 \uparrow \lambda_1}} \H^0(\Shbar,\omega(\lambda_1+\eta))]\\
\H^2(\textnormal{gr}^0 F_{\Shbar})=\textnormal{ker}&[\H^2(\Shbar,\omega(\tilde{\lambda}_1)^{\vee})
 \xrightarrow{\theta^{*}_{\tilde{\lambda}_0 \uparrow \tilde{\lambda}_1}} 
 \H^2(\Shbar,\omega(\tilde{\lambda}_0)^{\vee})]\\
\H^2(\textnormal{gr}^1 F_{\Shbar})=& \H^1(\Shbar,L_{M}(\tilde{\mu}_1)^{\vee}),
\end{align*}
where $\theta^{*}_{\tilde{\lambda}_0 \uparrow \tilde{\lambda}_1}$ is the dual differential operator of $\theta_{\tilde{\lambda}_0 \uparrow \tilde{\lambda}_1}$, as in \cite[Def 6.1.1]{Hodge-paper}.
\begin{proof}
Let $\BGG_{\Shbar}L(\lambda_1)$ be the complex over $\Shbar$ obtained by applying the de Rham realization functor $\Psi$ to $\BGG_{L(\lambda_1)^{\vee}}$ 
from \Cref{GL3-BGG}.
It is quasi-isomorphic to the de Rham complex $L(\lambda_1) \otimes \Omega^{\bullet}_{\Shbar}$,
and by the description of the graded pieces of \Cref{GL3-BGG} it has a descending $3$-step
 filtration $F^{\bullet}_{\Shbar}$ of the form 
 \begin{align*}
F^2_{\Shbar}&=[0 \to \omega(\lambda_0+\eta) 
\xrightarrow{\theta_{\lambda_0 \uparrow \lambda_1}} \omega(\lambda_1+\eta)] \\
\text{gr}^1 F_{\Shbar}&=[0\to L_{M}(\tilde{\mu}_1)^{\vee} \to 0] \\
\text{gr}^0 F_{\Shbar}&=[\omega(\tilde{\lambda}_1)^{\vee} 
\xrightarrow{\theta^{*}_{\tilde{\lambda}_0 \uparrow \tilde{\lambda}_1}} \omega(\tilde{\lambda}_0)^{\vee} \to 0],
 \end{align*}
 where the terms are in degrees $0,1,2$.
 Namely, to determine $F^2_{\Shbar}$ and $\text{gr}^1 F_{\Shbar}$ we use all the properties of \Cref{properties-functor}. To get the correct differential $\theta_{\lambda_0 \uparrow \lambda_1}$
 from \Cref{properties-functor}(3)
 we also use that all the $p$-restricted weights $\lambda$ are $p$-restricted for $M$, so that $W(\lambda)^{\vee}=W(-w_{0,M}\lambda)$.
 To determine $\text{gr}^0 F_{\Shbar}$ we also need to use the compatibility between Serre duality of maps of Verma modules and Serre duality of differential operators
in \cite[Prop 6.2.2]{Hodge-paper}.

 We prove that each of the graded pieces 
 has cohomology concentrated in degree $2$, which 
 implies the concentration of de Rham cohomology. For $F^2_{\Shbar}$
 we use the spectral sequence for the stupid filtration. 
 We know that for generic $\lambda_0$ the $\omega(\lambda_i+\eta)$ have cohomology concentrated 
 in $\H^0$ by \Cref{concentration-coherent}(1),
 and by \Cref{injectivity}
 we have that $\H^{0}(\omega(\lambda_0+\eta)) \to \H^0(\omega(\lambda_1+\eta))$ is injective.
 Thus, the spectral sequence degenerates on its first page, and the cohomology of $F^2_{\Shbar}$ is concentrated. 
 For $\text{gr}^0 F_{\Shbar}$, by Serre duality 
 the $\omega(\tilde{\lambda}_i)^{\vee}$ are concentrated in $\H^2$ for generic $\lambda_0$, and
 $\H^2(\Shbar,\omega(\tilde{\lambda}_1)^{\vee}) \to \H^2(\Shbar,\omega(\tilde{\lambda}_0)^{\vee})$ 
 is surjective from the compatibility  of differential 
 operators with Serre duality in \cite[Prop 6.1.2]{Hodge-paper},
 and the injectivity on $\H^0$. Finally, $\H^0(\Shbar,L_{M}(\tilde{\mu}_1)^{\vee})=0$ for generic $\lambda_0$ 
 since it embeds into 
 $\H^0(\Shbar,\omega(-w_{0,M}\tilde{\mu}_1))=0$ by
 \Cref{concentration-coherent}(3). 
 Similarly, by Serre duality $\H^2(\Shbar,L_{M}(\tilde{\mu}_1))^{\vee}=0$. 
\end{proof}
\end{theorem}

\begin{remark}
Note that in this case, if we assume that de Rham cohomology is concentrated in middle degree, then 
we would deduce that $\H^{0}(\Shbar,\omega(\lambda_0+\eta)) \to 
\H^0(\Shbar,\omega(\lambda_1+\eta))$ is injective. Otherwise we would have
$\H^1(F^2_{\Shbar}) \neq 0$, while $\H^0(\text{gr}^i F_{\Shbar})=0$ for all $i$, which would imply $\H^1_{\dR}(\Shbar,L(\lambda_1))_{\m} \neq 0$.
 Similarly, if we assume the injectivity of the theta linkage map, then concentration of 
de Rham cohomology in degree $2$ implies the concentration of $\H^{\bullet}(\Shbar,L_{M}(\tilde{\mu}_1)^{\vee})$ in degree $1$. 
\end{remark}

Given an eigenform $f \in \H^0(\Shbar,\omega(\lambda_1+\eta))$ with associated 
eigensystem $\mathfrak{m}$ we have that if $f$ does not come from $\H^0(\Shbar,\omega(\lambda_0+\eta))$
via the linkage map 
then $\H^2_{\dR}(L(\lambda_1))_{\mathfrak{m}} \neq 0$, or in other words $L(\lambda_1) \in W_{\dR}(\m)$.
However, conversely one has to keep in mind that an eigenform of $\H^2(F^2_{\Shbar})_{\mathfrak{m}}$ 
might not necessarily 
lift to an eigenform in $\H^0(\Shbar,\omega(\lambda_1+\eta))_{\mathfrak{m}}$, only to a generalized eigenform.  
\begin{example}(A non-generic weight) \label{non-generic-concentration}
Take $\lambda_0=(p-1,p-1,2) \in C_0$, which is almost at the boundary of $C_0$ and $C_1$. 
Then 
 $\H^{\bullet}(\text{gr}^1_{\Shbar} L(\lambda_1))=\H^{\bullet-1}(\Shbar,L_{M}(\tilde{\mu}_1)^{\vee})$ is still 
concentrated in degree $2$, by the more precise version of \Cref{concentration-coherent}(1) in \cite[Prop 7.3.2(1,2)]{Hodge-paper}. On the other hand, we have
$\H^1(F^2_{\Shbar} L(\lambda_1))=\Ker[\H^0(\Shbar,\omega(p,p,0))\xrightarrow{\theta_{\alpha_2}/H_{\alpha_2}}
\H^0(\Shbar,\omega(p+1,p,-1))]=\H^0(\Shbar,\omega(1,1,0))$ by \cite[Ex 5.3.6]{Hodge-paper},
so that $\H^0(\Shbar,\omega(1,1,0)) \subseteq \H^1_{\dR}(\Shbar,L(\lambda_1))$. However, 
the local Galois representations at $p$ appearing in $\H^0(\Shbar,\omega(1,1,0))$ should be quite degenerate. That is, 
in characteristic $0$ the weight (shifted by $\rho$) of the associated automorphic form would be $\{1,2,2\}$, which 
conjecturally would correspond to the Hodge--Tate weights of $\rho_{\mathfrak{m}}$. 
Therefore, it might be reasonable to expect that under some genericity conditions on $\overline{\rho}_{\mathfrak{m}}$ (i.e. 
lying in a certain open subset of the reduced Emerton-Gee stack),
$\H^{\bullet}_{\dR}(\Shbar,L(\lambda_1))_{\mathfrak{m}}$ would still be concentrated in middle degree. 
\end{example}

\subsection{The case of $\GSp_4$}
 Let $\Shbar$ be the Siegel threefold. Let $\lambda_0 \in C_0$. In this case since $w_{0}=-\text{id}$, $L(\lambda)$ is self-dual as a representation of $\text{Sp}_4$, and we will ignore 
 the twist by a central character. 
 We define the filtration $F^{i}_{\Shbar^{\tor}} \subseteq \BGG_{\Shbar^{\tor}} L(\lambda_i)$ in the derived category of $\Shbar^{\tor}$ as in the case of $\GL_3$.
  First we comment on the dual BGG complex for the $C_1$ alcove. 
\begin{remark}(The $C_1$ alcove) \label{remark-lambda1-gsp4}
   By the same method as in the case of $\GL_3$, if we assume that the map 
  $\H^0(\Shbar, \omega(\lambda_0+\eta))_{\mathfrak{m}} \to \H^0(\Shbar,\omega(\lambda_1+\eta))_{\mathfrak{m}}$ coming 
  from the differential in $F^3_{\Shbar^{\tor}} L(\lambda_1)=[0\to 0 \to \omega(\lambda_0+\eta)\to \omega(\lambda_1+\eta)]$ is injective, and 
  $\H^{2}(\Shbar^{\tor},L_{M}(\mu_1+\eta))_{\mathfrak{m}}=0$, then 
  $\H^{\bullet}_{\dR}(\Shbar^{\tor},L(\lambda_1))_{\mathfrak{m}}$ is concentrated in middle degree. 
  Here we are using the identification of canonical and subcanonical coherent cohomology localized at a non-Eisenstein 
  maximal ideal
  from \Cref{concentration-coherent}(2), so that Serre duality between the graded pieces 
  holds after localizing at $\mathfrak{m}$. 
  For Serre duality to be Hecke equivariant we need to twist the action on $\mathbb{T}$ on one of the sides 
  as in \cite[Prop 4.2.9]{higher-coleman}, but this twisted eigensystem is still non-Eisenstein. 
  Note 
  that a priori we don't know if the previous map is the linkage map $\theta_{\lambda_0 \uparrow \lambda_1}$ 
  we constructed. 
  Conversely, assume that de Rham cohomology is concentrated in middle degree. Then the map
  $\H^0(\Shbar, \omega(\lambda_0+\eta))_{\mathfrak{m}} \to \H^0(\Shbar,\omega(\lambda_1+\eta))_{\mathfrak{m}}$ is injective.
  This follows from the fact 
  that by \Cref{concentration-coherent} 
  $F^3_{\Shbar^\tor}$
  is concentrated in degrees $2,3$, $\text{gr}^{2}F_{\Shbar^\tor}$ is concentrated 
  in degree $3,4$, $\text{gr}^{1}F_{\Shbar^\tor}$ in degree $2,3$, and $\text{gr}^{0}_{\Shbar^\tor}$ in degree $3,4$. However, 
  it does not necessarily follow that $\H^{\bullet}(\Shbar^{\tor},L_{M}(\mu_1+\eta))_{\mathfrak{m}}$ is concentrated in degree $1$.
   
\end{remark}
 
  We now turn to $L(\lambda_2)$. Curiously, we will come close to proving the generic concentration in middle degree $3$, as opposed to the $C_1$ alcove. 
  We will in fact prove as much as it is necessary to upgrade the generic weak entailment of \cite{paper} to a generic entailment in \Cref{correct-entailment} ahead.
  However, as opposed to the $\GL_3$ case, first we need to prove some new results about concentration of coherent cohomology and the behaviour of a theta linkage map on higher cohomology.

  \subsubsection{A generic vanishing of coherent cohomology}
  In \Cref{random-vanishing} and \Cref{vanishing1} ahead we prove a vanishing result for coherent cohomology that we will need.
  We will use the notation $(a,b) \in X^*(T)$ from \cite[\S 2.1]{paper}.
  Let $\pi: \flag \to \Shbar^{\tor}$ be the flag Shimura variety as in 
  \cite{Hodge-paper}. For $\lambda \in X^*(T)$ there exist automorphic line bundles $\LL(\lambda)$ over $\flag$ 
  such that if $\lambda+\alpha$ is $M$-dominant, $\pi_* \LL(\lambda)=0$ and $R^1\pi_{*}\LL(\lambda)=\omega(-w_{0,M} \lambda-\alpha)^{\vee}$.
  On the other hand, if $\lambda$ is $M$-dominant then $\pi_*\LL(\lambda)=\omega(w_{0,M}\lambda)$
  and $R^1\pi_{*}\LL(\lambda)=0$. There is 
  a smooth map $\xi: \flag \to \Sbt:=[B\times B \backslash G]$ factoring through the stack of G-Zip flags \cite[Prop 4.4]{alexandre}. For $w \in W$ let $X_{w}\subseteq \flag$ be the pullback under $\xi$ of the Schubert variety 
  $[B \times B \backslash \overline{BwB}]$. We will use the fact that all the Schubert varieties for $\GSp_4$ are smooth except
   the codimension $1$ ones, namely for $w  \in \{s_{\alpha}s_{\beta}s_{\alpha}, s_{\beta}s_{\alpha}s_{\beta}\}$
   \cite[Thm 2.4]{smoothness-schubert}.
   For compactness we will use the notation $X_{\alpha}:=X_{s_{\beta}s_{\alpha}s_{\beta}}$.
   Over $\flag$ there is a full symplectic flag 
  $\LL \subseteq \omega \subseteq \LL^{\perp} \subseteq \mathcal{H}$
  and a full conjugate flag $F \subseteq \Ker(V)\subseteq F^{\perp} \subseteq \mathcal{H}$ that defines $\xi$.
   Then we have that $\LL(a,b)=\LL^{a}\otimes (\omega/\LL)^{b}$.
   We will use the following tower of strata 
  \begin{equation} \label{tower}
          X_{s_{\beta}} \subseteq  X_{s_{\alpha}s_{\beta}} \subseteq X_{\alpha} \subseteq \flag.
  \end{equation}
  First, $X_{\alpha}$ is a divisor of $\flag$ cut by a Hasse invariant which is a section of $\LL(-1,p)$, corresponding to the map
  $\LL\to \mathcal{H}/F^{\perp} \cong (\omega/\LL)^p$. Second, $ X_{s_{\alpha}s_{\beta}}$ is a divisor of $X_{\alpha}$ cut out by an element of $\LL(p-1,0)$ 
  corresponding to the map $\LL \to F^{\perp}/ \Ker(V)=\LL^p$.
   Finally, $X_{s_{\beta}}$ 
  is a divisor of $ X_{s_{\alpha}s_{\beta}}$ cut out by an element of $\LL(-p-1,0)$
  corresponding  to the map 
  $\LL \to \Ker(V)/F=\LL^{-p}$. 
We will need the following lemma in order to apply Serre duality on $X_{s_{\alpha}s_{\beta}}$.

\begin{lemma} \label{canonical-bundle}
\begin{enumerate}
\item Let $D$ be the reduced boundary divisor of $\flag$, and let $D_{s_{\alpha}s_{\beta}}:=D \cap X_{s_{\alpha}s_{\beta}}$. Then 
$D_{s_{\alpha}s_{\beta}}$ is a reduced divisor of $X_{s_{\alpha}s_{\beta}}$. 
\item 
The canonical bundle of $X_{s_{\alpha}s_{\beta}}$ is
$$
K_{X_{s_{\alpha}s_{\beta}}}=\LL(\alpha+2\beta)(-D_{s_{\alpha}s_{\beta}}).
$$
\end{enumerate}
\begin{proof}
For $(1)$ we show that the dimension of $D_{s_{\alpha}s_{\beta}}$ is at most $1$, which implies that it is a divisor. 
$X_{s_{\alpha}s_{\beta}}$ maps to the non-ordinary locus $\Shbar^{\text{n-ord}}$, with finite fibers away from the
zero-dimensional superspecial locus. Thus, it is enough to show that the intersection of $D$ and 
$\Shbar^{\text{n-ord}}$ is one-dimensional. For that we show the extension of the $p$-rank one locus 
$\Shbar^{=1,\tor} \subseteq \Shbar^{\text{n-ord}}$, whose complement has dimension $1$, does not intersect $D$ on 
an irreducible component. 
This is because the irreducible components of $D$ of positive dimension parametrize semiabelian schemes 
 $0 \to \mathbb{G}_m \to A \to E \to 0$ where $E$ is an elliptic curve. This irreducible component contains an 
 open dense subset where $E$ is ordinary, which forces $A$ to be ordinary, since Verschiebung acts invertibly on 
 the cotangent bundle of $\mathbb{G}_m$.
 Therefore, $D$ intersects $\Shbar^{=1}$ at dimension at most $1$.

 For $(2)$, first we prove the statement on the interior $Y_{s_{\alpha}s_{\beta}} \subseteq X_{s_{\alpha}s_{\beta}}$.
 Let $W$ be the kernel of $\mathfrak{g}/\mathfrak{p} \to -2\alpha-\beta$ 
 as a $B$-representation. Let $i: Y_{s_{\alpha}s_{\beta}} \hookrightarrow \flag^{\circ}$ be the map to the interior of the flag Shimura variety. 
 We claim that the map $\Sym^2 \omega \cong (\pi \circ i)^* \Omega^1_{\Shbar} \to i^* \Omega^1_{\flag^{\circ}} 
 \to \Omega^1_{Y_{s_{\alpha}s_{\beta}}}$ factors through $\Sym^2 \omega \to i^*F_{B}(W^{\vee})=\Sym^2 \omega/\LL^2$, 
 and it becomes an isomorphism. 
 Then 
  the formula for $K_{Y_{s_{\alpha}s_{\beta}}}$ follows by taking determinants. 
 It is enough to prove the claim on geometric points 
 $x \in Y_{s_{\alpha}s_{\beta}}$. Then $T_{Y_{s_{\alpha}s_{\beta}},x}$ is determined by deformations of $x$ along 
 square-zero thickenings $R \to \kappa(x)$. By Grothendieck--Messing theory, this data is the same as a choice
 of a full Hodge filtration $\tilde{\LL} \subseteq \tilde{\omega} \subseteq \tilde{\LL}^{\perp} \subseteq \mathcal{H} \otimes_{\kappa(x)} R$
 reducing to the full Hodge filtration at $x$, and such that $\tilde{\LL} \subseteq \Ker(V)$. 
 Here $\Ker(V) \subset \mathcal{H} \otimes_{\kappa(x)} R$ is the trivial lift of the conjugate filtration over $x$, 
 since the Verschiebung on $\mathcal{H} \otimes_{\kappa(x)} R$ is the base change from the one over $x$. 
 Similarly, $\tilde{F}^{\perp}=V^{-1}(\tilde{\LL}^{(p)})$ is the trivial lift of $F^{\perp}$, since 
 the augmentation ideal of $R$ is killed by Frobenius. 
  Therefore, we see that the data of 
 $T_{Y_{s_{\alpha}s_{\beta}},x}$ is equivalent to the one of the tangent bundle of the corresponding Schubert variety 
 $X(s_{\alpha}s_{\beta}) \subseteq G/B$.
 Thus, it suffices to show
 that the map $g: (\pi \circ i)^* \Omega^1_{G/P} \to i^* \Omega^1_{G/B} \to \Omega^1_{X(s_{\alpha}s_{\beta})}$ 
 factors through an isomorphism $h: i^*F_{G/B}(W^{\vee}) \cong \Omega^1_{X(s_{\alpha}s_{\beta})}$.
Here $F_{G/B}$ is the functor sending $B$-representations 
 to $G$-equivariant sheaves on $G/B$, and $i: X(s_{\alpha}s_{\beta}) \hookrightarrow G/B$ the inclusion. 
 Since $X(s_{\alpha}s_{\beta})$ is smooth,
$g$ is a map of $B$-equivariant vector bundles on $X(s_{\alpha}s_{\beta})$, so it is enough to check 
that $h$ is well-defined and an isomorphism on its fiber at $wB$ for $w \in \{1,s_{\alpha},s_{\beta},s_{\alpha}s_{\beta}\}$.
At the same time,
by considering a map $w: G/B \to G/B$
all of these are determined by its fiber at 
the identity. Then the claim follows from the classical fact that the fiber at the identity of $T_{X(s_{\alpha}s_{\beta})}$ 
contains the roots in $\Phi^{-} \cap (s_{\alpha}s_{\beta})^{-1}\Phi^{+}=\{-\alpha-\beta,-\beta\}$.

To extend it to the toroidal compactification,
let $U \subseteq X_{s_{\alpha}s_{\beta}}$ be the preimage under $\pi$ of the (extension of the) $p$-rank locus 
$\Shbar^{=1,\tor} \subseteq \Shbar^{=1}$. On the stratum of $D$ parametrizing semiabelian schemes with a torus of rank $1$ 
this is characterized by $E$ being supersingular. 
 Then $U$ is isomorphic to $\Shbar^{=1,\tor}$, since over it $\omega \cap \Ker(V)$
has rank $1$. Moreover, $D_{s_{\alpha}s_{\beta}}=D \cap U$.
This is because $D$ does not contain supersingular points, since the Ekedahl--Oort type of a point in the boundary is determined 
by the type of its abelian part. 
Therefore, it is enough 
to prove that $K_{U} \cong \LL(1,3)(-D \cap U)$. We can further reduce to the following claim. 
Let $Z \xhookrightarrow{i} Y \subseteq \Shbar^{\tor}$
be the $p$-rank $1$ locus and the union of the ordinary locus and the $p$-rank $1$ locus, respectively. 
Over $Z$, we have that $L:=\Ker(V)\cap \omega$ is a line bundle. 
We claim that $D':=D \cap Z$ is a reduced normal crossing divisor, and that the composition 
$g: \Sym^2 \omega \xrightarrow{\text{KS}} i^* \Omega^1_{Y}(\log D) \to \Omega^1_{Z}(\log D')$ factors through an isomorphism 
$\Sym^2 \omega/L^2 \cong \Omega^1_{Z}(\log D')$. This claim implies the statement about $K_{U}$ after taking determinants,
since the isomorphism 
$U \cong Z$ is given by the condition $\LL=L$. 
\'Etale locally $Y$ is given by an open subset of $\mathbb{A}^3$ with $D$ being defined by $\{x_1x_2x_3=0\}$. 
Let $D_1$ be an irreducible component of $D$, which in this model corresponds to $\{x_i=0\}$. On the other hand, 
$Z$ is defined by the Hasse invariant $H$. Then $D_i \cap Z$ is smooth since $Z$ intersects $D$ transversally 
\cite[\S 2.4]{MR4258631}, i.e.
 $dH$ is not in the span of $dx_{i}$ over $Z \cap D_i$. 
Moreover, from this model we see that the $D_i \cap Z$ intersect transversally. 
Therefore, $D'=\sum D_i \cap Z$ is a reduced normal crossing divisor. 
Since $\Ker(V)$
is parallel for the Gauss--Manin connection $\nabla$, we have that $\nabla(L) \subseteq \Ker(V) \otimes \Omega^1_{Z}(\log D')
\subseteq L^{\perp} \otimes \Omega^1_{Z}(\log D')$. By definition of the Kodaira--Spencer map $\KS$, this implies that $L^2$ is in the kernel of $g$. 
Since $Z$ intersects $D$ transversally and $D'$ is a normal crossing divisor the map 
$i^* \Omega^1_{Y}(\log D) \to \Omega^1_{Z}(\log D')$ is a surjection, this can be checked explicitly 
in the \'etale local description above. Therefore 
$\Sym^2 \omega/L^2 \to  \Omega^1_{Z}(\log D')$ is a surjection between vector bundles of the same rank, so it is an 
isomorphism.

\end{proof}
\end{lemma}
  
  We will also need a lemma on an upper bound for $W(\mathfrak{m})$, assuming that $F(\lambda_0) \in W(\mathfrak{m})$.

  \begin{lemma} \label{weight-elimination}
  Let $\lambda_0 \in C_0$, and $\m \subseteq \mathbb{T}$ a maximal ideal. Suppose that $F(\lambda_0) \in W(\m)$. Then for $\lambda_0$ generic $W(\m)$ can only contain at most another Serre weight in the lowest alcove.
  \begin{proof}
  We use the result of \cite[Thm 5.4.4]{heejong-lee} that says that for $\lambda_0$ generic $W(\m) \subseteq W^{?}(\overline{\rho}^{\text{ss}}_{\m})$. 
  \footnote{Although they use a different global setup, their proof is local, 
  based on the study of the crystalline lifts of $\overline{\rho}_{\m}$.}
  The right-hand side is Herzig's conjectural set of Serre weights \cite{Herzig-Tilouine}. One can check via its explicit 
  combinatorial recipe that $C_0 \cap W^{?}(\overline{\rho}^{\text{ss}}_{\m})$ contains at most two elements for generic $\lambda_0$. 
  \end{proof}
  \end{lemma}

    For $E$ a vector bundle over $\Shbar^{\tor}$ we use the notation $E^*=E^{\vee}\otimes \omega(\eta)$. We will use that Serre duality on $\Shbar^{\tor}$ or on any smooth $X_{w}$ is Hecke equivariant after twisting 
  the Hecke action on one of the sides \cite[Prop 4.2.9]{higher-coleman}. Then if $\m$ is non-Eisenstein, $F(\lambda_0) \in W(\mathfrak{m})$ and $\H^i(\Shbar^{\tor},E)_{\mathfrak{m}} \neq 0$, then 
  $\H^{3-i}(\Shbar^{\tor},E^*)_{\m'} \neq 0$ for a maximal ideal $\m'$ satisfying $F(-w_{0}\lambda_0) \in W(\m')$. Since $-w_0 \lambda_0$ and $\lambda_0$ only differ by a central character we will ignore 
  the difference between $\m$ and $\m'$. We also note that in our notation 
  the canonical bundle of $\flag$ is $\LL(4,2)(-D)$. 
   We now prove the vanishing of a coherent cohomology group localized at $\m$.
  The method consists of restricting 
  along the tower of strata \eqref{tower}, which will shift the weight of some of the automorphic line bundles in such a way 
  that when they are lifted back to $\flag$, their cohomology can be controlled by \Cref{concentration-coherent}.
  We also use all the dualities available, and the knowledge of the ample cone for $\LL(\lambda)$ due to \cite{alexandre}, 
  which we use to prove the vanishing for a cohomology group on $X_{s_{\beta}}$. 
  Finally, we use \Cref{weight-elimination} to rule out that the non-zero eigenclasses that might arise from this process 
  can all occur at the same time. 

  \begin{prop} \label{random-vanishing}
  Let $\lambda_0=(a,b) \in C_0$, and $\mathfrak{m}$ non-Eisenstein such that $F(\lambda_0) \in W(\mathfrak{m})$. Let $\mu_2=s_{\beta,0} \cdot \lambda_2=(p+b-1,a-p+1)$.
  Then 
  $$
\H^2(\Shbar^{\tor},L_{M}(\mu_2+\eta))_{\mathfrak{m}}=0
  $$
  for $\lambda_0 \in C_0$ generic. 
  \begin{proof}
  In all the statements we assume that we have chosen a sufficiently large constant $\epsilon$ such that $\lambda_0$ is $\epsilon$-generic. 
  First we prove the vanishing assuming that $F(p-a-5,b) \notin W(\m)$.
  Consider the exact sequence $0 \to L_{M}(\mu_2+\eta) \to W(\mu_2+\eta) \to W(\lambda_0+\eta) \to 0$. Taking the long exact sequence in cohomology after applying the functor $F_{P} : \Rep(P)\to \text{Coh}(\Shbar^{\tor})$ we 
  deduce that $\H^2(\Shbar^{\tor},L_{M}(\mu_2+\eta))_{\mathfrak{m}}=\H^2(\Shbar^{\tor},\omega(\mu_2+\eta))_{\mathfrak{m}}$, since $\H^{i}(\Shbar^{\tor},\omega(\lambda_0+\eta))=0$ for $i=1,2$ due to \Cref{concentration-coherent}(1,2).
  Then, $\H^2(\Shbar^{\tor},\omega(\mu_2+\eta))_{\mathfrak{m}}=\H^2(\flag,\LL(a-p+4,p+b+2))_{\mathfrak{m}}$. 
  We will decompose this using the long exact sequences obtained by restricting along the tower \eqref{tower}.
  Consider the exact sequence 
  $$
 \H^2(\flag,\LL(a-p+5,b+2))_{\m} \to \H^2(\flag,\LL(a-p+4,p+b+2))_{\m} \to  \H^2(X_{\alpha},\LL(a-p+4,p+b+2))_{\m}
  $$
  given by restricting to the divisor $X_{\alpha}$. The term on the left is equal to
   $\H^2(\Shbar^{\tor},\omega(b+2,a-p+5))_{\m}$. Since its weight is a Weyl translate of a weight in $C_0$, 
  if this is non-zero, then $F(p-a-5,b)\in W(\m)$. This follows by the concentration of $\H^{\bullet}(\Shbar^{\tor},\omega(b+2,a-p+5))_{\m}$ in degree $2$ 
  due to \Cref{concentration-coherent}(4), and by the argument 
  in \cite[Prop 5.1]{paper}. Concretely, since $w_{0}=-\text{id}$, $(p-b-5,b)$ is the $\rho$-dominant 
  Weyl translate of $(b+2,a-p+5)-\eta$. By our assumption on $W(\m)$ we can assume that the term on the left vanishes, 
  so we reduce to the vanishing of $\H^2(X_{\alpha},\LL(a-p+4,p+b+2))_{\m}$.

It fits in the exact sequence 
$$
\H^1(X_{s_{\alpha}s_{\beta}},\LL(a+3,p+b+2))_{\m} \to \H^2(X_{\alpha},\LL(a-p+4,p+b+2))_{\m} \to
\H^2(X_{\alpha},\LL(a+3,p+b+2))_{\m}.
$$
We claim that the term on the right vanishes, for that it is enough to show the vanishing of 
$\H^2(\flag,\LL(a+3,p+b+2))_{\m}$ and $\H^3(\flag,\LL(a+4,b+2))_{\m}=\H^2(\Shbar^{\tor},\omega(a+3,b+3))_{\m}$, 
which vanish by \Cref{concentration-coherent}(1,2). Thus, we have reduced to the vanishing of $\H^1(X_{s_{\alpha}s_{\beta}},\LL(a+3,p+b+2))_{\m}$.
 By Serre duality and \Cref{canonical-bundle} this is equivalent to the vanishing of 
 $\H^1(X_{s_{\alpha}s_{\beta}},\LL(-a-3,-p-b-2) \otimes K_{X_{s_{\alpha}s_{\beta}}})_{\m}=
 \H^1(X_{s_{\alpha}s_{\beta}},\LL(-a-2,-p-b+1)(-D_{s_{\alpha}s_{\beta}}))_{\m}$.
 Moreover, 
 $\H^1(X_{s_{\alpha}s_{\beta}},\LL(-a-2,-p-b+1)(-D_{s_{\alpha}s_{\beta}}))_{\m} \cong \H^1(X_{s_{\alpha}s_{\beta}},\LL(-a-2,-p-b+1))_{\m}$, 
 since its cone is $\H^{\bullet}(D_{s_{\alpha}s_{\beta}},\LL(a-2,-p-b+1))_{\m}$, whose eigensystems we can lift to boundary eigensystems in 
 $\Shbar^{\tor}$, which vanish by \Cref{concentration-coherent}(2).
Therefore, we have reduced to the vanishing of 
$\H^1(X_{s_{\alpha}s_{\beta}},\LL(-a-2,-p-b+1))_{\m}$.
  By considering the restriction exact sequences along 
 $X_{s_{\alpha}s_{\beta}} \subseteq X_{\alpha} \subseteq \flag$ it is enough to prove the vanishing of the following $4$
 spaces. First $\H^1(\flag,\LL(-a-2,-p-b+1))_{\m}=\H^0(\Shbar^{\tor},\omega(-a-3,-p-b+2))_{\m}$ which vanishes by \Cref{concentration-coherent}(3), second 
 $\H^2(\flag,\LL(-a-1,-2p-b+1))_{\m}=\H^2(\Shbar^{\tor},\omega(2p+b+1,a+5))^{\vee}_{\m}$, 
 third $\H^2(\flag,\LL(-a-p-1,-p-b+1))_{\m}=\H^2(\Shbar^{\tor},\omega(-p-b+1,-a-p-1))_{\m}=\H^1(\Shbar^{\tor},
 \omega(p+a+4,p+b+2))^{\vee}_{\m}$, and fourth 
 $\H^3(\flag,\LL(-a-p,-2p-b+1))_{\m}=\H^1(\Shbar^{\tor},\omega(2p+b+1,p+a+4))^{\vee}_{\m}$. 
The last three vanish by the more general version of \Cref{concentration-coherent}(1) which allows for non-$p$-restricted weights \cite[Thm 2.12]{paper}.

   Now we prove the desired vanishing 
  assuming that $F(p-a-4,b-1), F(p-a-4,b+1) \notin W(\m)$. Then the proposition follows from \Cref{weight-elimination}.
  We have that $\H^2(\Shbar^{\tor},L_{M}(\mu_2+\eta))_{\mathfrak{m}}$ is Serre dual to 
  $\H^1(\Shbar^{\tor},L_{M}(p-a-1,-b-p+1))_{\m}$. Reasoning as before the latter is equal to 
  $\H^1(\flag,\LL(-b-p+1,p-a-1))_{\m}$. It fits in the exact sequence 
  $$
  \H^1(\flag,\LL(-b-p+2,-a-1))_{\m} \to \H^1(\flag,\LL(-b-p+1,p-a-1))_{\m} \to \H^1(X_{\alpha},\LL(-b-p+1,p-a-1))_{\m},
  $$
  and the term on the left is Serre dual to 
  $\H^2(\Shbar^{\tor},\omega(p+b+1,a+4))_{\m}$, which vanishes by \Cref{concentration-coherent}(1,2).
  The term on the right fits in the exact sequence 
  $$
   \H^0(X_{s_{\alpha}s_{\beta}}, \LL(-b,p-a-1))_{\m} \to \H^1(X_{\alpha},\LL(-b-p+1,p-a-1))_{\m} \to \H^1(X_{\alpha},\LL(-b,p-a-1))_{\m}. 
  $$
  We claim that the term on the right vanishes, it is enough to show that for 
  $\H^1(\flag,\LL(-b,p-a-1))_{\m}=\H^1(\Shbar^{\tor},\omega(p-a-1,-b))_{\m}$ and 
  $\H^2(\flag,\LL(-b+1,-a-1))_{\m}=\H^1(\Shbar^{\tor},\omega(-b,-a))_{\m}$. The latter vanishes by 
  \Cref{concentration-coherent}(1,2) and Serre duality, and the non-vanishing of the first one
  implies that $F(p-a-4,b+1) \in W(\m)$ as before. 
  Therefore, we have reduced to the vanishing of $\H^0(X_{s_{\alpha}s_{\beta}}, \LL(-b,p-a-1))_{\m}$.
 It fits in the exact sequence 
 $$
\H^0(X_{s_{\alpha}s_{\beta}}, \LL(p-b+1,p-a-1))_{\m} \to \H^0(X_{s_{\alpha}s_{\beta}}, \LL(-b,p-a-1))_{\m} 
\to \H^0(X_{s_{\beta}},\LL(-b,p-a-1))_{\m},
 $$
 by considering the section of $\LL(-p-1,0)$ defining $X_{s_{\beta}}$.

  First we claim that $\H^0(X_{s_{\beta}},\LL(-b,p-a-1))_{\m}$ vanishes. 
  Let $D$ be the reduced boundary divisor of $\Shbar^{\tor}$.
 Then $D\cap X_{s_{\beta}}$ is empty. Namely, 
 $X_{s_{\beta}}$ maps with finite fibers to the supersingular locus of $\Shbar^{\tor}$,
 since a computation with the Dieudonne module of a $p$-rank $1$ abelian surface shows that $\LL=F$ cannot happen
 over $\Shbar^{=1,\tor}$.
  Then, the intersection of $D$ and the supersingular locus is empty, 
 as explained in the proof of \Cref{canonical-bundle}.
  It is enough to show that $\LL(-b,p-a-1)_{\mid X_{s_{\beta}}}$ has negative degree on $X_{s_{\beta}}$.
  Note that $X_{s_{\beta}}$ is defined by the condition $\LL=F$, and $F\cong(\omega/\LL)^{-p}$. This isomorphisms follows 
  from the Cartier isomorphism $\Ker(V)\cong \omega^{\vee,(p)}$ and the way $F$ is constructed. 
  Hence, $\LL(-b,p-a-1)_{\mid X_{s_{\beta}}}=\LL(-b-1,-a-1)_{\mid X_{s_{\beta}}}$. By \cite[Ex 4.30, Thm 5.10]{alexandre},
  $\LL(b+1,a+1)(-D)$ is ample on $\flag$. Thus, the restriction to $X_{s_{\beta}}$ of its inverse, 
  $ \LL(-b,p-a-1)_{\mid X_{s_{\beta}}}$ has negative degree. 
  Now we claim that $\H^0(X_{s_{\alpha}s_{\beta}}, \LL(p-b+1,p-a-1))_{\m}=0$. By considering the exact sequences along 
  $X_{s_{\alpha}s_{\beta}} \subseteq X_{\alpha} \subseteq \flag$ it is enough to prove the vanishing of the following $4$ groups. 
  First $\H^0(\flag,\LL(p-b+1,p-a-1))_{\m}$, which vanishes since $\pi_*\LL(p-b+1,p-a-1)=0$. Second, 
  $\H^1(\flag,\LL(p-b+2,-a-1))_{\m}=\H^0(\Shbar^{\tor},\omega(a,b-p-1)^{\vee})_{\m}$. It vanishes 
  by considering the exact sequence 
  $\omega(b-2,a-p+1)\to \omega(a,b-p-1)^{\vee} \to L_{M}(p-b+1,-a)$ and \Cref{concentration-coherent}(3).
  Third, $\H^1(\flag,\LL(-b+2,p-a-1))_{\m}=\H^1(\Shbar^{\tor},\omega(p-a-1,-b+2))_{\m}$, whose non-vanishing implies that 
  $F(p-a-4,b-1) \in W(\m)$ as before, so we can assume that it vanishes. Finally, $\H^2(\flag,\LL(-b+3,-a-1))_{\m}=
  \H^1(\Shbar^{\tor},\omega(-b+2,-a))_{\m}$, which vanishes after applying Serre duality and 
  \Cref{concentration-coherent}(1,2).

  \end{proof}
  \end{prop}

  \begin{remark}
  In fact, the proof shows that $\H^2(\Shbar^{\tor},L_{M}(\mu_2+\eta))_{\mathfrak{m}}=0$ for any non-Eisenstein $\m$.
  Otherwise, it would imply that $\m$ is modular for a lowest alcove weight, and that 
  $F(p-a-5,b)$ and either $F(p-a-4,b-1)$ or $F(p-a-4,b+1)$ are both in $W^{?}(\overline{\rho}_{\m})$, which is not possible. 
  \end{remark}
 
  In a similar way, we get a vanishing result for $\mu_1$.
   \begin{prop} \label{vanishing1}
  Let $\lambda_0=(a,b) \in C_0$, and $\mathfrak{m}$ non-Eisenstein. Let $\mu_1=s_{\beta,0} \cdot \lambda_1=(p-b-3,a-p+1)$.
  Then 
  $$
\H^2(\Shbar^{\tor},\omega(\mu_1+\eta))_{\mathfrak{m}}=0
  $$
  for $\lambda_0 \in C_0$ generic. 
  \begin{proof}
  We use the tower $X_{s_{\alpha}s_{\beta}} \subseteq X_{\alpha} \subseteq \flag$. We want to prove that $\H^2(\flag, \LL(a-p+4,p-b))_{\m}=0$. It fits in the exact sequence 
  $$
  \H^2(\flag,\LL(a-p+5,-b))_{\m} \to \H^2(\flag, \LL(a-p+4,p-b))_{\m} \to \H^2(X_{\alpha}, \LL(a-p+4,p-b))_{\m},
  $$
  and the first term is equal to $\H^2(\Shbar^{\tor},\omega(-b,a-p+5))_{\m}=\H^1(\Shbar^{\tor},\omega(p-a-2,b+3))^{\vee}_{\m}=0$ by \Cref{concentration-coherent}(2,3). We also have the exact sequence 
  $$
  \H^1(X_{s_{\alpha}s_{\beta}},\LL(a+3,p-b))_{\m} \to \H^2(X_{\alpha}, \LL(a-p+4,p-b))_{\m} \to \H^2(X_{\alpha}, \LL(a+3,p-b))_{\m}.
  $$
  We claim that the term on the right vanishes. This follows from the vanishing of $\H^2(\flag,\LL(a+3,p-b))_{\m}=\H^2(\Shbar^{\tor},\omega(p-b,a+3))_{\m}$ and 
  $\H^3(\flag,\LL(a+4,-b))_{\m}=\H^2(\Shbar^{\tor},\omega(a+3,1-b))_{\m}$. Thus we have reduced to proving the vanishing of $\H^1(X_{s_{\alpha}s_{\beta}},\LL(a+3,p-b))_{\m}$. 
  By Serre duality on this smooth surface this is equivalent to the vanishing of $\H^1(X_{s_{\alpha}s_{\beta}},\LL(-a-2,b-p+3))_{\m}$. We finish by checking that
  the following $4$ spaces that appear vanish. They all follow by \Cref{concentration-coherent} or \cite[Thm 2.12]{paper}.
  \begin{enumerate}
  \item $\H^1(\flag,\LL(-a-2,b-p+3))_{\m}=\H^0(\Shbar^{\tor},\omega(-a-3,b-p+4))_{\m}=0$.
  \item $\H^2(\flag,\LL(-a-1,b-2p+3))_{\m}=\H^2(\flag,\LL(a+5,2p-b-1))^{\vee}_{\m}=\H^2(\Shbar^{\tor},\omega(2p-b-1,a+5))^{\vee}_{\m}=0$.
  \item $\H^2(\flag,\LL(-a-p,b-p+3))_{\m}=\H^2(\Shbar^{\tor},\omega(b-p-3,-a-p))_{\m}=\H^1(\Shbar^{\tor},\omega(p+a+3,p-b+6))^{\vee}_{\m}=0$.
  \item $\H^3(\flag,\LL(-a-p+1,b-2p+3))_{\m}=\H^1(\flag,\LL(p+a+3,2p-b-1))^{\vee}_{\m}=\H^1(\Shbar^{\tor},\omega(2p-b-1,p+a+3))^{\vee}_{\m}=0$.
  \end{enumerate}
  \end{proof}

  \end{prop}

  \subsection{$p$-translating the lowest alcove BGG complex} \label{p-translation}
We explain the general idea of $p$-translating theta linkage maps and lowest alcove BGG complexes. 
As an application of the Deligne--Illusie theorem \cite{Deligne-Illusie} one gets a simple description of the Frobenius pushforward 
of a $p$-translated lowest alcove BGG complex. We use this to obtain some cohomological information of the BGG 
complex for $L(\lambda_2)$. In particular, we prove that a certain theta linkage map is injective in coherent degree $1$, 
a statement that seemed out of reach to us before. 

Let $W,V \in \Rep_{\Fpbar}(P)$. A simple observation is that as $(U\mathfrak{g},P)$-modules 
$\Ver_{P}(W^{(p)} \otimes V)=W^{(p)} \otimes \Ver_{P}V$, where $W^{(p)}$ carries a trivial $U\mathfrak{g}$ action. 
 Let $f: \Ver_{P} V^{\vee} \to \Ver_{P} W^{\vee}$ be a map, 
and $\chi \in X^*(T)$. We say that the $p$-translation of $f$ by $\chi$
is the map $\text{id} \otimes f: \Ver_{P} L_{M}(-pw_{0,M}\chi) \otimes W^{\vee} \to 
\Ver_{P} L_{M}(-pw_{0,M}\chi) \otimes V^{\vee}$, where we have used Steinberg's theorem 
$L_{M}(p\chi)=L_{M}(\chi)^{(p)}$ and the above identification.
Similarly, for $\lambda_0 \in C_0$, let $\BGG_{\Shbar^{\tor},\lambda_0}$ be the lowest alcove BGG 
complex which is quasi-isomorphic 
to the de Rham complex of $V(\lambda_0)$ \cite[Thm 8.2.3]{Hodge-paper}. By construction $\BGG_{\Shbar^{\tor},\lambda_0}=\Psi^{\can}\BGG_{\lambda_0}$, 
where $\BGG_{\lambda_0}$ is a complex in $\mathcal{O}_{P}$ whose terms are direct sums of Verma modules, as in \cite{Lan-Polo}. 
Let $\BGG_{(p\chi,\lambda_0)}$ be the complex obtained by $p$-translating by $\chi$ each map of $\BGG_{\lambda_0}$. 
Then we define the $p$-translation of  $\BGG_{\Shbar^{\tor},\lambda_0}$ by $\omega(\chi)^{(p)}$ 
as $\Psi^{\can} \BGG_{(p\chi,\lambda_0)}$.

\begin{lemma} \label{basic-p-translate}
Let $V,W \in \Rep(P)$, $\chi \in X^*(T)$, and $f: \Ver_{P} V^{\vee} \to \Ver_{P} W^{\vee}$. 
Denote $\theta=\Psi^{\can}(f): \mathcal{W} \to \mathcal{V}$.
 Let $g$ be the $p$-translation of $f$ by $\chi$ and $F$ the absolute Frobenius on $\Shbar^{\tor}$,
  then 
 $$
F_* \Psi(g): L_M(\chi) \otimes F_* \mathcal{W}=F_*(L_{M}(\chi)^{(p)} \otimes \mathcal{W})
\to F_*(L_{M}(\chi)^{(p)} \otimes \mathcal{V})=L_M(\chi) \otimes F_* \mathcal{V}
 $$
 is equal to $\text{id}_{L_{M}(\chi)} \otimes F_* \theta$.
 
 \begin{proof}
 Given $V,W$ as above one needs to check that the isomorphism $\Ver_{P} (W^{(p)} \otimes V)=W^{(p)} \otimes \Ver_{P} V$ 
 is compatible with the canonical isomorphisms 
 $e^{\le n}_{V}: F_{P} \Ver^{\le n,\vee}(V) \cong P^{n}_{\Shbar} \otimes \mathcal{V}$ from 
 \cite[Thm 3.2.2]{Hodge-paper} and the canonical isomorphism 
 $P^n_{\Shbar} \otimes [F^* \mathcal{W} \otimes \mathcal{V}]=F^*\mathcal{W} \otimes P^n_{\Shbar} \otimes \mathcal{V}$, since the first two define $\Psi$.
 The latter isomorphism is defined as follows. Let $x,y \in \mathcal{I}_{\Delta}$, and $w \in \mathcal{W}, v \in \mathcal{V}$ local sections, 
 so that $x \otimes y$ is a local section of $P^n_{\Shbar}$. 
 Then the map sends $x \otimes y \otimes (w^{(p)} \otimes v)$ to $w^{(p)} \otimes [(x \otimes y) \otimes v]$. This is well defined 
 since $x^p \otimes 1=1 \otimes x^p$ in $P^n_{\Shbar}$. After examining the construction of the 
 isomorphisms $e^{\le n}$ the statement boils down to the fact the canonical connection on $F^* \mathcal{W}$
 agrees with the one obtained by virtue of $W^{(p)}$ being a $(U\mathfrak{g},P)$-module 
 \cite[Prop 2.2.12]{Hodge-paper}. This can be done by matching the fact that the $U\mathfrak{g}$ action on $W^{(p)}$ is trivial to the fact that $w^{(p)} \in F^* \mathcal{W}$ is a horizontal section. 
 \end{proof}
\end{lemma}

We remark that the above definitions and \Cref{basic-p-translate} have obvious analogues for the Borel and $\pi: \flag \to \Shbar^{\tor}$, which are proved in the exact same way. 
In fact, we will use it to $p$-translate a complex on $\flag$, since \Cref{BGG-lambda2-gsp4}(2) dictates that we will have to consider complex made of Weyl modules for $M$,
 which cannot be reached by $p$-translating with respect to $P$. We start by proving an analogue of the lowest alcove BGG decomposition for $\flag$.

\begin{lemma} \label{BGG-flag}
Let $\lambda=(a,b) \in C_0$ satisfying $a+4<p$. There exists a complex $\BGG_{\flag,\lambda_0}$ embedded quasi-isomorphically into $V(\lambda_0) \otimes \Omega^{\bullet}_{\flag}(\log D)$, of the following form.
With cohomological grading $\BGG^i_{\flag,\lambda_0}=\oplus_{w \in W : l(w)=i} \LL(-w\cdot \lambda)$, and the components of each differential are non-trivial theta linkage maps. 
\begin{proof}
We prove that there is a complex $\BGG_{B,\lambda}$ in $\mathcal{O}_{B}$ quasi-isomorphic to $V(\lambda)$ of the form $\BGG^{4-i}_{B,\lambda}=\oplus_{w \in W: l(w)=i} \Ver_{B}(w \cdot \lambda)$, with all the maps being 
non-trivial. Then applying the analogue of $\Psi$ for $\flag$ as in \cite[Thm 8.1.3]{Hodge-paper} implies the result. For that we adapt a proof in characteristic $0$ to this setting. 
As before we define 
$
\BGG_{B,\lambda}\coloneqq \Std_{B}(V(\lambda))_{\chi_{\lambda}}
$
as the isotypic component of $\chi_{\lambda}$, and it is quasi-isomorphic to $V(\lambda)$ since the 
latter is irreducible. 
By the assumption on $\lambda$ for any $k\le 4$, the weights of $\wedge^k (\mathfrak{g}/\mathfrak{b}) \otimes V(\lambda)$
are $p$-small, so it admits an exhaustive filtration $M^{\bullet}$ by $B$-modules 
such that $\text{gr}^i M^{\bullet}=\lambda_i$ for some $p$-small weights $\lambda_i$. Thus, we get a filtration 
$\Ver_B(M^\bullet)_{\chi_{\lambda}} \subset \BGG^k_{B,\lambda}$ with graded pieces 
$\Ver_{B}(\lambda_i)_{\chi_{\lambda}}$. The latter is non-zero precisely when $\lambda_i \in W_{\text{aff}} \cdot \lambda$.
Since both weights are $p$-small this happens if and only if 
$\lambda_i=\omega \cdot \lambda$ for some $w \in W$. For $w \cdot \lambda$
to be a weight of $\wedge^k (\mathfrak{g}/\mathfrak{b}) \otimes V(\lambda)$ one needs $l(w)=k$ in which case 
its multiplicity is $1$: this can be checked over $\mathbb{C}$ where it is a classical result
\cite{Humphreys}. 
Therefore, we get a filtration of $\BGG^k_{B,\lambda}$ whose graded pieces are precisely 
$\Ver_{B}(w \cdot \lambda)$ for $w \in W(k)$, where $W(k) \subset W$ are the elemnents of length $k$. 
Now we observe that  $w \cdot \lambda \le w' \cdot \lambda$ for $w,w' \in W(k)$ implies $w=w'$.
Let $w \in W(k)$. We claim 
that we can choose two different filtrations of $\wedge^k (\mathfrak{g}/\mathfrak{b}) \otimes V(\lambda)$
as a $B$-module
such that in one of them the weight $w \cdot \lambda$ appears before all the other weights 
$\{w' \cdot \lambda :  w'\neq w \in W(k)\}$
in the filtration, 
and another one where $w \cdot \lambda$ appears last. Let $N$ be subquotient 
of the original filtration such that $w \cdot \lambda$ is a submodule of $N$ and $w' \cdot \lambda$ is a quotient.
We may assume that $N$ contains no highest weight vectors apart from possibly $w \cdot \lambda 
$ or $w' \cdot \lambda$. Otherwise such a vector would be a $B$-submodule 
of $N$ (using that for $p$-small weights $B$-representations are equivalent 
to $U\mathfrak{b}$-modules, since the elements $\frac{x^n}{n!} \in U(B)$ with $n \ge p$ automatically act trivially),
and we can move it in front 
of $w \cdot \lambda$. Thus, by induction the weights of the graded pieces of $N$, except maybe the last one, 
are in $w \cdot \lambda +\mathbb{Z}^{\le 0}\Phi^{+}$. But $w' \cdot \lambda$ is not in that set,
by assumption. 
Therefore, there are maps 
$\Ver_B(w \cdot \lambda) \hookrightarrow \BGG^k_{B,\lambda} \twoheadrightarrow \Ver_B(w \cdot \lambda)$.
 The composition is a non-zero scalar, since the weight $w \cdot \lambda$
has multiplicity one in $\BGG^k_{B,\lambda}$.
This defines $\Ver_B(w \cdot \lambda)$ as a summand of $\BGG^k_{B,\lambda}$.
Finally, any such 
two summands for different $w \in W(k)$ are disjoint, since there are no non-zero maps $\Ver_{B}(w \cdot \lambda)
\to \Ver_B(w' \cdot \lambda)$ for $w\neq w' \in W(k)$ because we cannot have 
$w \cdot \lambda \le w' \cdot \lambda$ or vice versa.

None of the two maps out of $\Ver_{B}(w_0 \cdot \lambda)$ can be zero, since the complex has cohomology concentrated in degree $4$ and maps between Borel Verma modules are injective. 
For $k \ge 1$ denote $W(k)=\{w_k,w'_k\}$.
Suppose that a map $\Ver_{B}(w_3 \cdot \lambda) \to \Ver_{B}(w_2 \cdot \lambda)$ is zero,  then $\Ver_{B}(w'_3 \cdot \lambda) \to \Ver_{B}(w_2 \cdot \lambda)$ is also zero, by the property of being a complex. 
Since $\BGG_{B,\lambda}$ is Serre self-dual up to twist this implies that the maps out of $\Ver_{P}(w'_2 \cdot \lambda)$ are zero. Since the complex in degree $2$ this would imply that the two maps to 
$\Ver_{P}(w'_2 \cdot \lambda)$ are jointly surjective, but this is not possible since the weight of the target is strictly larger than the ones on the domain. We conclude that none of the maps from degree $1$ to $2$ are zero, so that 
by Serre duality none of them are. 
\end{proof}
\end{lemma}

Given $\lambda_0=(a,b) \in C_0$ let $\lambda'_0=(p-a-1,b+3) \in C_0$. Let $M^{\bullet}$ be the $p$-translation of $\BGG_{\flag,\lambda^{'}_0}$ by $\LL(0,p)$. It has the following shape 
$$
\begin{tikzcd}[sep=small]
 &  \LL(w_{0,M}\mu_0+\eta) \arrow[r] \arrow[rdd] & \LL(w_{0,M}\lambda_0+\eta) \arrow[r] \arrow[rdd] & \LL(w_{0,M}\lambda_1+\eta) \arrow[rd]\\
\LL(w_{0,M}\mu_1+\eta) \arrow[ru] \arrow[rd] & \oplus & \oplus & \oplus & \LL(w_{0,M}\lambda_2+\eta). \\
& \LL(w_{0,M}\mu_2+\eta) \arrow[r] \arrow[ruu] & \LL(w_{0,M}\mu_4+\eta) \arrow[r] \arrow[ruu] &  \LL(w_{0,M}\lambda_4+\eta) \arrow[ru]
\end{tikzcd}
$$
\begin{figure}[h]
  \centering
  \begin{tikzpicture}[scale=1.2]
    \coordinate (A) at (0,0);
    \coordinate (B) at (2,0);
    \coordinate (C) at (2,2);
    \coordinate (D) at (4,2);
    \coordinate (P) at (1,1);
    \coordinate (Q) at (3,1);
    
    \draw (0,0)--(4,0);
    \draw[dotted] (0,0)--(0,-2);
    \draw[dotted] (0,0)--(-2,-2);
    \draw[dotted] (0,0)--(2,-2);
    \draw[dotted] (0,0)--(2,2);

    \draw (A) -- (C); 
    \draw (B) -- (Q); 
  
    \draw (B) -- (P); 
    \draw (C) -- (Q); 
    \draw (B) -- (C); 
    \draw (A) --(2,-2);
    \draw (B) --(2,-2);
    \draw (B) --(1,-1);
    \draw (B) --(3,-1);
   
    \draw (2,-2)--(3,-1);
    \draw (A)--(0,-2);

    \draw (1,-1)--(0,-2);

    \draw (3,1)--(4,0);
    \draw (3,-1)--(4,0);

    \draw (A) --(-1,-1);

    \draw (-1,-1)--(0,-2);

    \node at (0,0) {$\bullet$};
    \node at (-0.35,0) {$-\rho$};

    \node at (1.5,0.25) {$\cdot$};
    \node[left] at (1.5,0.25) {$\scriptstyle \lambda'_0$};

    \node at (2+1.5,0.25) {$\cdot$};
    \node[left] at (2+1.5,0.25) {$\scriptstyle \lambda_4$};

     \node at (0.5,0.25) {$\cdot$};
    \node[right] at (0.5,0.25) {$\scriptstyle \lambda_0$};

         \node at (0.5,-0.25) {$\cdot$};
    \node at (0.75,-0.25) {$\scriptstyle \mu_0$};

     \node at (2+1.5,-0.25) {$\cdot$};
    \node[left] at (2+1.5,-0.25) {$\scriptstyle \mu_4$};

        \node at (2+0.25,-1.5) {$\cdot$};
    \node[above] at (2+0.25,-1.5) {$\scriptstyle \mu_2$};

        \node at (2-0.25,-1.5) {$\cdot$};
    \node[above] at (2-0.25,-1.5) {$\scriptstyle \mu_1$};

      \node at (2-0.25,1.5) {$\cdot$};
    \node[below] at (2-0.25,1.5) {$\scriptstyle \lambda_1$};

     \node at (2+0.25,1.5) {$\cdot$};
    \node[below] at (2+0.25,1.5) {$\scriptstyle \lambda_2$};

     \draw[->,green] (2-0.25,-1.5) -- (2+0.25,-1.5) ;
     \draw[->,green] (2+0.25,-1.5)-- (2+1.5,-0.25)  ;
      \draw[->,green] (2+1.5,-0.25)-- (2+1.5,0.25) ;
      \draw[->,green]   (2-0.25,-1.5)--(0.5,-0.25) ;
       \draw[->,green]   (0.5,-0.25) --(0.5,0.25) ;
          \draw[->,green]   (0.5,0.25)--(2-0.25,1.5) ;
          \draw[->,green] (2-0.25,1.5)--(2+0.25,1.5);
          \draw[->,green] (2+1.5,0.25)--(2+0.25,1.5);
                    \draw[->,green,dotted] (2+0.25,-1.5) --(0.5,0.25);
                     \draw[->,green,dotted] (0.5,-0.25) --(2+1.5,-0.25);
                \draw[->,green,dotted] (2+1.5,-0.25) --(2-0.25,1.5);
             \draw[->,green,dotted] (0.5,0.25)  --(2+1.5,0.25);

\end{tikzpicture} 
  \caption{The shape of $\pi_*M^{\bullet}$ up to shift by $\eta$. The weight $\lambda'_0$ is translated to $\lambda_4$.}
  \label{figure4}
  \end{figure}

Note that all the weights are anti $M$-dominant, so that all the line bundles appearing have no $R^1\pi_{*}$. Let $F$ be the absolute Frobenius on $\flag$. By the analogue of 
\Cref{basic-p-translate} for $\flag$ we have that 
\begin{equation} \label{translated-M}
F_{*}M^{\bullet}=\LL(0,1) \otimes F_* \BGG_{\flag,\lambda'_0}.
\end{equation}
We can understand the right-hand term in the tensor product due to \cite{Illusie-coefficients}.

\begin{lemma} \label{Hodge-filtration}
Let $\lambda \in C_0$. Then 
$$
F_* \BGG_{\flag,\lambda}=\oplus \textnormal{gr}^i_{\textnormal{Hodge}}(V(\lambda) \otimes \Omega^{\bullet}_{\flag}(\log D))
$$
for generic $\lambda \in C_0$, 
where on the left-hand side we are considering the Hodge filtration induced by Griffiths transversality and the Hodge filtration on $V(\lambda)$.
\begin{proof}
Let $A \to \flag$ be the universal abelian surface away from the boundary.
 Given any $r,s \ge 1$ satisfying $r+s<p-4$, \cite[Thm 4.7]{Illusie-coefficients} proves that $F_{*}[\H^s_{\dR}(A^r/\flag) \otimes \Omega^{\bullet}_{\flag}]$ 
is the direct sum of the graded pieces of its Hodge filtration. In particular, this proves the statement over $\Shbar$ for $\lambda=(1,0)$.
For a general $\lambda$ one uses the method of cutting out $V(\lambda)$ by applying an appropriate idempotent operator to $\H^s_{\dR}(A^r/\flag)$ for some $r,s$ as in \cite[Prop 3.7]{Lan-Suh-1}. The way this is done 
shows that this idempotent also cuts out the Hodge filtration of the de Rham complex of $V(\lambda)$. To extend the result 
to the toroidal compactification, Illusie's theorem is already stated in this generality, and the rest is carried out in \cite{Lan-Suh-non-compact}.
\end{proof}
\end{lemma}

Putting together \Cref{Hodge-filtration} and \eqref{translated-M} we can get some interesting cohomological consequences about theta linkage maps. 

\begin{prop} \label{injectivity3}
There is a unique non-trivial map $\Ver_{P} W(\mu_2)^* \to \Ver_{P} W(\mu_1)^{*}$, let 
 $\theta: \omega(\mu_1+\eta) \to \omega(\mu_2+\eta)$  be its associated theta linkage map. 
 Then the composition 
$$
\H^1(\Shbar^{\tor},L_{M}(\mu_1+\eta)) \hookrightarrow \H^1(\Shbar^{\tor},\omega(\mu_1+\eta)) \xrightarrow{\theta} \H^1(\Shbar^{\tor},\omega(\mu_2+\eta))
$$
is injective for generic $\lambda_0$.
\begin{proof}
For the first statement note that $W(\mu_2)^*$ is generated as a $P$-representation by
a highest weight vector of weight $\nu_2$, since it is a non-split extension of 
 $W(\nu_0)$ by $L_{M}(\nu_2)$. The highest weight of $W(\mu_1)^*$ is $\nu_1$, and $\nu_1-\nu_2$ is a multiple of 
 the simple root $\beta$, so there a non-trivial map must be given by some power of $x_{-\beta}$. 
 On the other hand such a map exists, as the pushforward of the corresponding map of Borel Verma modules 
 \cite[Lem 5.2.1]{Hodge-paper}.
Assume first that $\H^1(M^{\bullet})=0$ for generic $\lambda_0$. Consider the page $2$ spectral sequence associated to the stupid filtration of $M^{\bullet}$, 
with the convention that the differentials are horizontal, and we label row first and column second. 
Then
$E^{(0,2)}_{2}=\Ker[\H^0(\Shbar,\omega(\lambda_0+\eta)) \to \H^0(\Shbar,\omega(\lambda_4+\eta)) \oplus \H^0(\Shbar,\omega(\lambda_1+\eta))]=0$. 
 This is because the map $\H^0(\Shbar,\omega(\lambda_0+\eta)) \to \H^0(\Shbar,\omega(\lambda_4+\eta))$
is $\theta_{\lambda_4 \uparrow \lambda_4}$, being the unique non-trivial theta linkage map between these two weights, which is injective by \Cref{injectivity} for generic $\lambda_0$. Therefore 
$E^{(1,0)}_{\infty}=E^{(1,0)}_{2}=\Ker[\H^1(\Shbar^{\tor},\omega(\mu_1+\eta)) \to \H^1(\Shbar^{\tor},\omega(\mu_0+\eta)) \oplus \H^1(\Shbar^{\tor},\omega(\mu_2+\eta))   ]$, where one of the maps is $\theta$, and the other is the canonical projection 
fitting in the exact sequence $0 \to \H^1(\Shbar^{\tor},L_{M}(\mu_1+\eta)) \hookrightarrow \H^1(\Shbar^{\tor},\omega(\mu_1+\eta)) \to \H^1(\Shbar^{\tor},\omega(\mu_0+\eta))$. This follows from the non-triviality of the differentials in 
\Cref{BGG-flag}. Thus, $E^{(1,0)}_{\infty}=\Ker[\H^1(\Shbar^{\tor},L_{M}(\mu_1+\eta)) \to \H^1(\Shbar^{\tor},\omega(\mu_2+\eta))]$.
Since $E^{(0,1)}_{2}=0$ by \Cref{concentration-coherent}(3), 
we have that $E^{(1,0)}_{\infty}=\H^1(M^{\bullet})=0$, which proves the proposition. 

Now we prove that $\H^1(M^{\bullet})=0$ for generic $\lambda_0$. Using \Cref{Hodge-filtration} we reduce to proving that 
$\H^1(\flag,\LL(0,1) \otimes \textnormal{gr}^i_{\text{Hodge}}(V(\lambda'_0) \otimes \Omega^{\bullet}_{\flag}(\log D)))=0$. 
The graded pieces of the Hodge filtration have the following form (we omit the logarithmic differentials for ease of notation), where we are taking the Hodge filtration to be decreasing: 
$$
\text{gr}^i_{\text{Hodge}}=[\text{gr}^{i+4}V(\lambda'_0) \to \text{gr}^{i+3}V(\lambda'_0) \otimes \Omega^1_{\flag} \to \ldots \to \text{gr}^{i}V(\lambda'_0) \otimes \Omega^4_{\flag}],
$$
where $\text{gr}^{*}V(\lambda'_0)$ are the graded pieces of the Hodge filtration, which come from pullback from $\Shbar^{\tor}$. This complex fits in an exact triangle 
$C^{\bullet} \to \text{gr}^i_{\text{Hodge}} \to D^{\bullet}$, induced by the exact sequence $0 \to \pi^*\Omega^1_{\Shbar} \to \Omega^1_{\flag} \to \Omega^1_{\flag/\Shbar} \to 0$. Here 
$$
C^{\bullet}=[\text{gr}^{i+4}V(\lambda'_0) \to \text{gr}^{i+3}V(\lambda'_0) \otimes \pi^* \Omega^1_{\Shbar} \to \ldots \to \text{gr}^{i+1}V(\lambda'_0) \otimes \pi^*\Omega^3_{\Shbar} \xrightarrow{0}
\text{gr}^{i}V(\lambda'_0) \otimes \Omega^4_{\flag} ],
$$
where the last differential is zero, since both the connection and $\text{gr}^{i+1}V(\lambda'_0)$ are pulled back from $\Shbar$.
Hence, $C^{\bullet}$ is quasi-isomorphic to $\pi^* \text{gr}^i_{\text{Hodge},\Shbar} \oplus 
\text{gr}^{i}V(\lambda'_0) \otimes \Omega^4_{\flag}[-4]$. By the lowest alcove BGG decomposition \cite[Thm 5.9]{Lan-Polo} the first summand is either quasi-isomorphic to zero or to $\omega(w_i \cdot \lambda'_0 +\eta)[i-3]$, where 
$w_i \in W^M$ are labelled by length. By the projection formula, $\H^{i-2}(\flag,\LL(0,1) \otimes \omega(w_i \cdot \lambda'_0 +\eta))=\H^{i-2}(\Shbar^{\tor},\omega(1,0)\otimes \omega(w_i \cdot \lambda'_0 +\eta))$. 
The latter tensor product has a filtration with graded pieces $\omega(w_i \cdot \lambda'_0 +\eta+(1,0))$ and $\omega(w_i \cdot \lambda'_0 +\eta+(0,-1))$, and both have no cohomology in degree $i-2$ for generic $\lambda_0$, 
by \Cref{concentration-coherent}(4). Therefore, $\H^1(C^{\bullet})=0$. On the other hand, $D^{\bullet}$ is of the form 
$$
D^{\bullet}=[0 \to \text{gr}^{i+3}V(\lambda'_0) \otimes \Omega^1_{\flag/\Shbar} \to \text{gr}^{i+2}V(\lambda'_0) \otimes \Omega^1_{\flag/\Shbar} \otimes \pi^*\Omega^1_{\Shbar} \to \ldots ],
$$
and we can check that the first non-zero differential is equal to the first non-zero differential of $\pi^* \text{gr}^i_{\text{Hodge},\Shbar}$ tensored with $\Omega^1_{\flag/\Shbar}$. By the lowest alcove BGG decomposition again,
 this map is either injective or it has a kernel isomorphic to $\omega(w_3 \cdot \lambda'_0+\eta)$. Then $\H^0(\flag,\LL(0,1) \otimes \omega(w_3 \cdot \lambda'_0+\eta) \otimes \Omega^1_{\flag/\Shbar})=0$
since $\LL(0,1) \otimes \Omega^1_{\flag/\Shbar}=\LL(1,0)$ has no $\pi_*$. By considering the spectral sequence for the stupid filtration of $D^{\bullet}$ we conclude that $\H^1(D^{\bullet})=0$, so that 
$\H^1(M^{\bullet})=0$.

\end{proof}
\end{prop}

  \subsubsection{The generic entailment}
    We can now prove the main application of the paper. 
	\begin{theorem} \label{correct-entailment}
	Let $\Shbar$ be the Siegel threefold, and $\lambda_0 \in C_0$. Let 
  $\mathfrak{m}$ be non-Eisenstein. 
 Then 
	$$
  F(\lambda_0) \in W(\m) \implies F(\lambda_2) \in W(\m)
  $$
  holds for generic $\lambda_0 \in C_0$. Moreover, \Cref{etale-dR-SW}(2) holds restricting to $\m$ non-Eisenstein and Serre weights in $C_2$. 
  \begin{proof}
    First, we note that $\H^{\bullet}_{\dR}(\Shbar^{\tor},L(\lambda_0))_{\m}$  (where $\m$ is non-Eisenstein) is concentrated in degree $3$ by \cite{Lan-Suh-non-compact}, 
    and by the integral \'etale-crystalline comparison theorem so is 
    the corresponding \'etale cohomology. Therefore $F(\lambda_0) \in W(\m) \iff \H^0(\Shbar,\omega(\lambda_0+\eta))_{\m} \neq 0$ as in \cite[Prop 5.1]{paper}, even if we are not assuming that $\m$ is generic.
   We start with the assumption that $F(\lambda_0) \in W(\m)$.
   Given $\lambda \in X_1(T)$ let $\chi_{\dR}L(\lambda)$ be the Euler characteristic of $\H^{\bullet}_{\dR}(\Shbar^{\tor},L(\lambda))_{\m}$ with middle degree counted positively. We define $\chi_{\et}F(\lambda)$ analogously for 
   \'etale cohomology. By the proof of \Cref{concentration-implies-dR-et}(1) we have that $\chi_{\dR} V(\lambda)=\chi_{\et} V(\lambda)_{\Fpbar}$. The integral \'etale-crystalline comparison theorem gives us 
   $\chi_{\dR} V(\lambda_0)=\chi_{\et}F(\lambda_0)$ for generic $\lambda_0$. Then we use the short exact sequences $0 \to L(\lambda_1) \to V(\lambda_1) \to V(\lambda_0) \to 0$ and $0 \to L(\lambda_2) \to V(\lambda_2) \to L(\lambda_1) \to 0$
   to conclude that $\chi_{\dR} L(\lambda_2)=\chi_{\et}F(\lambda_2)$ for generic $\lambda_0$. In particular, to prove the entailment it is enough for us to prove that $\chi_{\dR} L(\lambda_2)>0$. 

   By applying \Cref{properties-functor} to \Cref{BGG-lambda2-gsp4} we get a descending $4$-step filtration $F^{\bullet} $ of the logarithmic de Rham complex of 
   $\underline{L}(\lambda_2)$ over $\Shbar^{\tor}$ of the following form (where we use 
   all the properties of the theorems above). First 
   $F^3$ sits 
  in an exact triangle
   \begin{equation} \label{hey}
   [0 \to 0 \to \omega(\lambda_1+\eta) \xrightarrow{\theta_{\lambda_1 \uparrow \lambda_2}}  \omega(\lambda_2+\eta)] \to F^3 \to \omega(\lambda_0+\eta)[-3].
   \end{equation}
   Second, $\text{gr}^2$ sits in an exact triangle 
   \begin{equation} \label{hey2}
   C^{\bullet}:=[0 \to \omega(\mu_1+\eta) \xrightarrow{\theta} \omega(\mu_2+\eta) \to 0 ] \to \textnormal{gr}^2 \to \omega(\mu_0+\eta)[-2],
   \end{equation}
   where $\theta$  is the map from \Cref{injectivity3}.
   Finally, $\text{gr}^i$ for $i=0,1$ are Serre dual to $\text{gr}^{3-i}$, so that it is enough to show that $\chi \H^{\bullet}(F^3)_{\m} + \chi \H^{\bullet}(\text{gr}^2)_{\m}>0$. First, $\H^{\bullet}(F^3)_{\m}$ is concentrated in degree $3$. 
   The left-hand term in \eqref{hey} has cohomology concentrated in degree $3$ by considering the spectral sequence from the stupid filtration, 
   and using \Cref{concentration-coherent}(1,2) and \Cref{injectivity}. The right-hand term has cohomology concentrated in degree $3$
   by \Cref{concentration-coherent}(1,2), and by assumption $\H^0(\Shbar^{\tor},\omega(\lambda_0+\eta))_{\m} \neq 0$. Thus, $\chi \H^{\bullet}(F^3)_{\m}>0$. We now prove that $\chi \H^{\bullet}(\text{gr}^2)_{\m} \ge 0$, which concludes the proof 
   of the entailment. 

   Note that $\H^2(\Shbar^{\tor},L_M(\mu_2+\eta))_{\m}=\H^2(\Shbar^{\tor},\omega(\mu_2+\eta))_{\m}$ since 
   $\H^{i}(\Shbar^{\tor},\omega(\lambda_0+\eta))_{\m}=0$ for $i=1,2$. 
   Therefore, by \Cref{vanishing1}, \Cref{random-vanishing} and 
   \Cref{concentration-coherent} we have that $\H^{\bullet}(\Shbar^{\tor},\omega(\mu_i+\eta))_{\m}$ are concentrated in degree $1$ 
   for $i=1,2$. Thus,
   $\H^{\bullet}(C^{\bullet})_{\m}$  is concentrated in degree $2$ and $3$, and $\H^{2}(\text{gr}^2)_{\m}=\H^2(C^{\bullet})_{\m}=
   \Ker[\theta: \H^{1}(\Shbar^{\tor},\omega(\mu_1+\eta))_{\m} \to \H^{1}(\Shbar^{\tor},\omega(\mu_2+\eta))_{\m}]$.
    By \Cref{concentration-coherent}(4) there is an exact sequence 
   $$ 
   0 \to \H^1(\Shbar^{\tor},L_{M}(\mu_1+\eta))_{\m} \to \H^1(\Shbar^{\tor},\omega(\mu_1+\eta))_{\m} \to \H^1(\Shbar^{\tor},\omega(\mu_0+\eta))_{\m}.
   $$
   Then, by \Cref{injectivity3} we conclude that 
   $
   \H^{2}(\text{gr}^2)_{\m}=\Ker(\theta) \hookrightarrow \H^1(\Shbar^{\tor},\omega(\mu_0+\eta))_{\m}
   $
   via the map above. This means that the only possible negative summand in $\chi \H^{\bullet}(\text{gr}^2)_{\m}$ 
   gets cancelled by the right-hand term in \eqref{hey2}, so that $\chi \H^{\bullet}(\text{gr}^2)_{\m} \ge 0$. 
   Finally, for the statement about \Cref{etale-dR-SW}(2) we have to prove that $F(\lambda_2) \in W(\m) \iff L(\lambda_2 ) \in W_{\dR}(\m)$ for $\m$ non-Eisenstein. Here we have proved that $\H^{\bullet}_{\dR}(\Shbar^{\tor},L(\lambda_2))_{\m}$
   is concentrated in middle degree except for possibly $\H^{2}_{\dR}(\Shbar^{\tor},L(\lambda_2))_{\m} \subseteq \H^1(\Shbar^{\tor},\omega(\mu_0+\eta))_{\m}$ or 
   $\H^2(\Shbar^{\tor},\omega(\nu_0+\eta))_{\m} \twoheadrightarrow \H^{4}_{\dR}(\Shbar^{\tor},L(\lambda_2))_{\m}$. If $\H^1(\Shbar^{\tor},\omega(\mu_0+\eta))_{\m}=0=\H^2(\Shbar^{\tor},\omega(\nu_0+\eta))^{\vee}_{\m}$
   then we conclude by \Cref{concentration-implies-dR-et}, otherwise we would have $F(\lambda_0) \in W(\m)$ by \cite[Prop 5.1]{paper}, 
   and the result follows from the above. 

  \end{proof}
	\end{theorem}
   
  We remark that all the properties 
  in \Cref{BGG-lambda2-gsp4} are used in the proof of \Cref{correct-entailment}.

  \begin{remark}
  	One could try to say something about a filtration on $\BGG_{L(\lambda_3)}$, but one can see from \Cref{BGG-lambda2-gsp4} 
  that it becomes non-obvious how to define a meaningful filtration, and  
  the computations get increasingly hard. One plausible option is that $F^3 \BGG_{\Shbar} L(\lambda_3)
  =[0 \to \omega(\lambda_1+\eta) \to \omega(\lambda_2+\eta) \to \omega(\lambda_3+\eta)]$,
  which would be a piece of the $p$-translation by $\omega(p,p)$ of a lowest alcove dual BGG complex. Then the method in \Cref{p-translation}
  proves the generic concentration in degree $3$ of this piece. 
  In particular, it implies that $\Ker[\H^0(\Shbar, \omega(\lambda_2+\eta)) \to \H^0(\Shbar,\omega(\lambda_3+\eta))]=
  \text{Im}[\H^0(\Shbar, \omega(\lambda_1+\eta)) \to \H^0(\Shbar,\omega(\lambda_2+\eta))]$.
  \end{remark}

\bibliographystyle{alpha}
\bibliography{bib}

@article {Lan-Suh-1,
    AUTHOR = {Lan, Kai-Wen and Suh, Junecue},
     TITLE = {Vanishing theorems for torsion automorphic sheaves on compact
              {PEL}-type {S}himura varieties},
   JOURNAL = {Duke Math. J.},
  FJOURNAL = {Duke Mathematical Journal},
    VOLUME = {161},
      YEAR = {2012},
    NUMBER = {6},
     PAGES = {1113--1170},
      ISSN = {0012-7094},
   MRCLASS = {11G18 (11F75 14F17 14F30 14G35)},
  MRNUMBER = {2913102},
MRREVIEWER = {Yifeng Liu},
       DOI = {10.1215/00127094-1548452},
       URL = {https://doi.org/10.1215/00127094-1548452},
}

@incollection {Tilouine-Polo,
    AUTHOR = {Polo, Patrick and Tilouine, Jacques},
     TITLE = {Bernstein-{G}elfand-{G}elfand complexes and cohomology of
              nilpotent groups over {$\Bbb Z_{(p)}$} for representations
              with {$p$}-small weights},
      NOTE = {Cohomology of Siegel varieties},
   JOURNAL = {Ast\'{e}risque},
  FJOURNAL = {Ast\'{e}risque},
    NUMBER = {280},
      YEAR = {2002},
     PAGES = {97--135},
      ISSN = {0303-1179},
   MRCLASS = {17B56 (17B10 17B50 20G30 20J06)},
  MRNUMBER = {1944175},
MRREVIEWER = {James E. Humphreys},
}

@article {Lan-Polo,
    AUTHOR = {Lan, Kai-Wen and Polo, Patrick},
     TITLE = {Dual {BGG} complexes for automorphic bundles},
   JOURNAL = {Math. Res. Lett.},
  FJOURNAL = {Mathematical Research Letters},
    VOLUME = {25},
      YEAR = {2018},
    NUMBER = {1},
     PAGES = {85--141},
      ISSN = {1073-2780},
   MRCLASS = {11G18 (11F55 11G15 17B50 20G30)},
  MRNUMBER = {3818616},
MRREVIEWER = {Peter Bruin},
       DOI = {10.4310/mrl.2018.v25.n1.a5},
       URL = {https://doi.org/10.4310/mrl.2018.v25.n1.a5},
}

@inproceedings {katz-theta,
    AUTHOR = {Katz, Nicholas M.},
     TITLE = {A result on modular forms in characteristic {$p$}},
 BOOKTITLE = {Modular functions of one variable, {V} ({P}roc. {S}econd
              {I}nternat. {C}onf., {U}niv. {B}onn, {B}onn, 1976)},
    SERIES = {Lecture Notes in Math., Vol. 601},
     PAGES = {53--61},
 PUBLISHER = {Springer, Berlin-New York},
      YEAR = {1977},
   MRCLASS = {14D20 (10D05 14K15)},
  MRNUMBER = {463169},
MRREVIEWER = {V. V. Shokurov},
}

@misc{Hamann-Lee,
      title={Torsion Vanishing for Some {S}himura Varieties}, 
      author={Linus Hamann and Si Ying Lee},
      year={2023},
      eprint={2309.08705},
      archivePrefix={arXiv},
      primaryClass={math.NT}
}

@article {Gee-Herzig-Savitt,
    AUTHOR = {Gee, Toby and Herzig, Florian and Savitt, David},
     TITLE = {General {S}erre weight conjectures},
   JOURNAL = {J. Eur. Math. Soc. (JEMS)},
  FJOURNAL = {Journal of the European Mathematical Society (JEMS)},
    VOLUME = {20},
      YEAR = {2018},
    NUMBER = {12},
     PAGES = {2859--2949},
      ISSN = {1435-9855},
   MRCLASS = {11F80 (11F75)},
  MRNUMBER = {3871496},
MRREVIEWER = {Nguy\cftil{e}n Qu\^{o}c Th\'{a}ng},
       DOI = {10.4171/JEMS/826},
       URL = {https://doi.org/10.4171/JEMS/826},
}

@article {Herzig-Tilouine,
    AUTHOR = {Herzig, Florian and Tilouine, Jacques},
     TITLE = {Conjecture de type de {S}erre et formes compagnons pour {$\rm
              GSp_4$}},
   JOURNAL = {J. Reine Angew. Math.},
  FJOURNAL = {Journal f\"{u}r die Reine und Angewandte Mathematik. [Crelle's
              Journal]},
    VOLUME = {676},
      YEAR = {2013},
     PAGES = {1--32},
      ISSN = {0075-4102},
   MRCLASS = {11F70 (11F80)},
  MRNUMBER = {3028753},
MRREVIEWER = {Anne-Marie H. Aubert},
       DOI = {10.1515/CRELLE.2011.190},
       URL = {https://doi.org/10.1515/CRELLE.2011.190},
}

@misc{heejong-lee,
      title={Emerton--{G}ee stacks, {S}erre weights, and {B}reuil--{M}\'ezard conjectures for $\mathrm{GSp}_4$}, 
      author={Heejong Lee},
      year={2023},
      eprint={2304.13879},
      archivePrefix={arXiv},
      primaryClass={math.NT}
}

@article {BGG,
    AUTHOR = {Bernstein, I. N. and Gelfand, I. M. and Gelfand, S. I.},
     TITLE = {Structure of representations that are generated by vectors of
              highest weight},
   JOURNAL = {Funkcional. Anal. i Prilo\v{z}en.},
  FJOURNAL = {Akademija Nauk SSSR. Funkcional\cprime nyi Analiz i ego Prilo\v{z}enija},
    VOLUME = {5},
      YEAR = {1971},
    NUMBER = {1},
     PAGES = {1--9},
      ISSN = {0374-1990},
   MRCLASS = {16A64 (20G05)},
  MRNUMBER = {291204},
MRREVIEWER = {A. J. Coleman},
}

@article {cone-conjecture-gsp4,
    AUTHOR = {Goldring, Wushi and Koskivirta, Jean-Stefan},
     TITLE = {Automorphic vector bundles with global sections on ${G}$-{Z}ip-schemes},
   JOURNAL = {Compos. Math.},
  FJOURNAL = {Compositio Mathematica},
    VOLUME = {154},
      YEAR = {2018},
    NUMBER = {12},
     PAGES = {2586--2605},
      ISSN = {0010-437X},
   MRCLASS = {11F33 (11F46 14D07 14G35 14L30)},
  MRNUMBER = {3870455},
MRREVIEWER = {Rolf Berndt},
       DOI = {10.1112/S0010437X18007467},
       URL = {https://doi.org/10.1112/S0010437X18007467},
}

@article {alexandre,
    AUTHOR = {Alexandre, Thibault},
     TITLE = {Vanishing results for the coherent cohomology of automorphic
              vector bundles over the {S}iegel variety in positive
              characteristic},
   JOURNAL = {Algebra Number Theory},
  FJOURNAL = {Algebra \& Number Theory},
    VOLUME = {19},
      YEAR = {2025},
    NUMBER = {1},
     PAGES = {143--193},
      ISSN = {1937-0652},
   MRCLASS = {14G35},
  MRNUMBER = {4836460},
       DOI = {10.2140/ant.2025.19.143},
       URL = {https://doi.org/10.2140/ant.2025.19.143},
}

@misc{stacks-project,
    shorthand    = {Stacks},
    author       = {The {Stacks Project Authors}},
    title        = {\textit{Stacks Project}},
    howpublished = {\url{https://stacks.math.columbia.edu}},
    year         = {2018},
  }

@book {Janzten-book,
    AUTHOR = {Jantzen, Jens Carsten},
     TITLE = {Representations of algebraic groups},
    SERIES = {Mathematical Surveys and Monographs},
    VOLUME = {107},
   EDITION = {Second},
 PUBLISHER = {American Mathematical Society, Providence, RI},
      YEAR = {2003},
     PAGES = {xiv+576},
      ISBN = {0-8218-3527-0},
   MRCLASS = {20G05 (17B10)},
  MRNUMBER = {2015057},
}

@book {Humphreys,
    AUTHOR = {Humphreys, James E.},
     TITLE = {Representations of semisimple {L}ie algebras in the {BGG}
              category {${O}$}},
    SERIES = {Graduate Studies in Mathematics},
    VOLUME = {94},
 PUBLISHER = {American Mathematical Society, Providence, RI},
      YEAR = {2008},
     PAGES = {xvi+289},
      ISBN = {978-0-8218-4678-0},
   MRCLASS = {17B10},
  MRNUMBER = {2428237},
MRREVIEWER = {Serge M. Skryabin},
       DOI = {10.1090/gsm/094},
       URL = {https://doi.org/10.1090/gsm/094},
}

@article {Jantzen-decomp,
    AUTHOR = {Jantzen, Jens C.},
     TITLE = {Darstellungen halbeinfacher {G}ruppen und kontravariante
              {F}ormen},
   JOURNAL = {J. Reine Angew. Math.},
  FJOURNAL = {Journal f\"{u}r die Reine und Angewandte Mathematik. [Crelle's
              Journal]},
    VOLUME = {290},
      YEAR = {1977},
     PAGES = {117--141},
      ISSN = {0075-4102},
   MRCLASS = {20G15},
  MRNUMBER = {432775},
MRREVIEWER = {James E. Humphreys},
       DOI = {10.1515/crll.1977.290.117},
       URL = {https://doi.org/10.1515/crll.1977.290.117},
}

@article {Lan-Suh-non-compact,
    AUTHOR = {Lan, Kai-Wen and Suh, Junecue},
     TITLE = {Vanishing theorems for torsion automorphic sheaves on general
              {PEL}-type {S}himura varieties},
   JOURNAL = {Adv. Math.},
  FJOURNAL = {Advances in Mathematics},
    VOLUME = {242},
      YEAR = {2013},
     PAGES = {228--286},
      ISSN = {0001-8708},
   MRCLASS = {11G18 (14F17 14F20 14F30)},
  MRNUMBER = {3055995},
MRREVIEWER = {Liang Xiao},
       DOI = {10.1016/j.aim.2013.04.004},
       URL = {https://doi.org/10.1016/j.aim.2013.04.004},
}

@article{Le-Hung-Lin,
author={Le Hung, Bao V. and Zhongyipan Lin},
title={{B}reuil-{M}\'ezard cycles for low rank groups. {I}n preparation.},
year={}
}

@article {Illusie-coefficients,
    AUTHOR = {Illusie, Luc},
     TITLE = {R\'{e}duction semi-stable et d\'{e}composition de complexes de de
              {R}ham \`a coefficients},
   JOURNAL = {Duke Math. J.},
  FJOURNAL = {Duke Mathematical Journal},
    VOLUME = {60},
      YEAR = {1990},
    NUMBER = {1},
     PAGES = {139--185},
      ISSN = {0012-7094},
   MRCLASS = {14F40 (14F30)},
  MRNUMBER = {1047120},
MRREVIEWER = {Michel Gros},
       DOI = {10.1215/S0012-7094-90-06005-3},
       URL = {https://doi.org/10.1215/S0012-7094-90-06005-3},
}

@article{paper,
author={Ortiz, Martin},
title={A generic entailment for {GS}p$_4$},
year={2024},
}

@misc{Deding-unitary,
      title={Positivity of automorphic vector bundles on unitary {S}himura varieties}, 
      author={Deding Yang},
      year={2025},
      eprint={2503.08119},
      archivePrefix={arXiv},
      primaryClass={math.NT},
      url={https://arxiv.org/abs/2503.08119}, 
}

@article {general-cone-conjecture0,
    AUTHOR = {Goldring, Wushi and Koskivirta, Jean-Stefan},
     TITLE = {Automorphic vector bundles with global sections on {$G$}-{
              {Z}ip}{$^\mathcal{Z}$}-schemes},
   JOURNAL = {Compos. Math.},
  FJOURNAL = {Compositio Mathematica},
    VOLUME = {154},
      YEAR = {2018},
    NUMBER = {12},
     PAGES = {2586--2605},
      ISSN = {0010-437X},
   MRCLASS = {11F33 (11F46 14D07 14G35 14L30)},
  MRNUMBER = {3870455},
MRREVIEWER = {Rolf Berndt},
       DOI = {10.1112/S0010437X18007467},
       URL = {https://doi.org/10.1112/S0010437X18007467},
}

@article {Kisin-integral-model,
    AUTHOR = {Kisin, Mark},
     TITLE = {Integral models for {S}himura varieties of abelian type},
   JOURNAL = {J. Amer. Math. Soc.},
  FJOURNAL = {Journal of the American Mathematical Society},
    VOLUME = {23},
      YEAR = {2010},
    NUMBER = {4},
     PAGES = {967--1012},
      ISSN = {0894-0347},
   MRCLASS = {11G18 (14G35)},
  MRNUMBER = {2669706},
MRREVIEWER = {Jeffrey D. Achter},
       DOI = {10.1090/S0894-0347-10-00667-3},
       URL = {https://doi.org/10.1090/S0894-0347-10-00667-3},
}

@article {AndersenH1,
    AUTHOR = {Andersen, Henning Haahr},
     TITLE = {Representation theory via cohomology of line bundles},
   JOURNAL = {Transform. Groups},
  FJOURNAL = {Transformation Groups},
    VOLUME = {28},
      YEAR = {2023},
    NUMBER = {3},
     PAGES = {1033--1058},
      ISSN = {1083-4362},
   MRCLASS = {20G05 (20G10)},
  MRNUMBER = {4633003},
MRREVIEWER = {Anton Cox},
       DOI = {10.1007/s00031-022-09769-x},
       URL = {https://doi.org/10.1007/s00031-022-09769-x},
}

@article {Andersen2,
    AUTHOR = {Andersen, Henning Haahr},
     TITLE = {The first cohomology group of a line bundle on {$G/B$}},
   JOURNAL = {Invent. Math.},
  FJOURNAL = {Inventiones Mathematicae},
    VOLUME = {51},
      YEAR = {1979},
    NUMBER = {3},
     PAGES = {287--296},
      ISSN = {0020-9910},
   MRCLASS = {14M15 (14F05 14L30 20G05)},
  MRNUMBER = {530635},
MRREVIEWER = {S. I. Gel\cprime fand},
       DOI = {10.1007/BF01389921},
       URL = {https://doi.org/10.1007/BF01389921},
}

@article {caraiani-scholze-non-compact,
    AUTHOR = {Caraiani, Ana and Scholze, Peter},
     TITLE = {On the generic part of the cohomology of non-compact unitary
              {S}himura varieties},
   JOURNAL = {Ann. of Math. (2)},
  FJOURNAL = {Annals of Mathematics. Second Series},
    VOLUME = {199},
      YEAR = {2024},
    NUMBER = {2},
     PAGES = {483--590},
      ISSN = {0003-486X},
   MRCLASS = {11R39 (14G35 14G45)},
  MRNUMBER = {4713019},
       DOI = {10.4007/annals.2024.199.2.1},
       URL = {https://doi.org/10.4007/annals.2024.199.2.1},
}

@article {caraiani-scholze-compact,
    AUTHOR = {Caraiani, Ana and Scholze, Peter},
     TITLE = {On the generic part of the cohomology of compact unitary
              {S}himura varieties},
   JOURNAL = {Ann. of Math. (2)},
  FJOURNAL = {Annals of Mathematics. Second Series},
    VOLUME = {186},
      YEAR = {2017},
    NUMBER = {3},
     PAGES = {649--766},
      ISSN = {0003-486X},
   MRCLASS = {11F75 (11G18 11R23 14G35)},
  MRNUMBER = {3702677},
MRREVIEWER = {Nguy\cftil{e}n Qu\^{o}c Th\'{a}ng},
       DOI = {10.4007/annals.2017.186.3.1},
       URL = {https://doi.org/10.4007/annals.2017.186.3.1},
}

@article {Diamond-Sasaki,
    AUTHOR = {Diamond, Fred and Sasaki, Shu},
     TITLE = {A {S}erre weight conjecture for geometric {H}ilbert modular
              forms in characteristic {$p$}},
   JOURNAL = {J. Eur. Math. Soc. (JEMS)},
  FJOURNAL = {Journal of the European Mathematical Society (JEMS)},
    VOLUME = {25},
      YEAR = {2023},
    NUMBER = {9},
     PAGES = {3453--3536},
      ISSN = {1435-9855},
   MRCLASS = {11F33 (11F41 11F80)},
  MRNUMBER = {4634676},
MRREVIEWER = {Andrew James Graham},
       DOI = {10.4171/jems/1265},
       URL = {https://doi.org/10.4171/jems/1265},
}

@misc{prismatic-realization,
      title={The prismatic realization functor for {S}himura varieties of abelian type}, 
      author={Naoki Imai and Hiroki Kato and Alex Youcis},
      year={2025},
      eprint={2310.08472},
      archivePrefix={arXiv},
      primaryClass={math.NT},
      url={https://arxiv.org/abs/2310.08472}, 
}

@misc{F-gauges,
title={Prismatic {F}-gauges},
author={Bhargav Bhatt},
year={2023},
}

@incollection {Andersen-random,
    AUTHOR = {Andersen, Henning Haahr},
     TITLE = {A new proof of old character formulas},
 BOOKTITLE = {Invariant theory ({D}enton, {TX}, 1986)},
    SERIES = {Contemp. Math.},
    VOLUME = {88},
     PAGES = {193--207},
 PUBLISHER = {Amer. Math. Soc., Providence, RI},
      YEAR = {1989},
   MRCLASS = {20G10 (22E45)},
  MRNUMBER = {999992},
       DOI = {10.1090/conm/088/999992},
       URL = {https://doi.org/10.1090/conm/088/999992},
}

@article {mod-p-Harish-Chandra,
    AUTHOR = {Mirkovi\'{c}, Ivan and Rumynin, Dmitriy},
     TITLE = {Centers of reduced enveloping algebras},
   JOURNAL = {Math. Z.},
  FJOURNAL = {Mathematische Zeitschrift},
    VOLUME = {231},
      YEAR = {1999},
    NUMBER = {1},
     PAGES = {123--132},
      ISSN = {0025-5874},
   MRCLASS = {17B50 (17B35)},
  MRNUMBER = {1696760},
MRREVIEWER = {Kenneth A. Brown},
       DOI = {10.1007/PL00004719},
       URL = {https://doi.org/10.1007/PL00004719},
}

@misc{higher-coleman,
      title={Higher {C}oleman Theory}, 
      author={George Boxer and Vincent Pilloni},
      year={2021},
      eprint={2110.10251},
      archivePrefix={arXiv},
      primaryClass={math.NT},
      url={https://arxiv.org/abs/2110.10251}, 
}

@article{Hodge-paper,
author={Ortiz, Martin},
title={Theta operators on {H}odge type {S}himura varieties},
year={2026},
}

@misc{modular-O,
      title={On modular Soergel bimodules, Harish-Chandra bimodules, and category O}, 
      author={Ivan Losev},
      year={2023},
      eprint={2302.05782},
      archivePrefix={arXiv},
      primaryClass={math.RT},
      url={https://arxiv.org/abs/2302.05782}, 
}

@misc{Quan-O,
      title={Equivariant Koszul Duality, Modular Category $\mathcal{O}$, and Periodic Kazhdan--Lusztig Polynomials}, 
      author={Simon Riche and Quan Situ},
      year={2025},
      eprint={2511.18518},
      archivePrefix={arXiv},
      primaryClass={math.RT},
      url={https://arxiv.org/abs/2511.18518}, 
}

@misc{generic-vanishing-abelian,
      title={On the generic part of the cohomology of Shimura varieties of abelian type}, 
      author={Xiangqian Yang and Xinwen Zhu},
      year={2025},
      eprint={2505.04329},
      archivePrefix={arXiv},
      primaryClass={math.NT},
      url={https://arxiv.org/abs/2505.04329}, 
}

@article {smoothness-schubert,
    AUTHOR = {Billey, Sara and Postnikov, Alexander},
     TITLE = {Smoothness of {S}chubert varieties via patterns in root
              subsystems},
   JOURNAL = {Adv. in Appl. Math.},
  FJOURNAL = {Advances in Applied Mathematics},
    VOLUME = {34},
      YEAR = {2005},
    NUMBER = {3},
     PAGES = {447--466},
      ISSN = {0196-8858,1090-2074},
   MRCLASS = {14M15 (17B20)},
  MRNUMBER = {2123545},
MRREVIEWER = {Philip\ A.\ Foth},
       DOI = {10.1016/j.aam.2004.08.003},
       URL = {https://doi.org/10.1016/j.aam.2004.08.003},
}

@article {MR4258631,
    AUTHOR = {van der Geer, Gerard and Looijenga, Eduard},
     TITLE = {Lifting {C}hern classes by means of {E}kedahl-{O}ort strata},
   JOURNAL = {Tunis. J. Math.},
  FJOURNAL = {Tunisian Journal of Mathematics},
    VOLUME = {3},
      YEAR = {2021},
    NUMBER = {3},
     PAGES = {469--480},
      ISSN = {2576-7658,2576-7666},
   MRCLASS = {14G35 (11G18 14C25)},
  MRNUMBER = {4258631},
MRREVIEWER = {Amir\ D\v{z}ambi\'{c}},
       DOI = {10.2140/tunis.2021.3.469},
       URL = {https://doi.org/10.2140/tunis.2021.3.469},
}

@article {Deligne-Illusie,
    AUTHOR = {Deligne, Pierre and Illusie, Luc},
     TITLE = {Rel\`evements modulo {$p^2$} et d\'ecomposition du complexe de
              de {R}ham},
   JOURNAL = {Invent. Math.},
  FJOURNAL = {Inventiones Mathematicae},
    VOLUME = {89},
      YEAR = {1987},
    NUMBER = {2},
     PAGES = {247--270},
      ISSN = {0020-9910,1432-1297},
   MRCLASS = {14F40 (14C30)},
  MRNUMBER = {894379},
MRREVIEWER = {Thomas\ Zink},
       DOI = {10.1007/BF01389078},
       URL = {https://doi.org/10.1007/BF01389078},
}

@misc{le2025serreweightconjecturesmathrmgsp4,
      title={Serre weight conjectures for $\mathrm{GSp}_4$}, 
      author={Daniel Le and Bao V. Le Hung and Heejong Lee},
      year={2025},
      eprint={2510.04805},
      archivePrefix={arXiv},
      primaryClass={math.NT},
      url={https://arxiv.org/abs/2510.04805}, 
}

\end{document}